\documentclass[a4paper,12pt,amsart,frenchb]{article}

\usepackage{amsmath,amsbsy,amsfonts,amssymb}
\usepackage[french]{babel}
\newcommand{\ass}[2]{\vskip0.3cm\noindent
{\bf {#1}}. { \sl {#2}}\vskip0.3cm\noindent
}

\begin{document}

  \title{ Fronts d'onde des repr\'esentations temp\'er\'ees et  de r\'eduction unipotente pour $SO(2n+1)$ }
\author{J.-L. Waldspurger}
\date{28 septembre 2017}
 \maketitle
 
 {\bf Abstract} Let $G$ be a special orthogonal group $SO(2n+1)$ defined over a $p$-adic field $F$. Let $\pi$ be an admissible irreducible representation of $G(F)$ which is tempered and of unipotent reduction. We prove that $\pi$ has a wave front set. In some particular cases, for instance if $\pi$ is of the discrete series, we give a method to compute this wave front set. 
 \bigskip

{\bf Introduction}

Soit $F$ un corps local non archim\'edien et de caract\'eristique nulle et soit $n\geq1$ un entier. On suppose $p> 6n+4$, o\`u $p$ est la caract\'eristique r\'esiduelle de $F$.  Le groupe sp\'ecial orthogonal $SO(2n+1)$ a deux formes possibles d\'efinies sur $F$. Une forme d\'eploy\'ee que nous notons $G_{iso}$ et une forme non quasi-d\'eploy\'ee, qui est une forme int\'erieure du pr\'ec\'edent et que nous notons $G_{an}$. Soit $\sharp=iso$ ou $an$ et soit $\pi$ une repr\'esentation admissible irr\'eductible de $G_{\sharp}(F)$ dans un espace complexe $E$. Pour tout sous-groupe parahorique $K$ de $G_{\sharp}(F)$, notons $K^{u}$ son radical pro-$p$-unipotent et $E^{K^{u}}$ le sous-espace des \'el\'ements de $E$ fix\'es par $K^{u}$. De $\pi$ se d\'eduit une repr\'esentation de $K/K^{u}$ dans $E^{K^{u}}$. Le groupe $K/K^{u}$ s'identifie au groupe des points sur le corps r\'esiduel ${\mathbb F}_{q}$ de $F$ d'un groupe alg\'ebrique connexe d\'efini sur ${\mathbb F}_{q}$. Lusztig a d\'efini la notion de repr\'esentation unipotente d'un tel groupe. On dit que $\pi$ est de r\'eduction unipotente si et seulement s'il existe $K$ comme ci-dessus de sorte que $E^{K^{u}}$ soit non nul et que la repr\'esentation de $K/K^{u}$ dans $E^{K^{u}}$ soit unipotente. 
 
 Soit $\pi$ une repr\'esentation admissible irr\'eductible de $G_{\sharp}(F)$. Notons $\mathfrak{g}_{\sharp}$ l'alg\`ebre de Lie de $G_{\sharp}$. D'apr\`es Harish-Chandra,, dans un voisinage de l'origine dans $\mathfrak{g}_{\sharp}(F)$, le caract\`ere de $\pi$, descendu par l'exponentielle \`a $\mathfrak{g}_{\sharp}(F)$, est combinaison lin\'eaire de transform\'ees de Fourier d'int\'egrales orbitales nilpotentes. Fixons une cl\^oture alg\'ebrique $\bar{F}$ de $F$ et notons $\bar{{\cal N}}(\pi)$ l'ensemble des orbites nilpotentes $\boldsymbol{{\cal O}}$ dans $\mathfrak{g}_{\sharp}(\bar{F})$ qui v\'erifient la condition suivante: il existe une orbite nilpotente ${\cal O}$ dans $\mathfrak{g}_{\sharp}(F)$, qui est incluse dans $\boldsymbol{{\cal O}}$ et qui intervient avec un coefficient non nul dans le d\'eveloppement du caract\`ere de $\pi$. On dit que $\pi$ admet un front d'onde si $\bar{{\cal N}}(\pi)$ admet un plus grand \'el\'ement (pour l'ordre usuel sur les orbites nilpotentes). Si c'est le cas, on appelle ce plus grand \'el\'ement le front d'onde de $\pi$. Le th\'eor\`eme principal de l'article est le suivant.
 
 \ass{Th\'eor\`eme}{Soit $\sharp=iso$ ou $an$. Alors toute repr\'esentation admissible irr\'eductible de $G_{\sharp}(F)$, qui est temp\'er\'ee et de r\'eduction unipotente, admet un front d'onde. }

Pour tout entier $N\in {\mathbb N}$, notons ${\cal P}^{symp}(2N)$ l'ensemble des partitions symplectiques de $2N$ (une partition est dite symplectique si tout entier impair y intervient avec une multiplicit\'e paire). Pour une telle partition $\lambda$,  notons $Jord^{bp}(\lambda)$ l'ensemble (sans multiplicit\'es) des entiers pairs strictement positifs qui interviennent dans $\lambda$. Notons $\boldsymbol{{\cal P}^{symp}}(2N)$ l'ensemble des couples $(\lambda,\epsilon)$ o\`u $\lambda\in {\cal P}^{symp}(2N)$ et $\epsilon\in \{\pm 1\}^{Jord^{bp}(\lambda)}$. Notons $\mathfrak{Irr}_{quad}(2n)$ l'ensemble des quadruplets $(\lambda^+,\epsilon^+,\lambda^-,\epsilon^-)$ pour lesquels il existe deux entiers $n^+$ et $n^-$ de sorte que $n^++n^-=n$, $(\lambda^+,\epsilon^+)\in \boldsymbol{{\cal P}^{symp}}(2n^+)$ et $(\lambda^-,\epsilon^-)\in \boldsymbol{{\cal P}^{symp}}(2n^-)$. A un tel quadruplet 
$(\lambda^+,\epsilon^+,\lambda^-,\epsilon^-)$, on peut associer un indice $\sharp=iso$ ou $an$ et une repr\'esentation admissible irr\'eductible $\pi(\lambda^+,\epsilon^+,\lambda^-,\epsilon^-)$ de $G_{\sharp}(F)$, qui est temp\'er\'ee et de r\'eduction unipotente. L'indice $\sharp$ est d\'etermin\'e par une formule simple rappel\'ee en 1.5.  Indiquons bri\`evement quel est le param\`etre de Langlands de cette repr\'esentation. Notons $W_{F}$ le groupe de Weil de $F$ et $W_{DF}=W_{F}\times SL(2,{\mathbb C})$ le groupe de Weil-Deligne. Un  param\`etre  de Langlands est un couple $(\rho,\chi)$, o\`u $\rho$ est un homomorphisme de $W_{DF}$ dans $Sp(2n;{\mathbb C})$ et $\chi$ est un caract\`ere du groupe des composantes connexes du commutant dans $Sp(2n;{\mathbb C})$ de l'image de $\rho$. Dans le cas d'une repr\'esentation  $\pi(\lambda^+,\epsilon^+,\lambda^-,\epsilon^-)$, la restriction de $\rho$ \`a $W_{F}$ est la somme directe de $2n^+$ fois le caract\`ere trivial et de $2n^-$ fois l'unique caract\`ere non ramifi\'e d'ordre $2$. Le commutant de l'image de cette restriction est un groupe $Sp(2n^+;{\mathbb C})\times Sp(2n^-;{\mathbb C})$. Les classes de conjugaison d'\'el\'ements unipotents dans ce groupe sont param\'etr\'ees par ${\cal P}^{symp}(2n^+)\times {\cal P}^{symp}(2n^-)$. La restriction de $\rho$ \`a $SL(2;{\mathbb C})$ prend ses valeurs dans ce groupe et l'image d'un unipotent non trivial de $SL(2;{\mathbb C})$ est param\'etr\'e par  $(\lambda^+,\lambda^-)$. On voit que le groupe des composantes connexes du commutant dans $Sp(2n;{\mathbb C})$ de l'image de $\rho$ est isomorphe \`a $({\mathbb Z}/2{\mathbb Z})^{Jord^{bp}(\lambda^+)}\times ({\mathbb Z}/2{\mathbb Z})^{Jord^{bp}(\lambda^-)}$. Le couple $(\epsilon^+,\epsilon^-)$ s'identifie \`a un caract\`ere de ce groupe, qui n'est autre que le caract\`ere $\chi$ du couple $(\rho,\chi)$.

On note  $\mathfrak{Irr}_{quad}^{bp}(2n)$ le sous ensemble des $(\lambda^+,\epsilon^+,\lambda^-,\epsilon^-)\in \mathfrak{Irr}_{quad}(2n)$ tels que tous les termes de $\lambda^+$ et $\lambda^-$ soient pairs. 
On a prouv\'e en \cite{W4} 3.4 que, pour d\'emontrer le th\'eor\`eme, il suffisait de prouver que, pour tout $(\lambda^+,\epsilon^+,\lambda^-,\epsilon^-)\in \mathfrak{Irr}_{quad}^{bp}(2n)$, la  repr\'esentation  $\pi(\lambda^+,\epsilon^+,\lambda^-,\epsilon^-)$ admettait un front d'onde (cela r\'esulte d'un argument trivial d'induction). 

Pour une  repr\'esentation $\pi(\lambda^+,\epsilon^+,\lambda^-,\epsilon^-)$,  o\`u $(\lambda^+,\epsilon^+,\lambda^-,\epsilon^-)\in \mathfrak{Irr}_{quad}^{bp}(2n)$, on a un r\'esultat un  peu plus pr\'ecis. Dans \cite{W5}, on a  \'etudi\'e une certaine repr\'esentation d'un groupe de Weyl d\'efinie par Lusztig. En supposant, comme c'est ici le cas, que tous les termes de $\lambda^+$ sont pairs, on a associ\'e \`a $(\lambda^+,\epsilon^+) \in \boldsymbol{{\cal P}^{symp}}(2n^+)$  un  autre couple $(\lambda^{+,min},\epsilon^{+,min})\in \boldsymbol{{\cal P}^{symp}}(2n^+)$ (voir ci-dessous). De m\^eme, \`a $(\lambda^-,\epsilon^-)\in \boldsymbol{{\cal P}^{symp}}(2n^+)$, on associe  un  autre couple $(\lambda^{-,min},\epsilon^{-,min})\in \boldsymbol{{\cal P}^{symp}}(2n^-)$. Introduisons la r\'eunion usuelle de $\lambda^{+,min}$ et $\lambda^{-,min}$, que l'on note $\lambda^{+,min}\cup \lambda^{-,min}$. C'est une partition symplectique de $2n$. Notons ${\cal P}^{orth}(2n+1)$ l'ensemble des partitions orthogonales de $2n+1$ (une partition est orthogonale si et seulement si tout entier pair strictement positif y intervient avec multiplicit\'e paire). On sait bien que l'ensemble ${\cal P}^{orth}(2n+1)$ param\`etre les orbites nilpotentes dans $\mathfrak{g}_{\sharp}(\bar{F})$. Un front d'onde  est donc param\'etr\'e par un \'el\'ement de cet ensemble. D'autre part, \`a la suite de Spaltenstein, on d\'efinit une "dualit\'e" $d:{\cal P}^{symp}(2n)\to {\cal P}^{orth}(2n+1)$, cf. 2.6 (elle n'est ni injective, ni surjective, son image est le sous-ensemble des partitions "sp\'eciales" dans ${\cal P}^{orth}(2n+1)$). 

\ass{Th\'eor\`eme}{ Soit $(\lambda^+,\epsilon^+,\lambda^-,\epsilon^-)\in \mathfrak{Irr}^{bp}_{quad}(2n)$. Alors la repr\'esentation $\pi(\lambda^+,\epsilon^+,\lambda^-,\epsilon^-)$  admet un front d'onde. Celui-ci est param\'etr\'e par la partition $d(\lambda^{+,min}\cup \lambda^{-,min})$.}

La preuve de ce th\'eor\`eme reprend celle de \cite{W5}. Posons $\pi=\pi(\lambda^+,\epsilon^+,\lambda^-,\epsilon^-)$. L'existence d'un front d'onde pour $\pi$ se lit sur le caract\`ere de cette repr\'esentation. Celui-ci se calcule en fonction des repr\'esentations des diff\'erents  groupes finis $K/K^{u}$ dans $E^{K^{u}}$, avec les notations du premier paragraphe ci-dessus (en v\'erit\'e, le groupe fini est $K^{\dag}/K^{u}$, o\`u $K^{\dag}$ est le normalisateur de $K$ dans $G_{\sharp}(F)$). La construction de la repr\'esentation $\pi$  (qui est due \`a Lusztig) permet d'expliciter ces repr\'esentations de groupes finis. On les d\'ecrit \`a l'aide de repr\'esentations de groupes de Weyl $W_{m}$ de type $B_{m}$ ou $C_{m}$. Une vieille combinatoire tir\'ee de \cite{W1} permet alors de traduire l'existence d'un front d'onde et son calcul en un probl\`eme concernant exclusivement des repr\'esentations de tels groupes $W_{m}$, cf. 1.4. Les objets cruciaux qui interviennent ici sont les repr\'esentations $\boldsymbol{\rho}_{\lambda^+,\epsilon^+}$ et $\boldsymbol{\rho}_{\lambda^-,\epsilon^-}$ d\'efinies par Lusztig (ce ne sont pas ses notations) auxquelles on a fait allusion ci-dessus. Elles ne sont pas irr\'eductibles en g\'en\'eral et on conna\^{\i}t peu de choses sur leur d\'ecomposition en repr\'esentations irr\'eductibles. On sait toutefois que, disons dans la d\'ecomposition de $\boldsymbol{\rho}_{\lambda^+,\epsilon^+}$, il y a un \'el\'ement "minimal"  qui est la repr\'esentation $\rho_{\lambda^+,\epsilon^+}$ associ\'ee \`a $(\lambda^+,\epsilon^+)$ par la correspondance de Springer g\'en\'eralis\'ee. Dans \cite{W4}, cela nous a suffi pour traiter non pas la repr\'esentation $\pi$, mais son image par l'involution d'Aubert-Zelevinsky. Le point nouveau est le r\'esultat de \cite{W5} qui affirme (sous l'hypoth\`ese que tous les termes de $\lambda^+$ sont pairs) que la d\'ecomposition de $\boldsymbol{\rho}_{\lambda^+,\epsilon^+}$ admet aussi un \'el\'ement "maximal" pour un ordre convenable (cf. 4.1 pour un \'enonc\'e pr\'ecis). C'est $\rho_{\lambda^{+,min},\epsilon^{+,min}}\otimes sgn$, o\`u $sgn$ est le caract\`ere signature. Cette propri\'et\'e nous permet de conclure.

Les paragraphes 1 \`a 3 sont surtout consacr\'es \`a des rappels de r\'esultats ant\'erieurs. On a am\'elior\'e certains d'entre eux quand c'\'etait n\'ecessaire. Le th\'eor\`eme ci-dessus est d\'emontr\'e au paragraphe 4. 
 Dans le paragraphe 5, nous indiquons comment se calculent les partitions $\lambda^{+,min}$ et $\lambda^{-,min}$ (en fait leurs transpos\'ees) et nous donnons quelques exemples de fronts d'onde.

Je remercie A.-M. Aubert de m'avoir indiqu\'e une r\'ef\'erence utile. 
 \bigskip
 
 \section{Rappel pas tr\`es bref des r\'esultats de \cite{W4}}
 
 \bigskip
 
 \subsection{Partitions, notations}

 Soit $\lambda=(\lambda_{1},...,\lambda_{r})$ une suite finie de nombres r\'eels. Notons $l(\lambda)$ le plus grand entier $j\in \{1,...,r\}$ tel que $\lambda_{j}\not=0$. On identifie deux suites $\lambda$ et $\lambda'$ si $l(\lambda)=l(\lambda')$ et $\lambda_{j}=\lambda'_{j}$ pour tout $j\leq l(\lambda)$. 
Soit $\lambda$ une telle suite 
 et soit $k\in {\mathbb N}$. Quitte \`a adjoindre \`a $\lambda$ des termes nuls, on peut \'ecrire $\lambda=(\lambda_{1},...,\lambda_{r})$ avec $r\geq k$.  On pose $S_{k}(\lambda)=\lambda_{1}+...+\lambda_{k}$.  Evidemment, $S_{k}(\lambda)$ ne d\'epend plus de $k$ d\`es que  $k\geq l(\lambda)$. On pose $S(\lambda)=S_{l(\lambda)}(\lambda)$. On d\'efinit la somme $\lambda+\lambda'$ de deux suites $\lambda$ et $\lambda'$: $(\lambda+\lambda')_{j}=\lambda_{j}+\lambda'_{j}$ pour tout $j\geq1$.
 
Une partition est une suite finie d\'ecroissante d'entiers positifs ou nuls. On identifie  comme ci-dessus deux partitions qui ne diff\`erent que par des termes nuls.  Pour une partition $\lambda=(\lambda_{1},...,\lambda_{r})$ et pour un entier $i\geq1$, on note $mult_{\lambda}(i)$ le nombre d'indices $j$ tels que $\lambda_{j}=i$. On note $Jord(\lambda)$ l'ensemble des $i\geq 1$ tels que $mult_{\lambda}(i)>0$.   Pour tout  $N\in {\mathbb N}$, on note ${\cal P}(N)$ l'ensemble des partitions  $\lambda$ telles que $S(\lambda)=N$ et on note ${\cal P}_{2}(N)$ l'ensemble des couples $(\alpha,\beta)$ de partitions telles que $S(\alpha)+S(\beta)=N$.  On ordonne les \'el\'ements de ${\cal P}(N)$ de la fa\c{c}on usuelle: $\lambda\leq \lambda'$ si et seulement si $S_{k}(\lambda)\leq S_{k}(\lambda')$ pour tout $k\in {\mathbb N}$.   On d\'efinit la r\'eunion $\lambda\cup\lambda'$ de deux partitions $\lambda$ et $\lambda'$: pour tout entier $i\geq1$, $mult_{\lambda\cup \lambda'}(i)=mult_{\lambda}(i)+mult_{\lambda'}(i)$.  

Soit $\lambda$ une partition.  Pour tout $i\in {\mathbb N}$, on note $J(i)$ l'ensemble des $j\geq1$ tels que $\lambda_{j}=i$. Si $i=0$, on consid\`ere que $J(0)$ est l'intervalle infini $\{l(\lambda)+1,...\}$. Pour $i\in Jord(\lambda)$, $J(i)$ est non vide. On note $j_{min}(i)$, resp. $j_{max}(i)$, le plus petit, resp. grand, \'el\'ement de $J(i)$. On pose $j_{min}(0)=l(\lambda)+1$.

On note $W_{N}$ le groupe de Weyl d'un syst\`eme de racines de type $B_{N}$ ou $C_{N}$, avec la convention $W_{0}=\{1\}$. On note $sgn$ le caract\`ere signature usuel de $W_{N}$ et $sgn_{CD}$ le caract\`ere dont le noyau est le sous-groupe $W_{N}^D$ d'un syst\`eme de racines de type $D_{N}$. Les repr\'esentations irr\'eductibles de $W_{N}$ sont param\'etr\'ees par ${\cal P}_{2}(N)$. Pour $(\alpha,\beta)\in {\cal P}_{2}(N)$, on note $\rho(\alpha,\beta)$ la repr\'esentation param\'etr\'ee par $(\alpha,\beta)$. Les repr\'esentations irr\'eductibles de $W_{N}^D$ sont presque param\'etr\'ees par le quotient de ${\cal P}_{2}(N)$ par la relation d'\'equivalence $(\alpha,\beta)\equiv (\beta,\alpha)$. "Presque" parce qu'un couple de la forme $(\alpha,\alpha)$ param\`etre deux repr\'esentations irr\'eductibles. 

 Pour tout ensemble $E$, on note ${\mathbb C}[E]$ le ${\mathbb C}$-espace vectoriel de base $E$. Pour tout groupe fini $W$, on note $\hat{W}$ l'ensemble des classes d'\'equivalence de repr\'esentations irr\'eductibles de $W$. En identifiant une repr\'esentation \`a son caract\`ere, ${\mathbb C}[\hat{W}]$ est aussi l'espace des fonctions de $W$ dans ${\mathbb C}$ qui sont invariantes par conjugaison. 

\bigskip

\subsection{L'espace ${\cal R}$}

{\bf On fixe pour tout l'article un entier $n\geq1$}.  On note $\Gamma$ l'ensemble des quadruplets $\gamma=(r',r'',N^+,N^-)$ tels que
$$r'\in {\mathbb N},\,\,r''\in {\mathbb Z},\,\,N^+\in {\mathbb N},\,\,N^-\in {\mathbb N},\,\, r^{'2}+r'+N^++r^{''2}+N^-=n.$$
Pour un tel $\gamma$, on pose ${\cal R}(\gamma)={\mathbb C}[\hat{W}_{N^+}]\otimes {\mathbb C}[\hat{W}_{N^-}]$. On pose
$${\cal R}=\oplus_{\gamma\in \Gamma}{\cal R}(\gamma).$$

On d\'efinit un endomorphisme $\varphi\mapsto sgn\otimes \varphi$ de ${\cal R}$ de la fa\c{c}on suivante. Il respecte chaque sous-espace ${\cal R}(\gamma)$. Pour $\gamma$ comme ci-dessus, pour $\rho^+\in \hat{W}_{N^+}$ et $\rho^-\in \hat{W}_{N^-}$, on pose $sgn\otimes (\rho^+\otimes \rho^-)=(\rho^+\otimes sgn, \rho^-\otimes sgn)$.

On a d\'efini en \cite{W2} 1.10 un endomorphisme $\rho\iota$. Puisqu'il est essentiel \`a nos constructions, rappelons sa d\'efinition.  Soit $\gamma=(r',r'',N^+,N^-)\in \Gamma$ et $\varphi\in {\cal R}(\gamma) $. Posons $N=N^++N^-$. L'\'el\'ement $\rho\iota(\varphi)$ appartient \`a 
$$\oplus_{N_{1},N_{2}\in {\mathbb N},N_{1}+N_{2}=N}{\cal R}(r',(-1)^{r'}r'',N_{1},N_{2}).$$
 Soit $\delta=(r',(-1)^{r'},N_{1},N_{2})\in \Gamma$. D\'ecrivons  la composante $\rho\iota(\varphi)_{\delta}$ de $\rho\iota(\varphi)$ dans ${\cal R}(\delta)$.
    
  On d\'efinit un quadruplet d'entiers ${\bf a}=(a_{1}^+,a_{1}^-,a_{2}^+,a_{2}^-)$ par les formules suivantes:

${\bf a}=(0,0,0,1)$  si $0< r''\leq r'$ ou si  $r''=0$ et $r'$ est pair;

${\bf a}=(0,0,1,0)$ si  $-r'\leq r''<0$ ou si  $r''=0$ et $r'$ est impair;

${\bf a}=(0,1,0,0)$ si $r'<r''$;

${\bf a}=(1,0,0,0)$ si   $r''<-r'$.
 
  Notons ${\cal N}$ l'ensemble des quadruplets ${\bf N}=(N^+_{1},N^-_{1},N^+_{2},N^-_{2})$ d'entiers positifs ou nuls tels que
  $$N^+= N^+_{1}+N^+_{2},\,\,N^-=N^-_{1}+N^-_{2},\,\,N_{1}=N^+_{1}+N^-_{1},\,\,N_{2}=N^+_{2}+N^-_{2}.$$
 Pour un tel quadruplet, 
posons $W_{{\bf N}}=W_{N^+_{1}}\times W_{N^-_{1}}\times W_{N^+_{2}}\times W_{N^-_{2}}$. Ce groupe se plonge de fa\c{c}on \'evidente dans $W_{N_{1}}\times W_{N_{2}}$, resp. $W_{N^+}\times W_{N^-}$, et ces plongements sont bien d\'efinis \`a conjugaison pr\`es. On a donc des foncteurs de restriction $res^{W_{N^+}\times W_{N^-}}_{W_{{\bf N}}}$ et d'induction $ind^{W_{N_{1}}\times W_{N_{2}}}_{W_{{\bf N}}}$. On note $sgn_{CD}^{{\bf a}}$ le caract\`ere de $W_{{\bf N}}$ qui est le produit tensoriel des caract\`eres $sgn_{CD}^{a^+_{1}}$, $sgn_{CD}^{a^-_{1}}$, $sgn_{CD}^{a^+_{2}}$, $sgn_{CD}^{a^-_{2}}$ sur chacun des facteurs de $W_{{\bf N}}$. Alors 
$$\rho\iota(\varphi)_{\delta}=\sum_{{\bf N}\in {\cal N}}ind^{W_{N_{1}}\times W_{N_{2}}}_{W_{{\bf N}}}\left(sgn_{CD}^{{\bf a}}\otimes res^{W_{N^+}\times W_{N^-}}_{W_{{\bf N}}}( \varphi )\right).$$

\bigskip

\subsection{Correspondance de Springer g\'en\'eralis\'ee}
Soit $N\in {\mathbb N}$. On a d\'efini l'ensemble $\boldsymbol{{\cal P}^{symp}}(2N)$ dans l'introduction. La correspondance de Springer g\'en\'eralis\'ee dans le cas symplectique  est une bijection de $\boldsymbol{{\cal P}^{symp}}(2N)$ sur l'ensemble des couples $(k,\rho)$ o\`u

$k\in {\mathbb N}$ et $k(k+1)\leq 2N$;

$\rho\in \hat{W}_{N-k(k+1)/2}$.

Pour $(\lambda,\epsilon)\in \boldsymbol{{\cal P}^{symp}}(2N)$, on note $(k_{\lambda,\epsilon},\rho_{\lambda,\epsilon})$ le couple qui lui correspond et on pose $N_{\lambda, \epsilon}=N-k_{\lambda,\epsilon}(k_{\lambda,\epsilon}+1)/2$.  Rappelons comment on calcule $k_{\lambda,\epsilon}$. On note $i_{1}>...>i_{t}$ les entiers $i\in Jord^{bp}(\lambda)$ tels que $mult_{\lambda}(i)$ soit impair. On pose
$$M=\vert \{h=1,...,t; h\text{\,\,est\,\,pair\,\,et\,\,}\epsilon_{i_{h}}=-1\}\vert -\vert \{h=1,...,t; h\text{\,\,est\,\,impair\,\,et\,\,}\epsilon_{i_{h}}=-1\}\vert.$$
Alors, d'apr\`es \cite{W1} XI.3, on a 

(1) $k_{\lambda,\epsilon}=2M$ si $M\geq0$, $k_{\lambda,\epsilon}=-2M-1$ si $M<0$.

 On d\'efinit une autre repr\'esentation $\boldsymbol{\rho}_{\lambda,\epsilon}$ du m\^eme groupe $W_{N_{\lambda,\epsilon}}$, cf. \cite{W2} 5.1. En gros, $\rho_{\lambda,\epsilon}$ est l'action de $W_{N_{\lambda,\epsilon}}$ sur un sous-espace d\'etermin\'e par $\epsilon$ de l'espace de cohomologie de plus haut degr\'e d'une certaine vari\'et\'e alg\'ebrique, tandis que $\boldsymbol{\rho}_{\lambda,\epsilon}$ est l'action de $W_{N_{\lambda,\epsilon}}$ sur un sous-espace analogue de la somme de tous les espaces de cohomologie de cette vari\'et\'e.

Soit $(\lambda^+,\epsilon^+,\lambda^-,\epsilon^-)\in \mathfrak{Irr}_{quad}(2n)$. Pour $\zeta=\pm$, posons $2n^{\zeta}=S(\lambda^{\zeta})$,  $k^{\zeta}=k_{\lambda^{\zeta},\epsilon^{\zeta}}$, $N^{\zeta}=n^{\zeta}-k^{\zeta}(k^{\zeta}+1)/2$. On d\'efinit des entiers $r'\in{\mathbb N}$, $r''\in {\mathbb Z}$ par les formules suivantes

si $k^+\equiv k^-\,\,mod\,\,2{\mathbb Z}$, $r'=\frac{k^++k^-}{2}$, $r''=\frac{k^+-k^-}{2}$;

si  $k^+\not\equiv k^-\,\,mod\,\,2{\mathbb Z}$ et $k^+>k^-$, $r'=\frac{k^+-k^--1}{2}$, $r''=\frac{k^++k^-+1}{2}$;

si  $k^+\not\equiv k^-\,\,mod\,\,2{\mathbb Z}$ et $k^+<k^-$, $r'=\frac{k^--k^+-1}{2}$, $r''=-\frac{k^++k^-+1}{2}$.

Le quadruplet $\gamma=(r',r'',N^+,N^-)$ appartient \`a $\Gamma$. Puisque ${\cal R}(\gamma)={\mathbb C}[\hat{W}_{N^+}]\otimes {\mathbb C}[\hat{W}_{N^-}]$, on peut identifier $\boldsymbol{\rho}_{\lambda^+,\epsilon^+}\otimes \boldsymbol{\rho}_{\lambda^-,\epsilon^-}$ \`a un \'el\'ement de ${\cal R}(\gamma)$, a fortiori \`a un \'el\'ement de ${\cal R}$. Dans la suite $\boldsymbol{\rho}_{\lambda^+,\epsilon^+}\otimes \boldsymbol{\rho}_{\lambda^-,\epsilon^-}$  d\'esignera cet \'el\'ement. 

Pour $M\in {\mathbb N}$, on note ${\cal P}^{orth}(M)$ l'ensemble des partitions orthogonales de $M$. Pour une telle partition $\lambda$, on note $Jord^{bp}(\lambda)$ l'ensemble des entiers impairs $i\geq1$  tels que $mult_{\lambda}(i)>0$. On note $\boldsymbol{{\cal P}^{orth}}(M)$ l'ensemble des couples $(\lambda,\epsilon)$ o\`u $\lambda\in {\cal P}^{orth}(M)$ et $\epsilon\in \{\pm 1\}^{Jord^{bp}(\lambda)}/\{\pm 1\}$, le groupe $\{\pm 1\}$ s'envoyant diagonalement dans $\{\pm 1\}^{Jord^{bp}(\lambda)}$.

Soit $N\in {\mathbb N}$. La correspondance de Springer g\'en\'eralis\'ee dans le cas orthogonal impair est une bijection de $\boldsymbol{{\cal P}^{orth}}(2N+1)$ sur l'ensemble des couples $(k,\rho)$ tels que

$k\in {\mathbb N}$, $k$ est impair et $k^2\leq 2n+1$;

$\rho\in \hat{W}_{N-(k^2-1)/2}$.

On note $(k_{\lambda,\epsilon},\rho_{\lambda,\epsilon})$ le couple associ\'e \`a $(\lambda,\epsilon)\in \boldsymbol{{\cal P}^{orth}}(2N+1)$. On note $\boldsymbol{{\cal P}^{orth}}(2N+1)_{k=1}$ le sous-ensemble des $(\lambda,\epsilon)\in \boldsymbol{{\cal P}^{orth}}(2N+1)$ tels que $k_{\lambda,\epsilon}=1$.

La correspondance de Springer g\'en\'eralis\'ee dans le cas orthogonal pair est une bijection entre $\boldsymbol{{\cal P}^{orth}}(2N)$ et l'ensemble des couples $(k,\rho)$ tels que
 
 $k\in {\mathbb N}$, $k$ est pair et $k^2\leq 2n$;
 
 si $k>0$, $\rho\in \hat{W}_{N-k^2/2}$;
 
 si $k=0$, $\rho$ est une classe d'\'equivalence dans $\hat{W}_{N-k^2/2}$, deux repr\'esentations irr\'eductibles $\rho'$ et $\rho''$ \'etant  ici \'equivalentes si et seulement si $\rho'=\rho''$ ou $\rho'=\rho''\otimes sgn_{CD}$.
 
 On note $(k_{\lambda,\epsilon},\rho_{\lambda,\epsilon})$ le couple associ\'e \`a $(\lambda,\epsilon)\in \boldsymbol{{\cal P}^{orth}}(2N)$. On note $\boldsymbol{{\cal P}^{orth}}(2N)_{k=0}$ le sous-ensemble des $(\lambda,\epsilon)\in \boldsymbol{{\cal P}^{orth}}(2N)$ tels que $k_{\lambda,\epsilon}=0$. Quand $k_{\lambda,\epsilon}=0$, $\rho_{\lambda,\epsilon}$ n'est qu'une classe d'\'equivalence comme on vient de le dire. Autrement dit, $\rho_{\lambda,\epsilon}$ est param\'etr\'ee par un couple $(\alpha,\beta)\in {\cal P}_{2}(N)$ \`a l'ordre pr\`es. Si $\alpha=\beta$, on pose $\rho^+_{\lambda,\epsilon}=\rho^{-}_{\lambda,\epsilon}=\rho(\alpha,\beta)$. Si $\alpha\not=\beta$, on choisit $\alpha$ et $\beta$ de sorte que $\alpha> \beta$ pour l'ordre lexicographique. On pose $\rho^+_{\lambda,\epsilon}=\rho(\alpha,\beta)$ et $\rho^-_{\lambda,\epsilon}=\rho(\beta,\alpha)$.
 
 \bigskip 
 \subsection{Caract\'erisation du front d'onde}
 
 On a introduit les groupes $G_{iso}$ et $G_{an}$. Pour $\sharp=iso$ ou $an$, on note $Irr_{tunip,\sharp}$ l'ensemble des classes d'isomorphismes de repr\'esentations admissibles irr\'eductibles de $G_{\sharp}(F)$ qui sont temp\'er\'ees et de r\'eduction unipotente. On note $Irr_{tunip}$ la r\'eunion disjointe de $Irr_{tunip,iso}$ et $Irr_{tunip,an}$. On a d\'efini en \cite{W3} 1.5 un espace ${\cal R}^{par}$ et une application lin\'eaire $Rep:{\mathbb C}[Irr_{tunip}]\to {\cal R}^{par}$. A la suite de Lusztig, on a d\'efini en \cite{MW} 3.16  deux isomorphismes $Rep:{\cal R}\to {\cal R}^{par}$ et $k:{\cal R}\to {\cal R}^{par}$. On note ${\cal F}$ l'automorphisme  de ${\cal R}$ tel que $Rep\circ {\cal F}=k$. C'est une involution sur le calcul de laquelle nous reviendrons en 2.5. Pour $\pi\in Irr_{tunip}$, on note $\kappa_{\pi}$ l'\'el\'ement de ${\cal R}$ tel que $k(\kappa_{\pi})=Res(\pi)$. Soient $n_{1},n_{2}\in {\mathbb N}$ et $\rho_{1}\in \hat{W}_{n_{1}}$, $\rho_{2}\in \hat{W}_{n_{2}}$. Le quadruplet $\gamma=(0,0,n_{1},n_{2})$ appartient \`a $\Gamma$ et on a ${\cal R}(\gamma)={\mathbb C}[\hat{W}_{n_{1}}]\otimes {\mathbb C}[\hat{W}_{n_{2}}]$. Notons $\kappa_{\pi}(\gamma)$ la composante de $\kappa_{\pi}$ dans ${\cal R}(\gamma)$. C'est une combinaison lin\'eaire de repr\'esentations irr\'eductibles avec des coefficients complexes. On note $m_{\pi}(\rho_{1},\rho_{2})$ le coefficient de $\rho_{1}\otimes \rho_{2}$ dans cette combinaison lin\'eaire. 
 
 On pose $sgn_{iso}=1$, $sgn_{an}=-1$. Soit $\sharp=iso$ ou $an$, soit $\pi\in Irr_{tunip,\sharp}$, soient $n_{1},n_{2}\in {\mathbb N}$ tels que $n_{1}+n_{2}=n$ et soient $(\mu_{1},\eta_{1})\in \boldsymbol{{\cal P}^{orth}}(2n_{1}+1)_{k=1}$ et $(\mu_{2},\eta_{2})\in \boldsymbol{{\cal P}^{orth}}(2n_{2})_{k=0}$. Comme cas particulier de la d\'efinition ci-dessus, pour $\zeta=\pm$, on d\'efinit les "multiplicit\'es" $m_{\pi}(\rho_{\mu_{1},\eta_{1}}\otimes sgn,\rho_{\mu_{2},\eta_{2}}^{\zeta}\otimes sgn)$. On pose
 $$M_{\pi}(\mu_{1},\eta_{1};\mu_{2},\eta_{2})=m_{\pi}(\rho_{\mu_{1},\eta_{1}}\otimes sgn,\rho_{\mu_{2},\eta_{2}}^{+}\otimes sgn)+sgn_{\sharp}m_{\pi}(\rho_{\mu_{1},\eta_{1}}\otimes sgn,\rho_{\mu_{2},\eta_{2}}^{-}\otimes sgn).$$
 
 \ass{Proposition}{Soient $\sharp=iso$ ou $an$, $\pi\in Irr_{tunip,\sharp}$ et $\mu\in {\cal P}^{orth}(2n+1)$. Alors $\pi$ admet un front d'onde param\'etr\'e par $\mu$ si et seulement si les deux conditions suivantes sont v\'erifi\'ees.
 
 (i) Soient $n_{1},n_{2}\in {\mathbb N}$ tels que $n_{1}+n_{2}=n$ et $n_{2}\geq1$ si $\sharp=an$. Soient $(\mu_{1},\eta_{1})\in \boldsymbol{{\cal P}^{orth}}(2n_{1}+1)_{k=1}$ et $(\mu_{2},\eta_{2})\in \boldsymbol{{\cal P}^{orth}}(2n_{2})_{k=0}$. Supposons $M_{\pi}(\mu_{1},\eta_{1};\mu_{2},\eta_{2})\not=0$. Alors $\mu_{1}\cup \mu_{2}\leq \mu$.
 
 (ii) Il existe $n_{1},n_{2}\in {\mathbb N}$ tels que $n_{1}+n_{2}=n$ et $n_{2}\geq1$ si $\sharp=an$ et il existe  $(\mu_{1},\eta_{1})\in \boldsymbol{{\cal P}^{orth}}(2n_{1}+1)_{k=1}$ et $(\mu_{2},\eta_{2})\in \boldsymbol{{\cal P}^{orth}}(2n_{2})_{k=0}$ tels que  $M_{\pi}(\mu_{1},\eta_{1};\mu_{2},\eta_{2})\not=0$ et $\mu_{1}\cup \mu_{2}= \mu$.}
 
 Cf. \cite{W4} 3.7. Les notations de cette r\'ef\'erence \'etaient l\'eg\`erement diff\'erente, les multiplicit\'es \'etaient dans certains cas divis\'ees par $2$ mais cela ne change \'evidemment pas l'\'enonc\'e. D'autre part, dans \cite{W4}, la repr\'esentation $\pi$ \'etait d'une forme particuli\`ere, mais cela n'\'etait utilis\'e que pour d\'ecrire explicitement la fonction $\kappa_{\pi}$ dans le paragraphe 3.8 de loc. cit., cela n'intervient pas \`a ce point. 
 
 \bigskip
 
 \subsection{Les repr\'esentations $\pi(\lambda^+,\epsilon^+,\lambda^-,\epsilon^-)$}

 Soit $(\lambda^+,\epsilon^+,\lambda^-,\epsilon^-)\in \mathfrak{Irr}_{quad}(2n)$. En utilisant une construction de Lusztig, on a d\'efini en \cite{W3} 1.3 la repr\'esentation $\pi(\lambda^+,\epsilon^+,\lambda^-,\epsilon^-)$. Sa param\'etrisation de Langlands a \'et\'e rappel\'ee rapidement dans l'introduction. C'est une repr\'esentation  admissible, irr\'eductible et temp\'er\'ee de $G_{\sharp}(F)$, o\`u l'indice $\sharp$ est d\'etermin\'e par la formule
 $$(1)\qquad sgn_{\sharp}=(\prod_{i\in Jord^{bp}(\lambda^+)}\epsilon^+(i)^{mult_{\lambda^+}(i)})(\prod_{i\in Jord^{bp}(\lambda^-)}\epsilon^-(i)^{mult_{\lambda^-}(i)}),$$
 cf. 1.3 pour la d\'efinition de $sgn_{\sharp}$. Notons $D$ l'involution de Aubert-Zelevinski. Elle permute les repr\'esentations admissibles irr\'eductibles de $G_{\sharp}(F)$. On a l'\'egalit\'e
 
 (2) $Res\circ D(\pi(\lambda^+,\epsilon^+,\lambda^-,\epsilon^-))=Rep\circ\rho\iota(\boldsymbol{\rho}_{\lambda^+,\epsilon^+}\otimes \boldsymbol{\rho}_{\lambda^-,\epsilon^-}),$
 
 \noindent cf. \cite{W3} proposition 1.11.  
 
 L'espace ${\cal R}^{par}$ est somme directe finie d'espaces vectoriels ayant pour base les classes d'\'equivalence de repr\'esentations irr\'eductibles et unipotentes de groupes finis de la forme $SO(2n'+1;{\mathbb F}_{q})\times O(2n'';{\mathbb F}_{q})$, avec des notations compr\'ehensibles. Chacun de ces espaces est muni d'involutions du m\^eme type que $D$. L'espace 
 ${\cal R}^{par}$ est ainsi muni d'une involution $D^{par}$ et on a prouv\'e en \cite{W3} 1.7 l'\'egalit\'e $Res\circ D=D^{par}\circ Res$. Montrons que l'on a aussi
 
 (3)   $D^{par}\circ Rep(\varphi)=Rep(sgn\otimes \varphi)$ pour tout $\varphi\in {\cal R}$. 
 
 Preuve. Fixons $\gamma=(r',r'',N_{1},N_{2})\in \Gamma$, $\rho_{1}\in \hat{W}_{N_{1}}$, $\rho_{2}\in \hat{W}_{N_{2}}$ et consid\'erons l'\'el\'ement $\varphi=\rho_{1}\otimes \rho_{2}\in {\cal R}(\gamma)$. D'apr\`es Lusztig, le couple $(r',\rho_{1})$ param\`etre une repr\'esentation irr\'eductible $\pi_{1}$ d'un groupe fini $SO(2n_{1}+1,{\mathbb F}_{q})$, o\`u $n_{1}=N_{1}+r^{'2}+r'$ et ${\mathbb F}_{q}$ est le corps r\'esiduel de $F$. De m\^eme, le couple $(r'',\rho_{2})$ param\`etre une repr\'esentation irr\'eductible $\pi_{2}$ d'un groupe fini $SO(2n_{2},{\mathbb F}_{q})$, o\`u $n_{2}=N_{2}+r^{''2}$ (c'est la forme d\'eploy\'ee du groupe si $r''$ est pair, non d\'eploy\'ee si $r''$ est impair). Le terme $Rep(\varphi)$ n'est autre que $\pi_{1}\otimes \pi_{2}$. On d\'efinit usuellement une involution du groupe de Grothendieck des repr\'esentations de longueur finie de tels groupes finis (cf. \cite{C} 8.2 dans le cas d'un groupe connexe et \cite{DM} 3.10 dans le cas non connexe). C'est une somme altern\'ee de compos\'es de foncteurs de restriction et d'induction. D'apr\`es notre d\'efinition de \cite{W3} 1.7, $D^{par}(\pi_{1}\otimes \pi_{2})$ est le produit tensoriel des images de $\pi_{1}$ et $\pi_{2}$ par ces involutions multipli\'ees par des signes de sorte que ces images soient des repr\'esentations irr\'eductibles. D'autre part, pour tout $m\in {\mathbb N}$, on  d\'efinit une involution $D_{W_{m}}$ de ${\mathbb C}[\hat{W}_{m}]$ par une formule analogue: c'est  une somme altern\'ee de compos\'es de foncteurs de restriction et d'induction, cf. \cite{HL} corollaire 1. Les param\'etrages $(r',\rho_{1})\mapsto \pi_{1}$ et $(r'',\rho_{2})\mapsto \pi_{2}$ \'etant compatibles en un sens plus ou moins \'evident aux foncteurs de restriction et d'induction,   $D^{par}(\pi_{1}\otimes \pi_{2})$ est \'egal \`a l'image par $Rep$ de $\pm D_{W_{N_{1}}}(\rho_{1})\otimes D_{W_{N_{2}}}(\rho_{2})$, le signe \'etant choisi de sorte que ce terme soit le produit tensoriel de deux repr\'esentations irr\'eductibles. D'apr\`es \cite{HL} corollaire 1, on a $D_{W_{m}}(\rho)=\pm \rho\otimes sgn$  pour tout $m\in {\mathbb N}$ et tout $\rho\in \hat{W}_{m}$. Donc $D^{par}(\pi_{1}\otimes \pi_{2})$ est \'egal \`a l'image par $Rep$ de $(\rho_{1}\otimes sgn)\otimes(\rho_{2}\otimes sgn)$, c'est-\`a-dire de $sgn \otimes \varphi$. Cela prouve (3).

 Il est clair d'apr\`es sa d\'efinition  que l'endomorphisme $\rho\iota$ commute \`a la tensorisation $\varphi\mapsto sgn\otimes \varphi$. Alors la formule (2) se transforme en
 $$Res\circ \pi(\lambda^+,\epsilon^+,\lambda^-,\epsilon^-)=Rep\circ\rho\iota((\boldsymbol{\rho}_{\lambda^+,\epsilon^+}\otimes sgn)\otimes (\boldsymbol{\rho}_{\lambda^-,\epsilon^-}\otimes sgn)).$$ 
 En utilisant l'\'egalit\'e $Rep=k\circ {\cal F}$, on obtient finalement l'\'egalit\'e
 $$(4) \qquad \kappa_{\pi(\lambda^+,\epsilon^+,\lambda^-,\epsilon^-)}={\cal F}\circ\rho\iota((\boldsymbol{\rho}_{\lambda^+,\epsilon^+}\otimes sgn)\otimes (\boldsymbol{\rho}_{\lambda^-,\epsilon^-}\otimes sgn)).$$
 
 \bigskip
 
 \section{Symboles, partitions sp\'eciales, dualit\'e}
 
 \bigskip
 
 \subsection{Symboles}
 
Pour un couple $\Lambda=(X,Y)$ de sous-ensembles finis de ${\mathbb N}$, on d\'efinit le rang $rg(\Lambda)$ et le d\'efaut $def(\Lambda)$ par
$$rg(\Lambda)=S(X)+S(Y)-[(\vert X\vert +\vert Y\vert -1)^2/4],$$
o\`u $[.]$ d\'esigne la partie enti\`ere,
$$def(\Lambda)=sup(\vert X\vert -\vert Y\vert ,\vert Y\vert -\vert X\vert ).$$
On d\'efinit une relation d'\'equivalence entre couples de sous-ensembles finis de ${\mathbb N}$, engendr\'ee par $(X,Y)\sim (X',Y')$, o\`u
$$X'=\{x+1; x\in X\}\cup \{0\},\,\, Y'=\{y+1,y\in Y\}\cup\{0\}.$$
Le rang et le d\'efaut sont constants sur toute classe d'\'equivalence. On appelle symbole de d\'efaut impair une classe d'\'equivalence de couples $\Lambda=(X,Y)$ tels que $\vert X\vert >\vert Y\vert $ et $def(\Lambda)$ est impair. On appelle symbole de d\'efaut pair une classe d'\'equivalence de couples $\Lambda=(X,Y)$ tels que   $def(\Lambda)$ est pair (dans le cas pair, on n'impose pas $\vert X\vert \geq \vert Y\vert $). 

Soit $m\in {\mathbb N}$. On note $ S_{m,imp}$ l'ensemble des classes d'\'equivalence de symbole de d\'efaut impair et de rang $m$. Pour $\Lambda\in S_{m,imp}$, on pose $r(\Lambda)=(def(\Lambda)-1)/2$.  On note $S_{m,pair}$  l'ensemble des classes d'\'equivalence de symbole de d\'efaut pair et de rang $m$. Pour $\Lambda=(X,Y)\in S_{m,pair}$, on pose $r(\Lambda)=(\vert X\vert -\vert Y\vert )/2$. On a $def(\Lambda)=2\vert r(\Lambda)\vert $. 

{\bf Remarque.} La d\'efinition que l'on utilise ici des symboles de d\'efaut pair est diff\'erente de celle de \cite{W4} 1.2 o\`u l'on avait identifi\'e les couples $(X,Y)$ et $(Y,X)$. 
\bigskip

Notons $\Sigma_{m,imp}$ l'ensemble des triplets $(r,\alpha,\beta)$ o\`u $r\in {\mathbb N}$, $\alpha$ et $\beta$ sont des partitions et $r^2+r+S(\alpha)+S(\beta)=m$. Remarquons que, puisque les couples de partitions $(\alpha,\beta)$ v\'erifiant la relation pr\'ec\'edente param\`etrent les repr\'esentations irr\'eductibles de $W_{m-r^2-r}$, on peut identifier $\Sigma_{m,imp}$ \`a l'ensemble des couples $(r,\rho)$, o\`u $r\in {\mathbb N}$ v\'erifie $r^2+r\leq m$ et $\rho\in \hat{W}_{m-r^2-r}$. 
 On d\'efinit une application $symb:\Sigma_{m,imp}\to S_{m,imp}$ de la fa\c{c}on suivante. Soit $(r,\alpha,\beta)\in \Sigma_{m,imp}$. On  suppose que $\beta$ a $a$ termes pour un entier $a\geq0$ et que $\alpha$ en a $a+2r+1$. On pose $X=\alpha+\{a+2r,a+2r-1,...,0\}$, $Y=\beta+\{a-1,a-2,...,0\}$, $\Lambda=(X,Y)$. Alors, $symb(r,\alpha,\beta)=\Lambda$. Remarquons que $r=r(\Lambda)$. L'application $symb$ ainsi d\'efinie est une bijection de $\Sigma_{m,imp}$ sur $S_{m,imp}$.
 
 Notons $\Sigma_{m,pair}$ l'ensemble des triplets $(r,\alpha,\beta)$ o\`u $r\in {\mathbb Z}$, $\alpha$ et $\beta$ sont des partitions et $r^2+S(\alpha)+S(\beta)=m$.  On peut identifier $\Sigma_{m,pair}$ \`a l'ensemble des couples $(r,\rho)$, o\`u $r\in {\mathbb Z}$ v\'erifie $r^2\leq m$ et $\rho\in \hat{W}_{m-r^2}$. 
 On d\'efinit une application $symb:\Sigma_{m,pair}\to S_{m,pair}$ de la fa\c{c}on suivante. Soit $(r,\alpha,\beta)\in \Sigma_{m,pair}$.  On suppose que $\beta$ a $a$ termes et que $\alpha$ en a $a+2\vert r\vert $. Si $r\geq0$, on pose $X=\alpha+\{a+2r-1,a+2r-2,...,0\}$, $Y=\beta+\{a-1,a-2,...,0\}$. Si $r<0$, on pose $X=\beta+\{a-1,a-2,...,0\}$, $Y=\alpha+\{a+2\vert r\vert -1,a-2\vert r\vert -2,...,0\}$. On pose $\Lambda=(X,Y)$. Alors $symb(r,\alpha,\beta)=\Lambda$. Remarquons que $r=r(\Lambda)$. L'application $symb$ ainsi d\'efinie est une bijection de $\Sigma_{m,pair}$ sur $S_{m,pair}$. 
 
 Posons 
 $${\cal S}=\oplus_{n'+n''=n}{\mathbb C}[S_{n',imp}]\otimes {\mathbb C}[S_{n'',pair}].$$
D'apr\`es la construction de 1.2, l'espace ${\cal R}$ s'identifie \`a
 $$\oplus_{n'+n''=n}{\mathbb C}[\Sigma_{n',imp}]\otimes {\mathbb C}[\Sigma_{n'',pair}].$$
 Des bijections $symb$
 ci-dessus se d\'eduisent donc un isomorphisme encore not\'e $symb:{\cal R}\to {\cal S}$.
 
 \bigskip
 \subsection{Partitions sp\'eciales, cas symplectique}
 Soit $m\in {\mathbb N}$.  Une partition symplectique $\lambda\in {\cal P}^{symp}(2m)$ est dite sp\'eciale si $\lambda_{2j-1}$ et $\lambda_{2j}$ sont de m\^eme parit\'e pour tout $j\geq1$. On note ${\cal P}^{symp,sp}(2m)$  le sous-ensemble des partitions sp\'eciales. Soit $\lambda$ une telle partition sp\'eciale.  Consid\'erons l'ensemble  des \'el\'ements $i\in Jord^{bp}(\lambda)$ tels que $mult_{\lambda}(i)$ soit impair. S'il a un nombre pair d'\'el\'ements, on les note $i_{1}>i_{2}>...>i_{t}$. S'il a un nombre impair d'\'el\'ements,  on les note  $i_{1}>i_{2}>...>i_{t-1}$ et on pose $i_{t}=0$. Ainsi, $t$ est toujours pair. 
 On appelle intervalle de $\lambda$ un sous-ensemble de $Jord(\lambda)\cup \{0\}$ de l'une des formes suivantes
 
 $\{i\in Jord(\lambda)\cup\{0\}; i_{2h-1}\geq i\geq i_{2h}\}$ pour $h=1,...,t/2$;
 
 $\{i\}$ pour $i\in Jord^{bp}(\lambda)\cup\{0\}$ tel qu'il n'existe pas de $h=1,...,t/2$ de sorte que $ i_{2h-1}\geq i\geq i_{2h}$. 
 
Parce que $\lambda$ est sp\'eciale, on voit que les intervalles sont form\'es d'entiers pairs.
 On note $\tilde{Int}(\lambda)$ l'ensemble de ces intervalles. Il est ordonn\'e de fa\c{c}on naturelle:  $\Delta>\Delta'$ si et seulement si $i>i'$ pour tous $i\in \Delta$ et $i'\in \Delta'$. L'\'el\'ement minimal est celui qui contient $0$, on le note $\Delta_{min}$ et on pose $Int(\lambda)=\tilde{Int}(\lambda)-\{\Delta_{min}\}$. Pour $\Delta\in \tilde{Int}(\lambda)$, on note $J(\Delta)$ l'ensemble des $j\geq1$ tels que $\lambda_{j}\in \Delta$. C'est un intervalle de ${\mathbb N}$, qui est infini dans le cas $\Delta=\Delta_{min}$. On note $j_{min}(\Delta)$ le plus petit \'el\'ement de $J(\Delta)$ et, si $\Delta\not=\Delta_{min}$, on note $j_{max}(\Delta)$ le plus grand \'el\'ement de $J(\Delta)$. On v\'erifie que 
 
 $\{j_{min}(\Delta); \Delta\in \tilde{Int}(\lambda)\}$ est l'ensemble des $j\geq1$ tels que $j$ soit impair, $\lambda_{j}$ soit pair et $\lambda_{j-1}>\lambda_{j}$, avec la convention $\lambda_{0}=\infty$; 

 $\{j_{max}(\Delta); \Delta\in Int(\lambda)\}$ est l'ensemble des $j\geq1$ tels que $j$ soit pair, $\lambda_{j}$ soit pair et $\lambda_{j}>\lambda_{j+1}$.

 Par la correspondance de Springer, on associe \`a $(\lambda,1)$ une repr\'esentation irr\'eductible de $W_{m}$. Elle  est param\'etr\'ee par un couple $(\alpha(\lambda),\beta(\lambda))$.  On note $(X(\lambda),Y(\lambda))\in S_{m,imp}$ l'image de $(0,\alpha(\lambda),\beta(\lambda))$ par l'application $symb$. C'est un symbole sp\'ecial, c'est-\`a-dire que $\vert X(\lambda)\vert =\vert Y(\lambda)\vert +1$ et, si on note $X(\lambda)=(x_{1}>...>x_{a+1})$, $Y(\lambda)=(y_{1}>...>y_{a})$, on a $x_{1}\geq y_{1}\geq x_{2}\geq y_{2}\geq... \geq x_{a}\geq y_{a}\geq x_{a+1}$. On appelle famille de $\lambda$  l'ensemble des symboles $(X,Y)\in S_{m,imp}$ tels que, quitte \`a remplacer $(X,Y)$ et $(X(\lambda),Y(\lambda))$ par des symboles \'equivalents, on ait 
 
 (1) $X\cup Y=X(\lambda)\cup Y(\lambda)$, $X\cap Y=X(\lambda)\cap Y(\lambda)$. 
 
 On note $Fam(\lambda)$ la famille de $\lambda$.
  On montre que $S_{m,imp}$ est la r\'eunion disjointe des  $Fam(\lambda)$ quand $\lambda$ d\'ecrit l'ensemble ${\cal P}^{symp,sp}(2m)$. 
 
 Soit $\lambda\in {\cal P}^{symp,sp}(\lambda)$. On montre qu'il y a une unique bijection croissante $\Delta\mapsto x_{\Delta}$ de $\tilde{Int}(\lambda)$ sur $X(\lambda)-(X(\lambda)\cap Y(\lambda))$ et une unique bijection croissante $\Delta\mapsto y_{\Delta}$ de $Int(\lambda)$ sur $Y(\lambda)-(X(\lambda)\cap Y(\lambda))$. A un symbole $\Lambda=(X,Y) \in Fam(\lambda)$, on associe deux \'el\'ements $\tau\in ({\mathbb Z}/2{\mathbb Z})^{\tilde{Int}(\lambda)}$ et $\delta\in ({\mathbb Z}/2{\mathbb Z})^{Int(\lambda)}$ par les formules suivantes. On suppose les symboles choisis de sorte que (1) soit v\'erifi\'e. Alors, pour $\Delta\in \tilde{Int}(\lambda)$, resp. $\Delta\in Int(\lambda)$, on pose
 
 $$\tau(\Delta)=\vert \{\Delta'\in \tilde{Int}(\lambda); \Delta'\geq\Delta,\,\, x_{\Delta'}\in Y\}\vert $$
 $$+\vert \{\Delta'\in Int(\lambda); \Delta'>\Delta,\,\, y_{\Delta'}\in X\}\vert $$
 $$+ r(\Lambda)\,\,mod\,\,2{\mathbb Z};$$
 resp.
 $$\delta(\Delta)=\vert \{\Delta'\in Int(\lambda); \Delta'\geq\Delta,\,\, x_{\Delta'}\in Y\}\vert $$
 $$+\vert \{\Delta'\in Int(\lambda); \Delta'\geq\Delta,\,\, y_{\Delta'}\in X\}\vert \,\,mod\,\,2{\mathbb Z}.$$
Par cette construction, la famille  $Fam(\lambda)$ s'identifie \`a l'ensemble  des couples $(\tau,\delta)\in 
 ({\mathbb Z}/2{\mathbb Z})^{\tilde{Int}(\lambda)}\times ({\mathbb Z}/2{\mathbb Z})^{Int(\lambda)}$ tels que $\tau(\Delta_{min})=0$. On note ${\cal F}am(\lambda)$ cet ensemble. Pour  $(\tau,\delta)$ dans cet ensemble, provenant du symbole $\Lambda$, on pose $r(\tau,\delta)=r(\Lambda)$.  
 
 \bigskip
 
 \subsection{Partitions sp\'eciales, cas orthogonal impair}
  Soit $m\in {\mathbb N}$.  Une partition orthogonale $\lambda\in {\cal P}^{orth}(2m+1)$ est dite sp\'eciale si $\lambda_{2j}$ et $\lambda_{2j+1}$ sont de m\^eme parit\'e pour tout $j\geq1$. Il en r\'esulte que $\lambda_{1}$ est impair. On note ${\cal P}^{orth,sp}(2m+1)$  le sous-ensemble des partitions sp\'eciales. Soit $\lambda$ une telle partition sp\'eciale.  Les constructions du paragraphe pr\'ec\'edent s'appliquent.
  Par la correspondance de Springer, on associe \`a $(\lambda,1)$ une repr\'esentation irr\'eductible de $W_{m}$, puis un symbole appartenant \`a $S_{m,imp}$. Il est sp\'ecial. On d\'efinit   la famille de $\lambda$, que l'on note $Fam(\lambda)$.  On montre que $S_{m,imp}$ est la r\'eunion disjointe des  $Fam(\lambda)$ quand $\lambda$ d\'ecrit l'ensemble ${\cal P}^{orth,sp}(2m+1)$.  
  
  Remarquons que la conjonction des propri\'et\'es \'enonc\'ees ici et dans le paragraphe pr\'ec\'edent entra\^{\i}ne qu'il y a une bijection entre ${\cal P}^{symp,sp}(2m)$ et ${\cal P}^{orth,sp}(2m+1)$: $\lambda\in {\cal P}^{symp,sp}(2m)$ correspond \`a $\mu\in{\cal P}^{orth,sp}(2m+1)$ si et seulement si $Fam(\lambda)=Fam(\mu)$. En fait, nous utiliserons une autre bijection, la "dualit\'e", cf. 2.6. 
 
 \bigskip
 
 \subsection{Partitions sp\'eciales, cas orthogonal pair}
  Soit $m\in {\mathbb N}$.  Une partition orthogonale $\lambda\in {\cal P}^{orth}(2m)$ est dite sp\'eciale si $\lambda_{2j-1}$ et $\lambda_{2j}$ sont de m\^eme parit\'e pour tout $j\geq1$. On note ${\cal P}^{orth,sp}(2m)$  le sous-ensemble des partitions sp\'eciales. Soit $\lambda$ une telle partition sp\'eciale.  Consid\'erons l'ensemble  des \'el\'ements $i\in Jord^{bp}(\lambda)$ tels que $mult_{\lambda}(i)$ soit impair. On les note 
   $i_{1}>i_{2}>...>i_{t}$.  L'entier $t$ est forc\'ement  pair. 
 On appelle intervalle de $\lambda$ un sous-ensemble de $Jord(\lambda)$ de l'une des formes suivantes
 
 $\{i\in Jord(\lambda); i_{2h-1}\geq i\geq i_{2h}\}$ pour $h=1,...,t/2$;
 
 $\{i\}$ pour $i\in Jord^{bp}(\lambda)$ tel qu'il n'existe pas de $h=1,...,t/2$ de sorte que $ i_{2h-1}\geq i\geq i_{2h}$. 
 
Parce que $\lambda$ est sp\'eciale, on voit que les intervalles sont form\'es d'entiers impairs.
 On note $Int(\lambda)$ l'ensemble de ces intervalles. Comme dans le cas symplectique, il est ordonn\'e de fa\c{c}on naturelle. Pour $\Delta\in Int(\lambda)$, on definit $J(\Delta)$, $j_{min}(\Delta)$ et $j_{max}(\Delta)$ comme dans le cas symplectique.   On v\'erifie que 
 
 $\{j_{min}(\Delta); \Delta\in Int(\lambda)\}$ est l'ensemble des $j\geq1$ tels que $j$ soit impair, $\lambda_{j}$ soit impair et $\lambda_{j-1}>\lambda_{j}$, avec la convention $\lambda_{0}=\infty$; 

 $\{j_{max}(\Delta); \Delta\in Int(\lambda)\}$ est l'ensemble des $j\geq1$ tels que $j$ soit pair, $\lambda_{j}$ soit impair et $\lambda_{j}>\lambda_{j+1}$.

 Par la correspondance de Springer, on associe \`a $(\lambda,1)$ une repr\'esentation irr\'eductible de $W_{m}^D$. Elle  est param\'etr\'ee par un couple $(\alpha(\lambda),\beta(\lambda))$, qui n'est d\'etermin\'e qu'\`a l'ordre pr\`es.  On impose que $\alpha(\lambda)\geq \beta(\lambda)$ pour l'ordre lexicographique (s'il existe $j$ tel que $\alpha(\lambda)_{j}\not=\beta(\lambda)_{j}$, on a $\alpha(\lambda)_{j}>\beta(\lambda)_{j}$ pour le plus petit de ces entiers $j$). 
 On note $(X(\lambda),Y(\lambda))\in S_{m,pair}$ l'image de $(0,\alpha(\lambda),\beta(\lambda))$ par l'application $symb$. C'est un symbole sp\'ecial, c'est-\`a-dire que $\vert X(\lambda)\vert =\vert Y(\lambda)\vert $ et, si on note $X(\lambda)=(x_{1}>...>x_{a})$, $Y(\lambda)=(y_{1}>...>y_{a})$, on a $x_{1}\geq y_{1}\geq x_{2}\geq y_{2}\geq... \geq x_{a}\geq y_{a}$. On appelle famille de $\lambda$  l'ensemble des symboles $(X,Y)\in S_{m,pair}$ tels que, quitte \`a remplacer $(X,Y)$ et $(X(\lambda),Y(\lambda))$ par des symboles \'equivalents, on ait 
 
 (1) $X\cup Y=X(\lambda)\cup Y(\lambda)$, $X\cap Y=X(\lambda)\cap Y(\lambda)$. 
 
 On note $Fam(\lambda)$ la famille de $\lambda$.
 On montre que $S_{m,pair}$ est la r\'eunion disjointe des familles $Fam(\lambda)$ quand $\lambda$ d\'ecrit l'ensemble ${\cal P}^{orth,sp}(2m)$. 
 
 Soit $\lambda\in {\cal P}^{orth,sp}(\lambda)$. On montre qu'il y a une unique bijection croissante $\Delta\mapsto x_{\Delta}$ de $Int(\lambda)$ sur $X(\lambda)-(X(\lambda)\cap Y(\lambda))$ et une unique bijection croissante $\Delta\mapsto y_{\Delta}$ de $Int(\lambda)$ sur $Y(\lambda)-(X(\lambda)\cap Y(\lambda))$. A un symbole $\Lambda=(X,Y) \in Fam(\lambda)$, on associe deux \'el\'ements  $\tau,\delta\in ({\mathbb Z}/2{\mathbb Z})^{Int(\lambda)}$ par les m\^emes formules qu'en 2.2 (\`a ceci pr\`es qu'un $\tilde{Int}(\lambda)$ figurant dans ces derni\`eres est remplac\'e par $Int(\lambda)$). 
Par cette construction, la famille $Fam(\lambda)$ s'identifie \`a l'ensemble  des couples $(\tau,\delta)\in 
 ({\mathbb Z}/2{\mathbb Z})^{Int(\lambda)}\times ({\mathbb Z}/2{\mathbb Z})^{Int(\lambda)}$. On note ${\cal F}am(\lambda)$ cet ensemble. Pour  $(\tau,\delta)$ dans cet ensemble, provenant du symbole $\Lambda$, on pose $r(\tau,\delta)=r(\Lambda)$. Si $Int(\lambda)\not=\emptyset$, on note $\Delta_{min}$ son plus petit \'el\'ement et on v\'erifie que 
 
 (2) $\delta(\Delta_{min})\equiv r(\tau,\delta)\,\,mod\,\,2{\mathbb Z}$;
 
 \noindent  si $Int(\lambda)=\emptyset$, ${\cal F}am(\lambda)$ a un unique \'el\'ement $(\emptyset,\emptyset)$ et on a $r(\emptyset,\emptyset)=0$. 
 
 \bigskip
 
 \subsection{L'involution de Lusztig}
 Soit $m\in {\mathbb N}$. Soit $\lambda\in {\cal P}^{symp,sp}(2m)$. On note $\Lambda(\lambda)=(X(\lambda),Y(\lambda))$ le symbole sp\'ecial associ\'e \`a $\lambda$. On repr\'esente tout \'el\'ement de  la famille de $\lambda$ par un  symbole $\Lambda=(X,Y)$ v\'erifiant la condition 2.2(1).   Soient $\Lambda=(X,Y),\,\, \Lambda'=(X',Y')$ deux \'el\'ements de $Fam(\lambda)$. On pose
 $$<\Lambda,\Lambda'>=r(\Lambda)+r(\Lambda')+\vert X\cap X'\cap Y(\lambda)\vert +\vert Y\cap Y'\cap X(\lambda)\vert \,\,mod\,\,2{\mathbb Z}.$$
 Cela d\'efinit une application:
 $$<.,.>:Fam(\lambda)\times Fam(\lambda)\to {\mathbb Z}/2{\mathbb Z}.$$
 On d\'efinit un automorphisme ${\cal F}$ de l'espace ${\mathbb C}[Fam(\lambda)]$ par la formule
 $${\cal F}(\Lambda)=\vert Fam(\lambda)\vert ^{-1/2}\sum_{\Lambda'\in Fam(\lambda)}(-1)^{<\Lambda,\Lambda'>}\Lambda',$$
 les symboles \'etant ici identifi\'es  aux \'el\'ements de base de ${\mathbb C}[Fam(\lambda)]$. On v\'erifie qu'il est involutif. D'apr\`es ce que l'on a dit en 2.2, l'espace ${\mathbb C}[S_{m,imp}]$ est somme directe des sous-espaces ${\mathbb C}[Fam(\lambda)]$ quand $\lambda$ d\'ecrit ${\cal P}^{symp,sp}(2m)$. On note ${\cal F}$ l'automorphisme de ${\mathbb C}[S_{m,imp}]$ qui est la somme directe des automorphismes de ces sous-espaces que l'on vient de construire.
 
 Pour $\lambda\in {\cal P}^{orth,sp}(2m)$, on d\'efinit exactement de la m\^eme fa\c{c}on un automorphisme ${\cal F}$ de ${\mathbb C}[Fam(\lambda)]$. Puis, par somme directe, on en d\'eduit un automorphisme de ${\mathbb C}[S_{m,pair}]$. 
 
 Dans le cas orthogonal pair, on dispose d'une involution $\sigma$ de $Fam(\lambda)$: si $\Lambda=(X,Y)$, $\sigma(\Lambda)=(Y,X)$. Pour $\Lambda,\Lambda'\in Fam(\lambda)$, on v\'erifie la formule
 
 (1) $<\sigma(\Lambda),\Lambda'>\equiv r(\Lambda')+<\Lambda,\Lambda'>\,\,mod\,\,2{\mathbb Z}$.
 
 Rappelons que 
 $${\cal S}=\oplus_{n',n''\in {\mathbb N}, n'+n''=n}{\mathbb C}[S_{n',imp}]\otimes {\mathbb C}[S_{n'',pair}].$$
 On a d\'efini des automorphismes ${\cal F}$ de chacun des espaces qui interviennent ici. Par produit tensoriel et sommation, on en d\'eduit un automorphisme ${\cal F}$ de ${\cal S}$. On a d\'efini en 2.1 un isomorphisme $symb:{\cal R}\to {\cal S}$. Par celui-ci, on transporte l'automorphisme ${\cal F}$ de ${\cal S}$ en un automorphisme ${\cal F}$ de ${\cal R}$. C'est l'automorphisme de Lusztig introduit en 1.4.
 
 \bigskip
 
 \subsection{Dualit\'e}
 Soit $m\in {\mathbb N}$. On a introduit en 2.1 l'ensemble $\Sigma_{m,imp}$, que l'on voit ici comme un ensemble de couples $(r,\rho)$, o\`u $\rho\in \hat{W}_{m-r^2-r}$. On d\'efinit une involution de cet ensemble par $(r,\rho)\mapsto (r,\rho\otimes sgn)$. Transportons-la en une involution de $S_{m,imp}$ par la bijection $symb$. On note $d$ l'involution obtenue. Elle se calcule ainsi. Soit $\Lambda=(X,Y)\in S_{m,imp}$.  Fixons un entier $d$ plus grand que tous les termes de $X$ et $Y$. Posons
 $$X'=\{d,...,0\}-\{d-y;y\in Y\},\,\,Y'=\{d,...,0\}-\{d-x;x\in X\}.$$
 Alors $d(\Lambda)=(X',Y')$. Cette formule montre que $d$ conserve la d\'ecomposition en familles, c'est-\`a-dire que si $\Lambda$ et $\Lambda'$ sont dans une m\^eme famille, alors $d(\Lambda)$ et $d(\Lambda')$ sont aussi dans une m\^eme famille. On d\'efinit une application appel\'ee dualit\'e $d:{\cal P}^{symp,sp}(2m)\to {\cal P}^{orth,sp}(2m+1)$ ou $d: {\cal P}^{orth,sp}(2m+1)\to {\cal P}^{symp,sp}(2m)$ par la condition $Fam(d(\lambda))=d(Fam(\lambda))$.  Les deux applications sont inverses l'une de l'autre.
 
 Ces dualit\'es s'\'etendent en des applications $d:{\cal P}^{symp}(2m)\to {\cal P}^{orth,sp}(2m+1)$ ou $d: {\cal P}^{orth}(2m+1)\to {\cal P}^{symp,sp}(2m)$. Rappelons la d\'efinition de la premi\`ere, celle de la seconde \'etant similaire. Soit $\lambda\in {\cal P}^{symp}(2m)$. La correspondance de Springer associe au couple $(\lambda,1)\in \boldsymbol{{\cal P}^{symp}}(2m)$ une repr\'esentation $\rho_{\lambda,1}$ de $W_{m}$. Le couple $(0,\rho)$ appartient \`a $\Sigma_{m,imp}$. Il existe une unique partition symplectique sp\'eciale, que l'on note $sp(\lambda)$, dont la famille contient le  symbole $symb(0,\rho)$. En fait, on montre que $sp(\lambda)$ est la plus petite partition symplectique sp\'eciale $\lambda'$ telle que $\lambda\leq \lambda'$. On pose $d(\lambda)=d(sp(\lambda))$. Cette dualit\'e est d\'ecroissante: $\lambda\leq \lambda'$ entra\^{\i}ne $d(\lambda')\leq d(\lambda)$. 
 
 On peut remplacer $\Sigma_{m,imp}$ par $\Sigma_{m,pair}$ dans la construction ci-dessus. On obtient une dualit\'e $d$ qui est une involution de ${\cal P}^{orth,sp}(2m)$. Celle-ci s'\'etend en une application $d:{\cal P}^{orth}(2m)\to {\cal P}^{orth,sp}(2m)$, qui est d\'ecroissante.

 \subsection{Calcul de $d(\lambda)$}
Soient $m\in {\mathbb N}$ et  $\lambda\in {\cal P}^{symp}(2m)$. Pour $i\in Jord(\lambda)\cup\{0\}$, notons $J(i)$ l'ensemble des indices $j\geq1$ tels que $\lambda_{j}=i$. C'est un intervalle de ${\mathbb N}-\{0\}$, infini si $i=0$. On note$j_{min}(i)$ son plus petit \'el\'ement et, si $i\not=0$, $j_{max}(i)$ son plus grand \'el\'ement. Consid\'erons l'ensemble  des \'el\'ements $i$ de $Jord^{bp}(\lambda)$ tels que $mult_{\lambda}(i)$ soit impaire. Comme en 2.2, si cet ensemble a un nombre pair d'\'el\'ements, on les note $i_{1}>...>i_{t}$. S'il a un nombre impair d'\'el\'ements, on les note $i_{1}>...>i_{t-1}$ et on pose $i_{t}=0$. Pour $h=1,...t$, on v\'erifie que

 $j_{min}(i_{h})\equiv h\,\,mod\,\,2{\mathbb Z}$ et, si $i_{h}\not=0$, $j_{max}(i_{h})\equiv h\,\,mod\,\,2{\mathbb Z}$.

Consid\'erons les \'el\'ements de $Jord(\lambda)\cup\{0\}$ qui n'interviennent pas dans la suite $i_{1},...,i_{t}$, c'est-\`a-dire les $i\in Jord(\lambda)$ tels que $mult_{\lambda}(i)$ soit pair et  aussi $0$ dans le cas o\`u $i_{t}\not=0$.  Notons ${\cal J}(\lambda)$ cet ensemble. On d\'ecompose ${\cal J}(\lambda)$ en union disjointe ${\cal J}'(\lambda)\sqcup {\cal J}''(\lambda)$: ${\cal J}''(\lambda)$ est l'ensemble des $i\in {\cal J}(\lambda)$ tels qu'il existe $h=1,...t/2$ de sorte que $i_{2h-1}>i>i_{2h}$; ${\cal J}'(\lambda)$ est son compl\'ementaire. On v\'erifie que

 pour $i\in {\cal J}'(\lambda)$, $j_{min}(i)$ est impair et, si $i\not=0$, $j_{max}(i)$ est pair;

 pour $i\in {\cal J}''(\lambda)$, $j_{min}(i)$ est pair et, si $i\not=0$, $j_{max}(i)$ est impair.

Notons 

$P^+(\lambda)$ l'ensemble  des entiers  impairs $j\geq1$ tels que   $\lambda_{j}$ est pair et $\lambda_{j-1}>\lambda_{j}$, avec la convention $\lambda_{0}=\infty$;

$P^-(\lambda)$  l'ensemble  des entiers pairs $j\geq2$ tels que   $\lambda_{j}$ est pair et $\lambda_{j}>\lambda_{j+1}$; 

$Q^+(\lambda)$ l'ensemble des entiers pairs $j\geq2$ tels que  $\lambda_{j}$ est impair et $\lambda_{j-1}>\lambda_{j}$; 

$Q^-(\lambda)$ l'ensemble des entiers impairs $j\geq1$ tels que $\lambda_{j}$ est  impair et $\lambda_{j}>\lambda_{j+1}$. 

Ces ensembles sont disjoints. A l'aide des propri\'et\'es pr\'ec\'edentes, on voit que

$P^+(\lambda)$ est l'ensemble des $j_{min}(i)$ pour $i=i_{h}$ avec  $h$ impair, ou pour un \'el\'ement pair $i\in {\cal J}'(\lambda)$;

$P^-(\lambda)$ est l'ensemble des $j_{max}(i)$ pour  $i=i_{h}$ avec  $h$ pair et $i_{h}\not=0$ ou pour  un \'el\'ement pair non nul $i\in {\cal J}'(\lambda)$;

$Q^+(\lambda)$ est l'ensemble des $j_{min}(i)$ pour un \'el\'ement impair $i\in  {\cal J}''(\lambda)$; 

$Q^-(\lambda)$ est l'ensemble des $j_{max}(i)$ pour un \'el\'ement impair $i\in  {\cal J}''(\lambda)$. 

Ces ensembles sont disjoints. 
Les \'el\'ements de $P^+(\lambda)\cup P^-(\lambda)$ apparaissent presque tous par paires. Un \'el\'ement de $P^+(\lambda)$ de la forme $j_{min}(i)$ pour $i=i_{h}$ avec $h$ impair est suivi de l'\'el\'ement $j_{max}(i_{h+1})\in P^-(\lambda)$ sauf si $i_{h+1}=0$. Un \'el\'ement de $P^+(\lambda)$ de la forme $j_{min}(i)$ pour un \'el\'ement pair $i\in {\cal J}'(\lambda)$ est suivi de l'\'el\'ement $j_{max}(i)\in P^-(\lambda)$ sauf si $i=0$. Le plus petit \'el\'ement de $P^+(\lambda)$ est $j_{min}(i_{t-1})$ si $i_{t}=0$ ou $j_{min}(0)$ si $i_{t}\not=0$. Il n'est suivi d'aucun \'el\'ement de $P^-(\lambda)$. 
Il en r\'esulte que $\vert P^+\lambda)\vert =\vert P^-(\lambda)\vert +1$  et que, si on note ces ensembles
$$P^+(\lambda)=\{p^+_{1}<...<p^+_{a+1}\},\,\,P^-(\lambda)=\{p^-_{1}<...<p^-_{a}\},$$
on a les relations
$$p^+_{1}<p^-_{1}<p^+_{2}<p^-_{2}<...<p_{a}^+<p_{a}^-<p_{a+1}^+.$$
On voit de m\^eme que $\vert Q^+(\lambda)\vert =\vert Q^-(\lambda)\vert $ et que, si on note ces ensembles
$$Q^+(\lambda)=\{q_{1}^+<...<q_{b}^+\},\,\,Q^-(\lambda)=\{q_{1}^-<...<q_{b}^-\},$$
on a les relations
$$q^+_{1}<q^-_{1}<q^+_{2}<q^-_{2}<...<q^+_{b}<q^-_{b}.$$
Remarquons que, pour $j\in Q^+(\lambda)$, on a $\lambda_{j}=\lambda_{j+1}$. En effet, $\lambda_{j}$ est impair, donc $mult_{\lambda}(\lambda_{j})$ est pair. Mais on a aussi $\lambda_{j-1}>\lambda_{j}$,  d'o\`u l'\'egalit\'e cherch\'ee. De m\^eme, pour $j\in Q^-(\lambda)$, on a $j\geq2$ et $\lambda_{j-1}=\lambda_{j}$.

 D\'efinissons deux suites de nombres $\zeta(\lambda)=(\zeta(\lambda)_{1},\zeta(\lambda)_{2},...)$ et $s(\lambda)=(s(\lambda)_{1},s(\lambda)_{2},...)$ par les \'egalit\'es
$$\zeta(\lambda)_{j}=\left\lbrace\begin{array}{cc}1,& \text{ \,\,si\,\,}j\in P^+(\lambda),\\ -1,&\text{\,\,si\,\,}j\in P^-(\lambda),\\0,&\text{\,\,sinon;}\\ \end{array}\right.$$
$$s(\lambda)(j)=\left\lbrace\begin{array}{cc}1,& \text{ \,\,si\,\,}j\in Q^+(\lambda),\\ -1,&\text{\,\,si\,\,}j\in Q^-(\lambda),\\0,&\text{\,\,sinon.}\\ \end{array}\right.$$

\ass{Lemme}{(i)  On a l'\'egalit\'e $sp(\lambda)=\lambda+s(\lambda)$;

(ii) On a l'\'egalit\'e $^td(\lambda)=\lambda+\zeta(\lambda)$.}

Preuve. Posons $\nu=\lambda+s(\lambda)$. Montrons  que $\nu$ est une partition, c'est-\`a-dire que $\nu_{j}\geq \nu_{j+1}$ pour tout $j\geq1$. Puisque le couple $( \nu_{j},\nu_{j+1})$ s'obtient en ajoutant \`a $(\lambda_{j},\lambda_{j+1})$ un couple qui appartient \`a $\{-1,0,1\}\times\{-1,0,1\}$ et puisque $\lambda_{j}\geq \lambda_{j+1}$, la conclusion est claire sauf si le couple ajout\'e est $(-1,0)$, $(0,1)$ ou $(-1,1)$. Le premier cas se produit seulement si $j\in Q^-(\lambda)$. Dans ce cas on a $\lambda_{j}>\lambda_{j+1}$ par d\'efinition de $Q^-(\lambda)$ et alors $\lambda_{j}-1\geq \lambda_{j+1}$. Le deuxi\`eme cas se produit seulement si $j+1\in Q^+(\lambda)$. Dans ce cas, on a encore $\lambda_{j}>\lambda_{j+1}$ par d\'efinition de $Q^+(\lambda)$ et alors $\lambda_{j}\geq \lambda_{j+1}+1$. Le dernier cas se produit quand $j\in Q^-(\lambda)$ et $j+1\in Q^+(\lambda)$. On a encore $\lambda_{j}>\lambda_{j+1}$ De plus, $\lambda_{j}$ et $\lambda_{j+1}$ sont tous deux impairs. Donc $\lambda_{j}\geq \lambda_{j+1}+2$. Alors $\lambda_{j}-1\geq \lambda_{j+1}+1$.

 L'\'egalit\'e des nombres d'\'el\'ements de $Q^+(\lambda)$ et de $Q^-(\lambda)$ et la d\'efinition de $s(\lambda)$ entra\^{\i}nent que $S( \nu)=2n$. Une partition $\mu$ de $2n$ est symplectique et sp\'eciale si et seulement si, pour tout entier $j\geq1$ impair, $\mu_{j}$ et $\mu_{j+1}$ sont de m\^eme parit\'e et, si ces nombres sont impairs, alors ils sont \'egaux. Cela \'equivaut \`a: pour tout $j\geq1$ impair, si $\mu_{j}$ ou $\mu_{j+1}$ est impair, alors $\mu_{j}=\mu_{j+1}$. Montrons que $\nu$ v\'erifie cette condition. Soit  un entier $j\geq1$  impair, supposons $ \nu_{j}$ impair. L'entier $j$ n'appartient pas \`a $Q^+(\lambda)$ car il est pair. Il n'appartient pas \`a $Q^-(\lambda)$: sinon $\lambda_{j}$ serait impair et $ \nu_{j}=\lambda_{j}-1$  serait pair.  Donc $s(\lambda)_{j}=0$ et $ \nu_{j}=\lambda_{j}$. Puisque $j$ est impair, que $ \lambda_{j}= \nu_{j}$ est impair et que $j\not\in Q^-(\lambda)$, on a $\lambda_{j}=\lambda_{j+1}$. Cette \'egalit\'e entra\^{\i}ne que  $j+1$ n'appartient pas \`a $Q^+(\lambda)$. Il n'appartient pas non plus \`a $Q^-(\lambda)$ car $j+1$ est pair. Donc $s(\lambda)_{j+1}=0$, $ \nu_{j+1}=\lambda_{j+1}$  et on conclut $ \nu_{j}=\nu_{j+1}$. Une preuve analogue montre que, si $ \nu_{j+1}$ est impair, on a $ \nu_{j}=\nu_{j+1}$. Donc $ \nu$ est symplectique et sp\'eciale. 

Soit $j\geq1$. Par construction et d'apr\`es la description des ensembles $Q^+(\lambda)$ et $Q^-(\lambda)$, $S_{j}( \nu)=S_{j}(\lambda)$ sauf s'il existe un \'el\'ement impair $i\in  {\cal J}''(\lambda)$ tel que   $j_{min}(i)\leq j<j_{max}(i)$.  
 S'il existe un tel $i$, on a $S_{j}( \nu)=S_{j}(\lambda)+1$. Cela montre que $\lambda\leq  \nu$. Soit $\mu\in {\cal P}^{symp,sp}(2n)$ telle que $\lambda\leq \mu$. On a $S_{j}(\lambda)\leq S_{j}(\mu)$. Cela entra\^{\i}ne 

(1) $S_{j}( \nu)\leq S_{j}(\mu)$,

\noindent  sauf s'il existe  $i$ comme ci-dessus.  Supposons qu'il existe un tel $i$ et notons simplement $j^+=j_{min}(i)$, $j^-=j_{max}(i)$. On a $j^+\in Q^+(\lambda)$ et $j^-\in Q^-(\lambda)$. Par d\'efinition de $j_{min}(i)$, on a  $\lambda_{j^+-1}> \lambda_{j^+}$. Les entiers $\lambda_{1},...,\lambda_{j^+-1}$ sont tous les entiers strictement sup\'erieurs \`a $\lambda_{j^+}=i$ qui interviennent dans $\lambda$, compt\'es avec leurs multiplicit\'es. Puisque $\lambda$ est symplectique, l'entier  $S_{j^+-1}(\lambda)$ est pair.  Puisque $j^+\in Q^+(\lambda)$, $j^+$ est pair et on sait que $i$ est impair. Donc $S_{j^+}(\lambda)=S_{j^+-1}(\lambda)+i$ est impair et aussi $S_{j^+}(\lambda)+j^+$. Puisque $\lambda_{j'}=i$ est   impair pour $j'\in \{j^+,...,j\}$, on voit que $S_{j}(\lambda)+j$ est aussi impair. Supposons $j$ pair. Alors $S_{j}(\lambda)$ est impair. Or le fait que $\mu$ soit sp\'eciale entra\^{\i}ne que $S_{j}(\mu)$ est pair. L'in\'egalit\'e $S_{j}(\lambda)\leq S_{j}(\mu)$ est alors  stricte et on conclut $S_{j}(\nu)=S_{j}(\lambda)+1\leq S_{j}(\mu)$.  Supposons $j$ impair. On sait que $j^+$ est pair et que $j^-$ est impair par d\'efinition des ensembles $Q^+(\lambda)$ et $Q^-(\lambda)$. L'hypoth\`ese $j\in \{j^+,...,j^--1\}$ entra\^{\i}ne alors que $j-1\in \{j^+,...,j^--1\}$ et $j,j+1\in \{j^++1...,j^--1\}$. L'\'egalit\'e (1) est d\'emontr\'ee pour $j-1$ et pour $j+1$ puisque ces entiers sont pairs.  D'o\`u
$$S_{j-1}( \nu)\leq S_{j-1}(\mu)\text{\,\,et\,\,}S_{j+1}( \nu)\leq S_{j+1}(\mu).$$
De plus, puisque $j$ et $j+1$ appartiennent tous deux \`a $\{j^++1...,j^--1\}$, on a $\nu_{j}= i= \nu_{j+1}$. La seconde in\'egalit\'e ci-dessus se r\'ecrit
$$S_{j-1}( \nu)+2i\leq S_{j-1}(\mu)+\mu_{j}+\mu_{j+1}.$$
On additionne cette in\'egalit\'e avec la premi\`ere in\'egalit\'e ci-dessus et on obtient
$$S_{j-1}( \nu)+i\leq S_{j-1}(\mu)+(\mu_{j}+\mu_{j+1})/2.$$
Evidemment, $(\mu_{j}+\mu_{j+1})/2\leq \mu_{j}$, d'o\`u
$$S_{j-1}( \nu)+i\leq S_{j-1}(\mu)+\mu_{j}.$$
Le membre de gauche est $S_{j}( \nu)$, celui de droite $S_{j}(\mu)$. Cela ach\`eve de d\'emontrer (1).

L'in\'egalit\'e (1) signifie que $ \nu\leq \mu$. On a ainsi d\'emontr\'e  que $\nu$ \'etait la plus petite partition symplectique sp\'eciale $\mu$ telle que $\lambda\leq \mu$. Cette propri\'et\'e caract\'erise $sp(\lambda)$, ce qui d\'emontre le (i) de l'\'enonc\'e.

Prouvons maintenant que

(2) $\zeta(\lambda)=\zeta( \nu)+s(\lambda)$.

Par d\'efinition de ces suites, cela \'equivaut aux \'egalit\'es

(3) $P^+( \nu)=P^+(\lambda)\cup Q^-(\lambda)$,
$P^-( \nu)=P^-(\lambda)\cup Q^+(\lambda)$.

 La premi\`ere \'egalit\'e concerne des indices $j\geq1$ impairs. Soit un tel $j$. Supposons d'abord  $j\in P^+(\lambda)$. On a $ \nu_{j}=\lambda_{j}$ et ce terme est pair. On a de plus $\lambda_{j-1}> \lambda_{j}$. Si $j=1$, on a trivialement $ \nu_{j-1}> \nu_{j}$ et on conclut $j\in P^+( \nu)$. Supposons $j\geq2$. Certainement, $j-1\not\in Q^-(\lambda)$ puisque $j-1$ est pair. Donc $ \nu_{j-1}\geq \lambda_{j-1}$, d'o\`u $ \nu_{j-1}>  \nu_{j}$. Alors $j$ appartient \`a $P^+( \nu)$. Supposons maintenant $j\in Q^-(\lambda)$. Alors $ \nu_{j}=\lambda_{j}-1$ et $\lambda_{j}$ est impair, donc $ \nu_{j}$ est pair.   Comme on l'a vu, l'hypoth\`ese $j\in Q^-(\lambda)$ entra\^{\i}ne   $\lambda_{j-1}=\lambda_{j}$. Comme ci-dessus, $j-1$ n'appartient pas \`a $Q^-(\lambda)$ donc $ \nu_{j-1}\geq \lambda_{j-1}=\lambda_{j}>  \nu_{j}$. D'o\`u $j\in P^+ (\nu)$. Supposons enfin que $j\in P^+( \nu)$. En particulier $ \nu_{j}$ est pair. Si $\lambda_{j}$ est impair, on a n\'ecessairement $s(\lambda)_{j}\not=0$ et, puisque $j$ est impair, $j$ appartient \`a $Q^-(\lambda)$. Supposons $\lambda_{j}$ pair. Alors $s(\lambda)_{j}$ est pair donc nul. Si $j=1$, on a trivialement $\lambda_{j-1}> \lambda_{j}$ et $j$ appartient \`a $P^+(\lambda)$. Supposons $j\geq2$. Puisque $j\in P^+( \nu)$, on a $ \nu_{j-1}>  \nu_{j}$, autrement dit $\lambda_{j-1}+s(\lambda)_{j-1}>\lambda_{j}$. 
On n'a pas $j-1\in Q^+(\lambda)$ car   cette relation entra\^{\i}ne que $\lambda_{j-1}=\lambda_{j}$ est impair contrairement \`a l'hypoth\`ese. Donc $s(\lambda)_{j-1}\leq0$. L'in\'egalit\'e $\lambda_{j-1}+s(\lambda)_{j-1}>\lambda_{j}$ entra\^{\i}ne alors $\lambda_{j-1}> \lambda_{j}$, donc $j\in P^+(\lambda)$. Cela d\'emontre la premi\`ere \'egalit\'e de (3). La seconde se d\'emontre de fa\c{c}on analogue. Cela prouve (3), d'o\`u (2).

Dans le cas o\`u $\lambda$ est sp\'eciale, on a d\'efini l'ensemble d'intervalles $\tilde{Int}(\lambda)$.  On   voit que $P^+(\lambda)$ est l'ensemble des $j_{min}(\Delta)$ quand $\Delta$ d\'ecrit $\tilde{Int}(\lambda)$ et que $P^-(\lambda)$ est l'ensemble des $j_{max}(\Delta)$ pour $\Delta\in Int(\lambda) $. Alors $\zeta(\lambda)$ est la suite que l'on a d\'efinie en \cite{W4} 1.6. On a d\'emontr\'e dans cette r\'ef\'erence l'\'egalit\'e $^td(\lambda)=\lambda+\zeta(\lambda)$. Supprimons l'hypoth\`ese que $\lambda$ est sp\'eciale. Par d\'efinition, $d(\lambda)=d(sp(\lambda))$. D'o\`u
$$^td(\lambda)=^td(sp(\lambda))=sp(\lambda)+\zeta(sp(\lambda)).$$
Puisque $sp(\lambda)=\nu=\lambda+s(\lambda)$, l'\'egalit\'e (2) entra\^{\i}ne la deuxi\`eme  assertion de l'\'enonc\'e. $\square$

Soit maintenant $\lambda\in {\cal P}^{orth}(2m)$. On d\'efinit $P^+(\lambda)$ et $P^-(\lambda)$ en \'echangeant les conditions de parit\'e sur les $\lambda_{j}$. C'est-\`a-dire

$P^+(\lambda)$ l'ensemble  des entiers  impairs $j\geq1$ tels que   $\lambda_{j}$ est impair et $\lambda_{j-1}>\lambda_{j}$, avec la convention $\lambda_{0}=\infty$;

$P^-(\lambda)$  l'ensemble  des entiers pairs $j\geq2$ tels que   $\lambda_{j}$ est impair et $\lambda_{j}>\lambda_{j+1}$.

Dans ce cas, on a $\vert P^+(\lambda)\vert =\vert P^-(\lambda)\vert $. 
On d\'efinit la suite $\zeta$ comme plus haut.  Nous aurons besoin de l'analogue du (ii) du lemme ci-dessus, mais seulement dans le cas o\`u $\lambda$ est sp\'eciale. C'est-\`a-dire

(4) si $\lambda\in {\cal P}^{orth,sp}(2m)$, on a $^td(\lambda)=\lambda+\zeta(\lambda)$.

Cf. \cite{W4} 1.7. 

\bigskip

\section{Induction endoscopique}

   \bigskip
   
   \subsection{L'induite endoscopique de deux partitions sp\'eciales}
   Soient $n_{1},n_{2}\in {\mathbb N}$ tels que $n_{1}+n_{2}=n$ et soient $\lambda_{1}\in {\cal P}^{symp,sp}(2n_{1})$ et $\lambda_{2}\in {\cal P}^{orth,sp}(2n_{2})$. Pour un indice $j\geq1$, on dit que $\lambda_{1,j}$, resp. $\lambda_{2,j}$, est de bonne parit\'e si $\lambda_{1,j}$ est pair, resp. $\lambda_{2,j}$ est impair. Notons
   
   $J^+$ l'ensemble des $j\geq1$ tels que $j$ soit impair, $\lambda_{1,j}$ et $\lambda_{2,j}$ soient de bonne parit\'e et il existe $d=1,2$ de sorte que $\lambda_{d,j-1}>\lambda_{d,j}$ (avec toujours la convention $\lambda_{d,0}=\infty$);
 
   $J^-$ l'ensemble des $j\geq1$ tels que $j$ soit pair, $\lambda_{1,j}$ et $\lambda_{2,j}$ soient de bonne parit\'e et il existe $d=1,2$ de sorte que $\lambda_{d,j}>\lambda_{d,j+1}$.
   
   On v\'erifie que $\vert J^+\vert =\vert J^-\vert $ et que, si on note $j_{1}^+<...<j_{a}^+$ les \'el\'ements de $J^+$ et $j_{1}^-<...<j_{a}^-$ ceux de $J^-$, on a $j_{1}^+<j_{1}^-<j_{2}^+<...<j_{a}^+<j_{a}^-$. On d\'efinit une suite $\xi=(\xi_{1},\xi_{2},...)$ de nombres entiers par $\xi_{j}=1$ si $j\in J^+$, $\xi_{j}=-1$ si $j\in J^-$ et $\xi_{j}=0$ si $j\not\in J^+\cup J^-$. On pose
   $$\lambda=\lambda_{1}+\lambda_{2}+\xi.$$
   C'est une partition symplectique de $2n$, appel\'ee l'induite endoscopique de $\lambda_{1}$ et $\lambda_{2}$.
 
 Pour unifier les notations, on pose $\tilde{Int}(\lambda_{2})=Int(\lambda_{2})$. Pour $d=1,2$, posons $J_{d,min}=\{j_{min}(\Delta); \Delta\in \tilde{Int}(\lambda_{d})\}$, $J_{d,max}=\{j_{max}(\Delta);\Delta\in Int(\lambda_{d})\}$. On note ${\cal J}^+=J_{1,min}\cap J_{2,min}$, ${\cal J}^-=J_{1,max}\cap J_{2,max}$,
 $${\cal J}=J_{1,min}\cup J_{2,min}\cup J_{1,max}\cup J_{2,max}\cup\{\infty\}.$$
 Appelons intervalle relatif d'indices  un sous-ensemble de ${\mathbb N}-\{0\}$ de l'une des formes suivantes
 
 (1) $\{j\}$ pour $j\in {\cal J}^+\cup {\cal J}^-$;
 
 (2) $\{j,...,j'\}$ o\`u $j$ et $j'$ sont deux \'el\'ements cons\'ecutifs de ${\cal J}$ tels qu'il existe un unique $d=1,2$ de sorte que $\{j,...,j'\}\subset J(\Delta)$ pour un $\Delta\in \tilde{Int}(\lambda_{d})$. 
 
 Pour un  intervalle relatif d'indices $J$, on pose $D(J)=\{\lambda_{j}; j\in J\}$. On appelle intervalle  de $\lambda$ relatif \`a $(\lambda_{1},\lambda_{2})$ un sous-ensemble de $Jord(\lambda)\cup\{0\}$ de la forme $D(J)$. On note $\tilde{Int}_{\lambda_{1},\lambda_{2}}(\lambda)$ l'ensemble de ces intervalles relatifs. On montre qu'ils sont disjoints, form\'es de nombres pairs et que $\tilde{Int}_{\lambda_{1},\lambda_{2}}(\lambda)$ est une partition de $Jord^{bp}(\lambda)\cup\{0\}$. Pour un intervalle relatif $D$, on note $J(D)$ l'intervalle relatif d'indices $J$ tel que $D=D(J)$.  
  Les intervalles relatifs sont ordonn\'es de fa\c{c}on naturelle: $D>D'$ si et seulement si $i>i'$ pour tous $i\in D$, $i'\in D'$. L'intervalle minimal est celui qui contient $0$, on le note $D_{min}$ et on pose $Int_{\lambda_{1},\lambda_{2}}(\lambda)=\tilde{Int}_{\lambda_{1},\lambda_{2}}(\lambda)-\{D_{min}\}$. Pour $D\in \tilde{Int}_{\lambda_{1},\lambda_{2}}(\lambda)$, on note $j_{min}(D)$, resp. $j_{max}(D)$, le plus petit, resp. grand, \'el\'ement de $J(D)$ (on consid\`ere que $j_{max}(D_{min})=\infty$). 
    
  Montrons que
  
  (3) pour tout $j\in {\cal J}$, il existe un unique intervalle relatif $D$ tel que $j\in \{j_{min}(D),j_{max}(D)\}$.
  
  Preuve. L'unicit\'e est claire puisque, quand $D$ parcourt $ \tilde{Int}_{\lambda_{1},\lambda_{2}}(\lambda)$, les $J(D)$ sont disjoints. Pour $j=\infty$, on a $j=j_{max}(D_{min})$. Soit $j\in {\cal J}$ diff\'erent de $\infty$. Supposons par exemple $j$ pair, le cas $j$ impair \'etant similaire. La d\'efinition de ${\cal J}$ et cette hypoth\`ese de parit\'e impose qu'il existe $d=1,2$ et $\Delta_{d}\in Int(\lambda_{d})$ de sorte que $j=j_{max}(\Delta_{d})$. Pour fixer la notation, on suppose qu'il en est ainsi pour $d=1$. L'ensemble des $j'\in {\cal J}$ tels que $j'<j$ n'est pas vide: il contient $j_{min}(\Delta_{1})$. Notons $j^-$ le plus grand de ces \'el\'ements. On a donc $j_{min}(\Delta_{1})\leq j^-$ et $\{j^-,...,j\}$ est contenu dans $J(\Delta_{1})$. Si $\{j^-,...,j\}$ n'est contenu dans $J(\Delta_{2})$ pour aucun $\Delta_{2}\in \tilde{Int}(\lambda_{2})$, il existe par d\'efinition des intervalles relatifs un tel intervalle $D$ tel que $J(D)=\{j^-,...,j\}$ et on a $j=j_{max}(D)$. Supposons qu'il existe un $\Delta_{2}\in \tilde{Int}(\lambda_{2})$ de sorte que $\{j^-,...,j\}\subset J(\Delta_{2})$. Si $j=j_{max}(\Delta_{2})$, alors, par d\'efinition des intervalles relatifs, il existe un tel intervalle $D$ tel que $\{j\}=J(D)$ et on conclut. Supposons $j<j_{max}(\Delta_{2})$. On note $j^+$ le plus petit \'el\'ement de ${\cal J}$ qui soit strictement sup\'erieur \`a $j$. Comme pr\'ec\'edemment, on a $j^+\leq j_{max}(\Delta_{2})$, d'o\`u $\{j,...,j^+\}\subset J(\Delta_{2})$. S'il existait $\Delta'_{1}\in \tilde{Int}(\lambda_{1})$ v\'erifiant $\{j,...,j^+\}\subset J(\Delta'_{1})$, on aurait $\Delta'_{1}=\Delta_{1}$ puisque $j\in J(\Delta_{1})$ et aussi $j_{max}(\Delta'_{1})\geq j^+>j$. Cela contredit l'hypoth\`ese $j=j_{max}(\Delta_{1})$. Un tel $\Delta'_{1}$ n'existe donc pas et, par d\'efinition des intervalles relatifs, il existe  un tel intervalle $D$ tel que $J(D)=\{j,...,j^+\}$. Alors  $j=j_{min}(D)$. $\square$ 
 
 On d\'efinit une fonction $\chi_{\lambda_{1},\lambda_{2}}:\tilde{Int}_{\lambda_{1},\lambda_{2}}(\lambda)\to {\mathbb Z}/2{\mathbb Z}$ de la fa\c{c}on suivante. Soit $D\in \tilde{Int}_{\lambda_{1},\lambda_{2}}(\lambda)$. Si $\vert J(D)\vert =1$, $\chi_{\lambda_{1},\lambda_{2}}(D)=0$. Si $\vert J(D)\vert \geq2$, $J(D)$ est de la forme (2) ci-dessus et cette relation nous fournit un indice $d\in \{1,2\}$. On note  $\chi_{\lambda_{1},\lambda_{2}}$ l'image de $d$ dans ${\mathbb Z}/2{\mathbb Z}$.  
 Remarquons que l'on a $\chi_{\lambda_{1},\lambda_{2}}(D_{min})=0$.

 On d\'efinit l'ensemble
$P^+_{\lambda_{1},\lambda_{2}}(\lambda)$ form\'e des $j_{min}(D)$ qui sont impairs, pour $D\in \tilde{Int}_{\lambda_{1},\lambda_{2}}$ et l'ensemble $P^-_{\lambda_{1},\lambda_{2}}$ form\'es des $j_{max}(D)$ qui sont  pairs, pour $D\in Int_{\lambda_{1},\lambda_{2}}(\lambda)$.  On d\'efinit une suite $\zeta_{\lambda_{1},\lambda_{2}}(\lambda)=(\zeta_{\lambda_{1},\lambda_{2}}(\lambda)_{1},\zeta_{\lambda_{1},\lambda_{2}}(\lambda)_{2},...)$ par $\zeta_{\lambda_{1},\lambda_{2}}(\lambda)_{j}=1$ si $j\in P^+_{\lambda_{1},\lambda_{2}}(\lambda)$, $\zeta_{\lambda_{1},\lambda_{2}}(\lambda)_{j}=-1$ si $j\in P^-_{\lambda_{1},\lambda_{2}}(\lambda)$, $\zeta_{\lambda_{1},\lambda_{2}}(\lambda)_{j}=0$ si $j\not\in P^+_{\lambda_{1},\lambda_{2}}(\lambda)\cup P^-_{\lambda_{1},\lambda_{2}}(\lambda)$.  

\ass{Lemme}{ On a l'\'egalit\'e $\zeta(\lambda_{1})+\zeta(\lambda_{2})=\zeta_{\lambda_{1},\lambda_{2}}(\lambda)+\xi$.}

Preuve. Restreignons-nous d'abord \`a l'ensemble des $j\geq1$ impairs. Alors les fonctions ci-dessus sont les fonctions caract\'eristiques des ensembles $P^+(\lambda_{1})$, $P^+(\lambda_{2})$, $P_{\lambda_{1},\lambda_{2}}(\lambda)$ et $J^+$. Il s'agit donc de prouver les \'egalit\'es

(4) $P^+(\lambda_{1})\cup P^+(\lambda_{2})=P_{\lambda_{1},\lambda_{2}}^+(\lambda)\cup J^+$;

(5) $P^+(\lambda_{1})\cap P^+(\lambda_{2})=P_{\lambda_{1},\lambda_{2}}^+(\lambda)\cap J^+$.

Rappelons que, puisque $\lambda_{d}$ est sp\'eciale pour $d=1,2$, $P^+(\lambda_{d})$ est l'ensemble des $j_{min}(\Delta_{d})$ pour $\Delta_{d}\in Int(\lambda_{d})$.   Consid\'erons un $j$ appartenant \`a l'ensemble de gauche de (4). Pour fixer la notation, supposons $j\in P^+(\lambda_{1})$. Alors $j=j_{min}( \Delta_{1})$ pour un $ \Delta_{1}\in Int(\lambda_{1})$, en particulier $j$ appartient \`a l'ensemble ${\cal J}$. Si  $\lambda_{2,j}$ est  impair, $j$ appartient \`a $ J^+$ par d\'efinition de cet ensemble. Supposons  $\lambda_{2,j}$  pair. Soit $j^+$ le plus petit \'el\'ement du sous-ensemble des \'el\'ements de l'ensemble ${\cal J}$ qui sont strictement sup\'erieurs \`a $j$. Ce sous-ensemble  contenant $j_{max}(\Delta_{1})$ (ou il convient ici de consid\'erer que $j_{max}(\Delta_{1,min})= \infty$), $j^+$ existe et on a $j^+\leq j_{max}(\Delta_{1})$. L'ensemble $\{j,...,j^+\}$ est contenu dans $J(\Delta_{1})$ mais, puisque  $\lambda_{2,j}$ est de mauvaise parit\'e, il n'existe pas de $\Delta_{2}\in Int(\lambda_{2})$ tel que $\{j,...,j^+\}$ soit contenu dans $J(\Delta_{2})$. Par d\'efinition $\{j,...,j^+\}$ est alors \'egal \`a $J(D)$ pour un intervalle relatif $D$ et on a $j=j_{min}(D)$. Donc $j\in P_{\lambda_{1},\lambda_{2}}^+(\lambda)$. Inversement, consid\'erons un $j$ qui appartient \`a l'ensemble de droite de (4). Si $j\in J^+$, il est par d\'efinition de la forme $j_{min}(\Delta_{d})$ pour un $d=1,2$ et un $\Delta_{d}\in Int(\lambda_{d})$. C'est-\`a-dire $j\in P^+(\lambda_{d})$. Supposons $j\in P_{\lambda_{1},\lambda_{2}}^+(\lambda)$. Alors $j=j_{min}(D)$ pour un $D\in \tilde{Int}_{\lambda_{1},\lambda_{2}}(\lambda)$. Par d\'efinition des intervalles relatifs, $j$ appartient \`a ${\cal J}$. Puisque $j$ est impair, $j$ est forc\'ement de la forme $j_{min}(\Delta_{d})$ pour un $d=1,2$ et un $\Delta_{d}\in \tilde{Int}(\lambda_{d})$. C'est-\`a-dire $j\in P^+(\lambda_{d})$. Cela prouve (4).

Soit $j\in P^+(\lambda_{1})\cap P^+(\lambda_{2})$. Alors, pour $d=1,2$, $j$ est de  la forme $j_{min}(\Delta_{d})$ pour un $\Delta_{d}\in Int(\lambda_{d})$. En particulier, $\lambda_{d,j}$ est de la  bonne parit\'e. Par d\'efinition de $J^+$, on a $j\in J^+$. Cela implique que  $\lambda_{j}$ est pair. Donc il existe un intervalle relatif $D\in \tilde{Int}_{\lambda_{1},\lambda_{2}}(\lambda)$ tel que $j\in J(D)$.   Si $j=1$, on a forc\'ement $j=j_{min}(D)$ et $j\in P_{\lambda_{1},\lambda_{2}}^+(\lambda)$. Supposons $j\geq2$. Pour $d=1,2$, l'hypoth\`ese $j=j_{min}(\Delta_{d})$ implique que $\{j-1,j\}$ n'est contenu dans  $J(\Delta'_{d})$ pour aucun $\Delta'_{d}\in Int(\lambda_{d})$. Par d\'efinition des intervalles relatifs, $\{j-1,j\}$ n'est donc contenu
dans $J(D')$ pour aucun $D'\in \tilde{Int}_{\lambda_{1},\lambda_{2}}(\lambda)$. En particulier $j-1\not\in J(D)$, d'o\`u $j=j_{min}(D)$ et $j\in P_{\lambda_{1},\lambda_{2}}^+(\lambda)$. Inversement, soit $j\in 
P_{\lambda_{1},\lambda_{2}}^+(\lambda)\cap J^+$. Par d\'efinition de $J^+$, $\lambda_{1,j}$ et $\lambda_{2,j}$ sont de bonne parit\'e et il existe $d=1,2$ et $\Delta_{d}\in Int(\lambda_{d})$ de sorte que $j=j_{min}(\Delta_{d})$. Pour fixer la notation, on suppose que ce $d$ est \'egal \`a $1$. Donc $j\in P^+(\lambda_{1})$.  L'hypoth\`ese que $\lambda_{2,j}$ est de bonne parit\'e implique qu'il existe $\Delta_{2}\in Int(\lambda_{2})$ de sorte que $j\in J(\Delta_{2})$. Supposons d'abord que tous les \'el\'ements de ${\cal J}$ soient sup\'erieurs ou \'egaux \`a $j$. Dans ce cas, $j=j_{min}(\Delta_{2})$ et $j\in P^+(\lambda_{2})$. Supposons maintenant qu'il existe des \'el\'ements de ${\cal J}$ strictement inf\'erieurs \`a $j$, notons $j^-$ le plus grand d'entre eux.   L'hypoth\`ese $j\in P_{\lambda_{1},\lambda_{2}}^+(\lambda)$ signifie que $j=j_{min}(D)$ pour un intervalle relatif $D$. Donc 
$\{j^-,...,j\}$ n'est  de la forme $J(D')$ pour aucun $D'\in \tilde{Int}_{\lambda_{1},\lambda_{2}}(\lambda)$. Les entiers $j^-$ et $j$ sont deux \'el\'ements cons\'ecutifs de ${\cal J}$. Ces deux propri\'et\'es et la d\'efinition des intervalles relatifs entra\^{\i}nent que le nombre de $d$ pour lesquels il existe $\Delta'_{d}\in Int(\lambda_{d})$ tel que $\{j^-,...,j\}\subset J(\Delta'_{d})$ est pair. Pour $d=1$, il n'existe pas de tel $\Delta'_{1}$ car $j=j_{min}(\Delta_{1})$. Donc il n'existe pas non plus de tel $\Delta'_{2}$. En particulier $\{j^-,...,j\}\not\subset J(\Delta_{2})$. Puisque $\{j_{min}(\Delta_{2}),...,j\}\subset J(\Delta_{2})$, cela entra\^{\i}ne $j^-<j_{min}(\Delta_{2})$, et, puisque $j_{min}(\Delta_{2})\in {\cal J}$, la d\'efinition de $j^-$ entra\^{\i}ne $j\leq j_{min}(\Delta_{2})$, d'o\`u forc\'ement $j=j_{min}(\Delta_{2})$. Donc $j\in P^+(\lambda_{2})$. Cela prouve (5).

Un raisonnement analogue vaut en se restreignant \`a l'ensemble des entiers pairs  $j\geq2$. Cela prouve  le lemme. $\square$

On dit que $\lambda_{1}$ et $\lambda_{2}$ induisent r\'eguli\`erement $\lambda$ si et seulement si $\tilde{Int}_{\lambda_{1},\lambda_{2}}(\lambda)$ est la partition la plus fine de $Jord^{bp}(\lambda)\cup\{0\}$, c'est-\`a-dire si et seulement si tout intervalle relatif est r\'eduit \`a un seul \'el\'ement.  Dans ce cas, $\chi_{\lambda_{1},\lambda_{2}}$ est d\'efinie sur $Jord^{bp}(\lambda)\cup \{0\}$ et on a $\chi_{\lambda_{1},\lambda_{2}}(0)=0$.

\bigskip

\subsection{Une proposition d'existence}
Soient $n\in {\mathbb N}$ et $\lambda\in {\cal P}^{symp}(2n)$. Fixons une fonction $\chi:Jord^{bp}(\lambda)\cup\{0\}\to {\mathbb Z}/2{\mathbb Z}$ telle que $\chi(i)=0$ pour tout $i\in Jord^{bp}(\lambda)$ tel que $mult_{\lambda}(i)=1$ et telle que $\chi(0)=0$.

\ass{Proposition}{Il existe $n_{1},n_{2}\in {\mathbb N}$ tels que $n_{1}+n_{2}=n$ et il existe $\lambda_{1}\in {\cal P}^{symp,sp}(2n_{1})$ et $\lambda_{2}\in {\cal P}^{orth,sp}(2n_{2})$ tels que

(i) $\lambda_{1}$ et $\lambda_{2}$ induisent r\'eguli\`erement $\lambda$;

(ii) $d(\lambda_{1})\cup d(\lambda_{2})=d(\lambda)$;

(iii) $\chi_{\lambda_{1},\lambda_{2}}=\chi$.}

La preuve est identique \`a celle de \cite{W4} 1.11. On la refait car, dans cette r\'ef\'erence, on avait b\^etement suppos\'e que tous les termes de $\lambda$ \'etaient pairs. On utilise les notations de 2.7. 

Preuve. Notons $\mathfrak{ J}^+$ l'ensemble des $j\geq1$ tels que $j$ soit impair, $\lambda_{j}$ soit pair  et $\lambda_{j}>\lambda_{j+1}$. Notons $\mathfrak{J}^-$ l'ensemble des $j\geq2$ tels que $j$ et $\lambda_{j}$ soient pairs et $\lambda_{j-1}> \lambda_{j}$. On voit que $\mathfrak{ J}^+$ est l'ensemble des $j_{max}(i)$ pour $i=i_{h}$ avec $h$ impair ou pour $i\in {\cal J}''(\lambda)\cap Jord^{bp}(\lambda)$. De m\^eme,
 $\mathfrak{J}^-$ est l'ensemble des $j_{min}(i)$ pour  $i=i_{h}$ avec $h$ pair ou pour $i\in {\cal J}''(\lambda)\cap Jord^{bp}(\lambda)$. On en d\'eduit que $\mathfrak{ J}^+$ et $\mathfrak{J}^-$ ont le m\^eme nombre d'\'el\'ements et que, si on note $\mathfrak{J}^+=\{j_{1}^+<...<j_{c}^+\}$ et $\mathfrak{J}^-=\{j_{1}^-<...<j^-_{c}\}$, on a
$$j_{1}^+<j_{1}^-<j_{2}^+<j_{2}^-<....<j^+_{c}<j^-_{c}.$$
On note $\mathfrak{r}=(\mathfrak{r}_{1},\mathfrak{r}_{2},...)$ la suite de nombres d\'efinie par $\mathfrak{r}_{j}=1$ si $j\in \mathfrak{J}^+$, $\mathfrak{r}_{j}=-1$ si $j\in \mathfrak{J}^-$ et $\mathfrak{r}_{j}=0$ si $j\not\in \mathfrak{J}^+\cup\mathfrak{J}^-$. 

 Soit $d\in \{1,2\}$. Pour $j\geq1$, disons que $j$ et $j+1$ sont $d$-li\'es si et seulement si l'une des conditions suivantes est v\'erifi\'ee:

(1)(a) $\lambda_{j}=\lambda_{j+1}$ est pair et $\chi(\lambda_{j})=d+1$ (c'est-\`a-dire $\chi(\lambda_{j})\equiv d+1\,\,mod\,\,2{\mathbb Z}$);

(1)(b) $j\in \mathfrak{J}^+$;

(1)(c) $j+1\in \mathfrak{J}^-$;

(1)(d) $\lambda_{j}$ et $\lambda_{j+1}$ sont impairs et  $\lambda_{j}\in {\cal J}''$.

Remarquons que cette derni\`ere condition \'equivaut \`a

(1)(d') $\lambda_{j}$ et $\lambda_{j+1}$ sont impairs et  $\lambda_{j+1}\in {\cal J}''$.

En effet, si (1)(d) est v\'erifi\'ee, on a $i_{h}>\lambda_{j}>i_{h+1}$ pour un $h$ impair. Alors $i_{h}>\lambda_{j+1}\geq i_{h+1}$. Mais $\lambda_{j+1}\not=i_{h+1}$ puisque $\lambda_{j+1}$ est impair et $i_{h+1}$ est pair. Donc $i_{h}>\lambda_{j+1}>i_{h+1}$ et $\lambda_{j+1}\in {\cal J}''$. La r\'eciproque est similaire.

Pour deux entiers $1\leq j\leq j'$, disons qu'ils sont $d$-li\'es si et seulement si $k$ et $k+1$ sont $d$-li\'es pour tout $k=j,...,j'-1$. C'est une relation d'\'equivalence et les classes sont des intervalles de ${\mathbb N}-\{0\}$, \'eventuellement infinis. On  note $\tilde{\mathfrak{Int}}_{d}$ l'ensemble des classes d'\'equivalence ayant au moins deux \'el\'ements. Pour $\mathfrak{I}\in\tilde{ \mathfrak{Int}}_{d}$, on note $j_{min}(\mathfrak{I})$, resp. $j_{max}(\mathfrak{I})$, le plus petit, resp. grand, \'el\'ement de $\mathfrak{I}$ (avec $j_{max}(\mathfrak{I})=\infty$ si $\mathfrak{I}$ est infini). Pour $d=1,2$ d\'efinissons une fonction $p_{d}:{\mathbb N}-\{0\}\to {\mathbb Z}/2{\mathbb Z}$ par $p_{d}(j)=1$ s'il existe $\mathfrak{I}\in \tilde{\mathfrak{Int}}_{d}$ tel que $j\in \mathfrak{I}$, $p_{d}(j)=0$ sinon. Montrons que

(2) l'ensemble $\tilde{\mathfrak{Int}}_{d}$ est fini; il contient un \'el\'ement infini si et seulement si $d=1$;

\noindent on  note $\mathfrak{Int}_{1}$ l'ensemble $\tilde{\mathfrak{Int}}_{1}$ priv\'e de cet \'el\'ement infini et on pose $\mathfrak{Int}_{2}=\tilde{\mathfrak{Int}}_{2}$;

(3) pour $\mathfrak{I}\in \tilde{\mathfrak{Int}}_{d}$, $j_{min}(\mathfrak{I})$ est impair et $j_{max}(\mathfrak{I})$ est pair ou infini;

(4) pour $j\geq1$, on a
$$p_{1}(j)+p_{2}(j)=\left\lbrace\begin{array}{cc}
2,&{\text\,\,si\,\,} j\in \mathfrak{J}^+\cup \mathfrak{J}^-;\\
 1,&{\text \,\,si\,\,} \lambda_{j}{\text\,\, est\,\, pair\,\, et \,\,}j\not\in \mathfrak{J}^+\cup\mathfrak{J}^-;\\
 0,&{\text\,\,  si\,\,} \lambda_{j}{\text\,\, est \,\,impair\,\, et\,\,}  \lambda_{j}\in {\cal J}'(\lambda);\\ 
 2,&{\text\,\, si\,\,} \lambda_{j}{\text\,\, est\,\, impair\,\, et\,\,}  \lambda_{j}\in {\cal J}''(\lambda);\\ \end{array}\right.$$
 
 (5) $\mathfrak{J}^+$ est \'egal \`a l'ensemble des $j\geq1$ tels que $p_{1}(j)=p_{2}(j)=1$ et qu'il existe $d=1,2$ et un \'el\'ement de $\mathfrak{I}\in\mathfrak{Int}_{d}$ de sorte que $j=j_{min}(\mathfrak{I})$;
 
 (6) $\mathfrak{J}^-$ est \'egal \`a l'ensemble des $j\geq1$ tels que $p_{1}(j)=p_{2}(j)=1$ et qu'il existe $d=1,2$ et un \'el\'ement de $\mathfrak{I}\in\mathfrak{Int}_{d}$ de sorte que $j=j_{max}(\mathfrak{I})$.

Soit $l(\lambda)$ le plus grand entier $l$ tel que $\lambda_{l}>0$. Parce que $\chi(0)=0$, on voit que, pour $j>l(\lambda)$, $j$ et $j+1$ sont $1$-li\'es mais pas $2$-li\'es. Donc $\{l(\lambda)+1,...\}$ est contenu dans un intervalle infini $\mathfrak{I}_{1,min}\in \tilde{\mathfrak{Int}}_{1}$ tandis que, pour $j\geq l(\lambda)+2$, $\{j\}$ est une classe d'\'equivalence pour la $2$-liaison et $j$ n'est pas contenu dans un \'el\'ement de $\mathfrak{Int}_{2}$. Cela prouve (2).

Soit $\mathfrak{I}\in \tilde{\mathfrak{Int}}_{d}$. On pose simplement $j=j_{min}(\mathfrak{I})$. Montrons que $j $ est impair.   C'est \'evident si $j=1$. On suppose $j\geq2$. Par d\'efinition, $j$ et $j+1$ sont $d$-li\'es tandis que $j-1$ et $j$ ne le sont pas. Si (1)(b) ou (1)(c) est v\'erifi\'ee, $j$ est trivialement impair. 
Supposons v\'erifi\'ee (1)(a). On n'a pas $\lambda_{j-1}=\lambda_{j}$: sinon ces entiers seraient pairs, on aurait $\chi(\lambda_{j-1})=\chi(\lambda_{j})=d+1$ et $j-1$ et $j$  v\'erifieraient l'analogue de (1)(a) et seraient $d$-li\'es. Donc $\lambda_{j-1}>\lambda_{j}$. Alors $j$ est impair ou appartient \`a $\mathfrak{J}^-$. Or cette derni\`ere relation est exclue car elle entra\^{\i}ne que $j-1$ et $j$ v\'erifient l'analogue de (1)(c) et  sont $d$-li\'es. Donc $j$ est impair. Supposons maintenant que (1)(d)  soit v\'erifi\'ee. 
Supposons d'abord que $\lambda_{j-1}$ est impair.  Alors $j-1$ et $j$ v\'erifient l'analogue de (1)(d') et sont $d$-li\'es, ce qui n'est pas le cas. Donc $\lambda_{j-1}$ est pair et $\lambda_{j-1}>\lambda_{j}$. Alors $j-1$ est pair ou $j-1\in \mathfrak{J}^+$. Or cette derni\`ere relation est exclue car elle entra\^{\i}ne que $j-1$ et $j$ v\'erifient l'analogue de (1)(b) et sont $d$-li\'es. donc $j-1$ est pair et $j$ est impair. Cela montre que $j_{min}(\mathfrak{I})$ est impair. Une preuve analogue montre que $j_{max}(\mathfrak{I})$ est pair s'il n'est pas infini. Cela prouve (3).

Soit $j\in \mathfrak{J}^+$. Alors (1)(b) est v\'erifi\'e et $j$ et $j+1$ sont $d$-li\'es pour $d=1,2$. Donc  $p_{1}(j)=p_{2}(j)=1$. Soit maintenant $j\in \mathfrak{J}^-$. Alors $j$ est pair donc diff\'erent de $1$. L'analogue de (1)(c) pour le couple $(j-1,j)$ est v\'erifi\'ee et $j-1$ et $j$ sont $d$-li\'es pour $d=1,2$. Donc $p_{1}(j)=p_{2}(j)=1$. Supposons maintenant $\lambda_{j}$ pair mais $j\not\in \mathfrak{J}^+\cup \mathfrak{J}^-$. Supposons par exemple $j$ impair, le cas o\`u $j$ est pair se traitant de fa\c{c}on analogue. Puisque $j\not\in \mathfrak{J}^+$, on a $\lambda_{j}=\lambda_{j+1}$. Les entiers $j$ et $j+1$ sont $d$-li\'es pour l'unique $d$ tel que $\chi(\lambda_{j})=d+1$. Pour ce $d$,  $p_{d}(j)=1$. Soit $d'$ l'autre \'el\'ement de $\{1,2\}$. On doit prouver que $j$ n'appartient \`a aucun \'el\'ement de $\tilde{\mathfrak{Int}}_{d'}$. On vient de voir que $j$ et $j+1$ ne sont pas $d'$-li\'es.  Si $j$ appartenait \`a un \'el\'ement $\mathfrak{I}'\in \tilde{\mathfrak{Int}}_{d'}$, cet intervalle serait fini et $j$ serait \'egal \`a  $j_{max}(\mathfrak{I}')$. Mais alors $j$ serait pair d'apr\`es (3), contrairement \`a l'hypoth\`ese. Supposons maintenant  $\lambda_{j}$ impair, $j$ impair  et $\lambda_{j}\in {\cal J}'(\lambda)$. Cette derni\`ere condition implique  d'apr\`es 2.7 que $j_{max}(\lambda_{j})$ est pair, donc $j<j_{max}(\lambda_{j})$, donc $\lambda_{j+1}=\lambda_{j}$.
Pour $d=1,2$, les conditions (1)(a), (1)(b)  et (1)(c) ne sont pas v\'erifi\'ees: elles imposent que $\lambda_{j}$ ou $\lambda_{j+1}$ est pair. La condition (1)(d) ne l'est pas puisque  $\lambda_{j}\in {\cal J}'(\lambda)$. Donc $j$ et $j+1$ ne sont pas $d$-li\'es. Si $j=1$, $j$ n'appartient donc \`a aucun \'el\'ement de $\tilde{\mathfrak{Int}}_{d}$. Si $j>1$, les analogues des conditions (1)(a) et (1)(c) pour le couple $(j-1,j)$ ne sont pas v\'erifi\'ees: elles imposent que $\lambda_{j}$ est pair. L'analogue de (1)(b) n'est pas v\'erifi\'ee: elle impose $j-1$ impair donc $j$ pair.  L'analogue de (1)(d') n'est pas v\'erifi\'ee puisque $\lambda_{j}\in {\cal J}'(\lambda)$. Donc $j-1$ et $j$ ne sont pas $d$-li\'es. Donc  $p_{d}(j)=0$. Supposons maintenant  $\lambda_{j}$ impair, $j$ pair  et  $j\in {\cal J}'(\lambda)$.  Cette derni\`ere condition implique  d'apr\`es 2.7 que $j_{min}(\lambda_{j})$ est impair, donc $j_{min}(\lambda_{j})<j$, donc $\lambda_{j-1}=\lambda_{j}$. Des arguments analogues \`a ceux ci-dessus montrent que, pour $d=1,2$,   $p_{d}(j)=0$. Supposons  enfin que  $\lambda_{j}$ est  impair et que    $j\in {\cal J}''(\lambda)$.   Puisque $mult_{\lambda}(\lambda_{j})$ est paire, on a $\lambda_{j-1}=\lambda_{j}$ ou $\lambda_{j+1}=\lambda_{j}$. Dans le premier cas, $j-1$ et $j$ v\'erifient l'analogue de (1)(d') et  sont $d$-li\'es pour tout $d$. Dans le deuxi\`eme cas, $j$ et $j+1$ v\'erifient (1)(d) et sont $d$-li\'es pour tout $d$. Donc $p_{d}(j)=1$ pour tout $d$. Cela d\'emontre (4).  

 Soit $j\in \mathfrak{J}^+$. D'apr\`es (4), on a $p_{1}(j)=p_{2}(j)=1$, c'est-\`a-dire que, pour tout $d$, il existe $\mathfrak{I}_{d}\in\tilde{\mathfrak{Int}}_{d}$ tel que $j\in \mathfrak{I}_{d}$. Si $j=1$, on a forc\'ement $j=j_{min}(\mathfrak{I}_{d})$ pour tout $d$. Supposons $j>1$. On veut montrer que $j=j_{min}(\mathfrak{I}_{d})$ pour au moins un $d$, autrement dit que $j-1$ et $j$ ne sont pas $d$-li\'es pour au moins un $d$. Les analogues pour de couple $(j-1,j)$ des conditions (1)(b) et (1)(c) ne sont pas v\'erifi\'ees: elles impliquent que $j$ est pair, alors que $j$ est impair puisque $j\in \mathfrak{J}^+$. L'analogue de (1)(d) n'est pas v\'erifi\'ee, puisque $\lambda_{j}$ est pair. Donc $j-1$ et $j$ ne sont $d$-li\'es que si l'analogue de (1)(a) est v\'erifi\'ee. Mais cette analogue ne peut \^etre v\'erifi\'ee que pour un unique $d$. Cela d\'emontre que $\mathfrak{J}^+$ est contenu dans l'ensemble d\'ecrit en (5). Inversement, soit $j\geq1$, supposons que $p_{1}(j)=p_{2}(j)=1$ et qu'il existe $d=1,2$ et un \'el\'ement de $\mathfrak{I}\in\tilde{\mathfrak{Int}}_{d}$ de sorte que $j=j_{min}(\mathfrak{I})$.   Autrement dit, ou bien $j=1$, ou bien il existe $d$ tel que $j-1$ et $j$ ne sont pas $d$-li\'es. D'apr\`es (3), $j$ est impair. D'apr\`es (4), on a soit $j\in \mathfrak{J}^+\cup \mathfrak{J}^-$, soit $\lambda_{j}$ est impair et  $\lambda_{j}\in {\cal J}''(\lambda)$. Dans le premier cas, l'imparit\'e de $j$ entra\^{\i}ne $j\in \mathfrak{J}^+$, ce que l'on veut prouver. Supposons donc que $\lambda_{j}$  est  impair et $\lambda_{j}\in {\cal J}''(\lambda)$. D'apr\`es 2.7, cette condition entra\^{\i}ne que $j_{min}(\lambda_{j})$ est pair, donc $j_{min}(\lambda_{j})<j$ et $\lambda_{j-1}=\lambda_{j}$. Alors $j-1$ et $j$ v\'erifient l'analogue de (1)(d') et sont $d$-li\'es. Cela contredit l'hypoth\`ese. On a ainsi prouv\'e (5). La preuve de (6) est similaire.

  La relation (3) entra\^{\i}ne

(7) $p_{d}(j)=p_{d}(j+1)$ si $j$ est impair.

La d\'efinition de $\mathfrak{r}$ et l'assertion (4) entra\^{\i}nent

(8) $\mathfrak{r}_{j}\equiv p_{1}(j)+p_{2}(j)+1+\lambda_{j}\,\,mod\,\,2{\mathbb Z}$.

On va montrer qu'il existe des suites d'entiers positifs ou nuls $\lambda_{1}$ et $\lambda_{2}$ v\'erifiant les conditions suivantes, pour tout $j\geq1$:

(9) $\lambda_{1,j}+\lambda_{2,j}+\mathfrak{r}_{j}=\lambda_{j}$;

(10) pour $d=1,2$, $\lambda_{d,j}\equiv d+p_{d}(j)\,\,mod\,\,2{\mathbb Z}$;

(11) pour $d=1,2$, on a

(a) $\lambda_{d,j}=\lambda_{d,j+1}$ si $j$ est pair, $p_{d}(j)=1$ et il n'existe pas de $\mathfrak{I}\in \mathfrak{Int}_{d}$ tel que $j=j_{max}(\mathfrak{I})$ ou si $j$ est impair et $p_{d}(j)=0$;

(b) $\lambda_{d,j}>\lambda_{d,j+1}$ si $j$ est pair et il existe $\mathfrak{I}\in \mathfrak{Int}_{d}$ tel que $j=j_{max}(\mathfrak{I})$ (la condition que $j$ est pair est redondante d'apr\`es (3));

(c) $\lambda_{d,j}\geq \lambda_{d,j+1}$ si $j$ est impair et $p_{d}(j)=1$ ou si $j$ est pair et $p_{d}(j)=0$. 

On raisonne par r\'ecurrence descendante sur $j$. Pour $j\geq l(\lambda)+2$, on pose $\lambda_{1,j}=\lambda_{2,j}=0$. On a vu dans la preuve de (2) que $j$ \'etait contenu dans $\mathfrak{I}_{1,min}$ mais dans aucun \'el\'ement de $\mathfrak{Int}_{2}$. Donc $p_{1}(j)=1$ et $p_{2}(j)=0$. De plus, $j$ n'appartient pas \`a $\mathfrak{J}^+\cup\mathfrak{J}^-$ donc $\mathfrak{r}_{j}=0$. On voit alors que toutes les conditions ci-dessus sont v\'erifi\'ees.

On fixe $j$ et on suppose que l'on a fix\'e des termes $\lambda_{1,j'}$, $\lambda_{2,j'}$ pour $j'>j$ de sorte que les conditions ci-dessus soient v\'erifi\'ees pour ces $j'$. Pour $d=1,2$, on pose $\lambda_{d,j}=\lambda_{d,j+1}+e_{d}$, avec $e_{d}\in {\mathbb Z}$. Les conditions ci-dessus se traduisent en termes de ces entiers $e_{d}$. L'analogue de (9) \'etant v\'erifi\'ee pour $j+1$, cette condition (9) se traduit par

(12) $e_{1}+e_{2}=\lambda_{j}-\lambda_{j+1}+\mathfrak{r}_{j+1}-\mathfrak{r}_{j}$.

De m\^eme, la condition (10) se traduit par

(13) $e_{d}\equiv p_{d}(j)+p_{d}(j+1)\,\,mod\,\,2{\mathbb Z}$.

Remarquons que, si (12) est v\'erifi\'ee, la relation (8) entra\^{\i}ne

$e_{1}+e_{2}\equiv p_{1}(j)+p_{1}(j+1)+p_{2}(j)+p_{2}(j+1)\,\,mod\,\,2{\mathbb Z}$.

\noindent Donc (13) est v\'erifi\'ee pour un $d$ si et seulement si elle l'est pour les deux $d$. 

La condition (11) se traduit par $e_{d}=0$ dans le cas (a), $e_{d}>0$ dans le cas (b) et $e_{d}\geq0 $ dans le cas (c). Remarquons que, dans le cas (a), la condition $e_{d}=0$ est compatible avec (13), autrement dit $p_{d}(j)=p_{d}(j+1)$.  En effet, si $j$ est impair,  cette relation est toujours vraie d'apr\`es (5). Si $j$ est pair, la condition (11)(a) impose que $j$ et $j+1$ sont $d$-li\'es donc $p_{d}(j)=p_{d}(j+1)=1$. 

Supposons la condition (11)(a) v\'erifi\'ee pour un $d$, disons pour $d=1$. On n'a pas le choix pour $e_{1}$: on pose $e_{1}=0$.  La condition (12) impose $e_{2}=\lambda_{j}-\lambda_{j+1}+\mathfrak{r}_{j+1}-\mathfrak{r}_{j}$. Comme on vient de le dire, la condition (13) est v\'erifi\'ee pour $d=1$. Elle l'est donc aussi pour $d=2$. Il reste \`a v\'erifier les conditions provenant de (11) pour $d=2$. 

Supposons $j$ impair. Supposons d'abord que la condition (11)(a) soit v\'erifi\'ee pour $d=2$, auquel cas on doit v\'erifier que $e_{2}=0$. La condition (11)(a)  pour $j$ impair est que $p_{d}(j)=0$. Cette condition est v\'erifi\'ee pour $d=1,2$. D'apr\`es (4), $\lambda_{j}$ est impair et  $\lambda_{j}\in {\cal J}'(\lambda)$. D'apr\`es 2.7, $j_{max}(\lambda_{j})$ est pair, donc $j<j_{max}(\lambda_{j})$ et $\lambda_{j}=\lambda_{j+1}$.  Evidemment, $j,j+1\not\in \mathfrak{J}^+\cup\mathfrak{J}^-$, donc $\mathfrak{r}_{j}=\mathfrak{r}_{j+1}=0$. Alors $e_{2}=\lambda_{j}-\lambda_{j+1}+\mathfrak{r}_{j+1}-\mathfrak{r}_{j}=0$. 
La condition (11)(b) n'est pas v\'erifi\'ee pour $d=2$ puisque $j$ est impair. Supposons la condition (11)(c) v\'erifi\'ee pour $d=2$. On doit alors prouver que $e_{2}\geq0$. Puisque $j$ est impair, cette condition est que $p_{2}(j)=1$. On a aussi $p_{1}(j)=0$ puisque (11)(a) est v\'erifi\'ee pour $d=1$. D'apr\`es (7), on a aussi $p_{1}(j+1)=0$ et $p_{2}(j+1)=1$. Alors, d'apr\`es (4),  ni $j$, ni $j+1$  n'appartiennent  \`a $ \mathfrak{J}^+\cup\mathfrak{J}^-$. Donc $\mathfrak{r}_{j}=\mathfrak{r}_{j+1}=0$.  Donc $e_{2}=\lambda_{j}-\lambda_{j+1} \geq0$.  

Supposons plut\^ot $j$ pair. Supposons d'abord que la condition (11)(a) soit v\'erifi\'ee pour $d=2$, auquel cas on doit v\'erifier que $e_{2}=0$. Pour $j$ pair, la condition (11)(a) pour $d$  est que $p_{d}(j)=1$ et qu'il n'existe pas de $\mathfrak{I}\in \mathfrak{Int}_{d}$ tel que $j=j_{max}(\mathfrak{I})$. Cette condition est v\'erifi\'ee pour $d=1,2$. D'apr\`es (4), on a soit $j\in \mathfrak{J}^+\cup \mathfrak{J}^-$, soit $\lambda_{j}$  est impair et $\lambda_{j}\in {\cal J}''(\lambda)$. Dans le premier cas, la parit\'e de $j$ impose $j\in \mathfrak{J}^-$.  Mais alors la relation (6)  implique l'existence de $d$ et de $\mathfrak{I}\in \mathfrak{Int}_{d}$ tels que $j=j_{max}(\mathfrak{I})$, contrairement aux hypoth\`eses. Supposons donc que $\lambda_{j}$  soit impair et que $\lambda_{j}\in {\cal J}''(\lambda)$. D'apr\`es 2.7, $j_{max}(\lambda_{j})$ est impair, donc $j<j_{max}(\lambda_{j})$ et $\lambda_{j}=\lambda_{j+1}$. Evidemment, $j,j+1\not\in \mathfrak{J}^+\cup\mathfrak{J}^-$, donc $\mathfrak{r}_{j}=\mathfrak{r}_{j+1}=0$. Alors $e_{2}=\lambda_{j}-\lambda_{j+1}+\mathfrak{r}_{j+1}-\mathfrak{r}_{j}=0$. Supposons maintenant v\'erifi\'ee la condition (11)(b) pour $d=2$. On doit prouver que $e_{2}>0$. La condition est que  $p_{2}(j)=1$ et qu'il existe $\mathfrak{I}\in \mathfrak{Int}_{2}$ tel que $j=j_{max}(\mathfrak{I})$. On a aussi $p_{1}(j)=1$ puisque (11)(a) est v\'erifi\'ee pour $d=1$. D'apr\`es (6), on a $j\in \mathfrak{J}^-$. Cela entra\^{\i}ne $\mathfrak{r}_{j}=-1$.   Puisque $j+1$ est impair, on a $j+1\not\in \mathfrak{J}^-$ donc 
  $\mathfrak{r}_{j+1}\leq0$. Alors  $e_{2}=\lambda_{j}-\lambda_{j+1}+\mathfrak{r}_{j+1}-\mathfrak{r}_{j}\geq\lambda_{j}-\lambda_{j+1}+1>0$.  Supposons enfin v\'erifi\'ee la condition (11)(c) pour $d=2$, autrement dit $p_{2}(j)=0$. On doit v\'erifier que $e_{2}\geq0$. Puisque $p_{1}(j)=1$, on a $\lambda_{j}$ pair et $j\not\in \mathfrak{J}^+\cup \mathfrak{J}^-$ d'apr\`es (4). Le m\^eme raisonnement que dans le cas $j$ impair s'applique et on conclut $e_{2}\geq0$.
  
  On peut maintenant supposer que la condition (11)(a) n'est v\'erifi\'ee pour aucun $d$. Supposons la condition (11)(b) v\'erifi\'ee pour $d=1$. On choisit pour $e_{1}$ le plus petit entier strictement
positif v\'erifiant la relation (13). On a $e_{1}=1$ ou $2$. La condition r\'esultant de (11)(b) pour $d=1$ est $e_{1}>0$, elle est v\'erifi\'ee. On pose $e_{2}=\lambda_{j}-\lambda_{j+1}+\mathfrak{r}_{j+1}-\mathfrak{r}_{j}-e_{1}$. Comme pr\'ec\'edemment, il reste seulement \`a prouver que $e_{2}$ v\'erifie les conditions r\'esultant de (11) pour $d=2$. On a exclu la condition (11)(a). Supposons que la condition (11)(b) soit v\'erifi\'ee pour $d=2$. On doit montrer que $e_{2}>0$. Les conditions (11)(b) sont v\'erifi\'ees pour $d=1,2$, c'est-\`a-dire que $j$ est pair et qu'il existe $\mathfrak{I}_{d}\in \mathfrak{Int}_{d}$ de sorte que $j=j_{max}(\mathfrak{I}_{d})$. Autrement dit, $p_{d}(j)=1$ mais $j$ et $j+1$ ne sont pas $d$-li\'es. D'apr\`es (6), on a $j\in \mathfrak{J}^-$, donc  $\lambda_{j}$ est pair. Si $\lambda_{j+1}=\lambda_{j}$, $j$ et $j+1$ v\'erifient  (1)(a) pour un $d$ et sont $d$-li\'es contrairement \`a l'hypoth\`ese. Donc $\lambda_{j}>\lambda_{j+1}$. Puisque $j\in \mathfrak{J}^-$, on a aussi $\mathfrak{r}_{j}=-1$. Le nombre $j+1$ est impair donc n'appartient pas \`a $\mathfrak{J}^-$, d'o\`u $\mathfrak{r}_{j+1}\geq0$. On voit alors que $e_{2}=\lambda_{j}-\lambda_{j+1}+\mathfrak{r}_{j+1}-\mathfrak{r}_{j}-e_{1}$ est strictement positif sauf si les trois conditions suivantes sont v\'erifi\'ees: $\lambda_{j}=\lambda_{j+1}+1$, $\mathfrak{r}_{j+1}=0$ et $e_{1}=2$. Supposons ces conditions v\'erifi\'ees. Puisque $p_{1}(j)=1$ et $e_{1}=2$, la condition (13) pour $d=1$, qui est v\'erifi\'ee par d\'efinition de $e_{1}$, implique $p_{1}(j+1)=1$. Puisque $\lambda_{j}=\lambda_{j+1}+1$, $\lambda_{j+1}$ est impair. Puisque $\mathfrak{r}_{j+1}=0$, la relation (8) implique que $p_{2}(j+1)=1$. Alors, pour $d=1,2$, $j+1$ appartient \`a un \'el\'ement $\mathfrak{I}'_{d}\in \tilde{\mathfrak{Int}}_{d}$. Puisque $j$ et $j+1$ ne sont pas $d$-li\'es, on a forc\'ement $j+1=j_{min}(\mathfrak{I}'_{d})$. D'apr\`es (5), cela entra\^{\i}ne $j+1\in \mathfrak{J}^+$. Donc $\mathfrak{r}_{j+1}=1$ contrairement \`a l'hypoth\`ese. Cette contradiction conclut. Supposons maintenant que la condition (11)(c) soit v\'erifi\'ee pour $d=2$. On doit montrer que $e_{2}\geq0$. On a toujours la condition (11)(b) pour $d=1$, c'est-\`a-dire  que $j$ est pair, que $p_{1}(j)=1$ mais que $j$ et $j+1$ ne sont pas $1$-li\'es. La condition (11)(c) pour $d=2$ dit que $p_{2}(j)=0$. Alors $j$ et $j+1$ ne sont pas non plus $2$-li\'es. D'autre part, la relation (4) entra\^{\i}ne que $\lambda_{j}$ est pair et que $j\not\in \mathfrak{J}^+\cup \mathfrak{J}^-$. D'o\`u $\mathfrak{r}_{j}=0$. On ne peut pas avoir $\lambda_{j}=\lambda_{j+1}$ sinon la relation (1)(a) serait v\'erifi\'ee pour un $d$ et $j$ et $j+1$ seraient $d$-li\'es, ce qui n'est pas le cas. On n'a pas $j+1\in \mathfrak{J}^-$ puisque $j+1$ est impair. Donc $\mathfrak{r}_{j+1}\geq0$. On voit alors que $e_{2}=\lambda_{j}-\lambda_{j+1}+\mathfrak{r}_{j+1}-\mathfrak{r}_{j}-e_{1}$ est  positif ou nul sauf si les m\^emes conditions que ci-dessus sont v\'erifi\'ees: $\lambda_{j}=\lambda_{j+1}+1$, $\mathfrak{r}_{j+1}=0$ et $e_{1}=2$. Ces conditions sont exclues par le  m\^eme raisonnement que ci-dessus. D'o\`u $e_{2}\geq0$.

Il nous reste \`a traiter le cas o\`u (11)(c) est v\'erifi\'ee pour $d=1,2$. On choisit pour $e_{1}$ le plus petit entier 
positif  ou nul v\'erifiant la relation (13). On a $e_{1}=0$ ou $1$. La condition r\'esultant de (11)(c) pour $d=1$ est $e_{1}\geq0$, elle est v\'erifi\'ee. On pose $e_{2}=\lambda_{j}-\lambda_{j+1}+\mathfrak{r}_{j+1}-\mathfrak{r}_{j}-e_{1}$. Comme pr\'ec\'edemment, il reste seulement \`a prouver que $e_{2}$ v\'erifie la condition r\'esultant de (11)(c) pour $d=2$, c'est-\`a-dire $e_{2}\geq0$. 

Supposons d'abord $j$ impair. Les conditions (11)(c) pour $d=1,2$ disent que $p_{1}(j)=p_{2}(j)=1$. D'apr\`es (7), on a aussi $p_{1}(j+1)=p_{2}(j+1)=1$. La relation (13) pour $d=1$ implique $e_{1}=0$. Si ni $j$, ni $j+1$ n'appartiennent \`a $\mathfrak{J}^+\cup\mathfrak{J}^-$, on a $\mathfrak{r}_{j}=\mathfrak{r}_{j+1}=0$ et $e_{2}=\lambda_{j}-\lambda_{j+1}\geq0$. Si un seul des \'el\'ements $j$ et $j+1$ appartiennent \`a $\mathfrak{J}^+\cup\mathfrak{J}^-$, on a par parit\'e $j\in \mathfrak{J}^+$ et $j+1\not\in\mathfrak{J}^+\cup\mathfrak{J}^-$, ou $j+1\in \mathfrak{J}^-$ et $j\not\in\mathfrak{J}^+\cup\mathfrak{J}^-$. Alors $\mathfrak{r}_{j+1}-\mathfrak{r}_{j}=-1$. Mais l'hypoth\`ese $j\in \mathfrak{J}^+$ ou $j+1\in \mathfrak{J}^-$ implique $\lambda_{j}>\lambda_{j+1}$. Alors $e_{2}=\lambda_{j}-\lambda_{j+1}-1\geq0$. Enfin supposons que $j$ et $j+1$ appartiennent  tous deux \`a $\mathfrak{J}^+\cup\mathfrak{J}^-$. La parit\'e impose $j\in \mathfrak{J}^+$ et $j+1\in\mathfrak{J}^-$. Alors $\mathfrak{r}_{j+1}-\mathfrak{r}_{j}=-2$. Mais les hypoth\`eses $j\in \mathfrak{J}^+$ et $j+1\in\mathfrak{J}^-$ imposent non seulement $\lambda_{j}>\lambda_{j+1}$ mais aussi que $\lambda_{j}$ et $\lambda_{j+1}$ sont pairs. Donc $\lambda_{j}\geq \lambda_{j+1}+2$. Alors $e_{2}=\lambda_{j}-\lambda_{j+1}-2\geq0$.

Supposons maintenant $j$ pair. Les conditions (11)(c) pour $d=1,2$ disent que $p_{1}(j)=p_{2}(j)=0$. D'apr\`es (4), $\lambda_{j}$ est impair donc $\mathfrak{r}_{j}=0$. On n'a pas $j+1\in \mathfrak{J}^-$ puisque $j+1$ est impair. Donc $\mathfrak{r}_{j+1}\geq0$. On voit alors que $e_{2}=\lambda_{j}-\lambda_{j+1}+\mathfrak{r}_{j+1}-e_{1}$ est positif ou nul sauf si les trois conditions suivantes sont v\'erifi\'ees: $\lambda_{j}=\lambda_{j+1}$, $\mathfrak{r}_{j+1}=0$ et $e_{1}=1$. Supposons ces conditions v\'erifi\'ees. D'apr\`es (13) pour $d=1$, on a $p_{1}(j+1)=1$. Puisque $\lambda_{j}=\lambda_{j+1}$, $\lambda_{j+1}$ est impair. L'\'egalit\'e $\mathfrak{r}_{j+1}=0$ et la relation (8) entra\^{\i}nent alors $p_{2}(j+1)=1$. Pour $d=1,2$, $j+1$ appartient donc \`a un \'el\'ement $\mathfrak{I}_{d}\in \tilde{\mathfrak{Int}}_{d}$. Puisque $p_{d}(j)=0$, $j$ et $j+1$ ne sont pas $d$-li\'es, donc $j+1=j_{min}(\mathfrak{I}_{d})$. Mais alors, (5) nous dit que $j+1$ appartient \`a $\mathfrak{J}^+$, donc $\mathfrak{r}_{j+1}=1$ contrairement \`a l'hypoth\`ese. Cette contradiction conclut. Cela ach\`eve la preuve de l'existence de nos suites $\lambda_{1}$ et $\lambda_{2}$.

Fixons donc de telles suites $\lambda_{1}$ et $\lambda_{2}$. La condition (11) entra\^{\i}ne que ce sont des partitions, c'est-\`a-dire qu'elles sont d\'ecroissantes. Montrons que

(14) il existe des entiers positifs ou nuls $n_{1}$ et $n_{2}$ tels que $n_{1}+n_{2}=n$, que $\lambda_{1}$ appartienne \`a ${\cal P}^{symp,sp}(2n_{1})$ et que $\lambda_{2}$ appartienne \`a ${\cal P}^{orth,sp}(2n_{2})$.

Si les deux derni\`eres conditions sont v\'erifi\'ees, on a forc\'ement $n_{1}+n_{2}=n$. En effet, la relation (9) implique que $S(\lambda_{1})+S(\lambda_{2})+S(\mathfrak{r})=S(\lambda)$ et on a $S(\mathfrak{r})=0$. Pour prouver les deux derni\`eres conditions, on doit prouver que, pour $d=1,2$ et $k\geq1$, les termes $\lambda_{d,2k-1}$ et $\lambda_{d,2k}$ sont de m\^eme parit\'e et que, quand cette parit\'e est celle de $d$, on a $\lambda_{d,2k-1}=\lambda_{d,2k}$. La premi\`ere condition r\'esulte de (10) et (7). Si $\lambda_{d,2k-1}\equiv d\,\,mod\,\,2{\mathbb Z}$, la condition (10) impose $p_{d}(2k-1)=0$. Alors les conditions de (11)(a) sont v\'erifi\'ees pour $j=2k-1$, d'o\`u $\lambda_{d,2k-1}=\lambda_{d, 2k}$. Cela prouve (14).

Gr\^ace \`a (14), on d\'efinit comme en 3.1 les ensembles d'intervalles $\tilde{Int}(\lambda_{1})$, $\tilde{Int}(\lambda_{2})$, les ensembles $J^+$ et $J^-$ et la fonction $\xi$. Montrons que

(15) on a $\{J(\Delta); \Delta\in \tilde{Int}(\lambda_{d})\}=\tilde{\mathfrak{Int}}_{d}$ pour $d=1,2$;  on a $J^+=\mathfrak{J}^+$, $J^-=\mathfrak{J}^-$ et $\xi=\mathfrak{r}$. 

Soit $d=1,2$. La r\'eunion des $J(\Delta)$ quand $\Delta$ d\'ecrit $\tilde{Int}(\lambda_{d})$ est l'ensemble des $j\geq1$ tels que $\lambda_{d,j}$ soit de bonne parit\'e. D'apr\`es (10), c'est l'ensemble des $j\geq1$ tels que $p_{d}(j)=1$. Cet ensemble d'indices est donc d\'ecoup\'e de deux fa\c{c}ons en intervalles: les $J(\Delta)$ pour $\Delta\in \tilde{Int}(\lambda_{d})$ et les $\mathfrak{I}\in \tilde{\mathfrak{Int}}_{d}$. Pour prouver que ces d\'ecoupages co\"{\i}ncident, il suffit de prouver que les ensembles d'\'el\'ements maximaux de ces intervalles co\"{\i}ncident (en admettant ici que l'\'el\'ement maximal d'un intervalle infini est $\infty$). C'est-\`a-dire qu'il suffit de prouver l'\'egalit\'e
$$\{j_{max}(\Delta);\Delta\in \tilde{Int}(\lambda_{d})\}=\{j_{max}(\mathfrak{I}); \mathfrak{I}\in \tilde{\mathfrak{Int}}_{d}\}.$$
L'infini intervient dans les deux ensembles pour $d=1$ et n'intervient dans aucun d'eux pour $d=2$  (d'apr\`es (2) pour l'ensemble de droite). On \'elimine ces termes. Pour $j\geq1$, $j$ n'intervient dans ces ensembles  que si $j$ est pair (d'apr\`es (3) pour celui de droite) et  $\lambda_{d,j}\equiv d+1\,\,mod\,\,2{\mathbb Z}$ autrement dit $p_{d}(j)=1$. Supposons ces conditions v\'erifi\'ees. Alors $j$ intervient dans l'ensemble de gauche si et seulement si $\lambda_{d,j}> \lambda_{d,j+1}$. Si $j$ intervient dans l'ensemble de droite, la condition (11)(b) est v\'erifi\'ee et l'in\'egalit\'e pr\'ec\'edente l'est aussi. Si $j$ n'intervient pas dans l'ensemble de droite, la condition (11)(a) est v\'erifi\'ee et l'in\'egalit\'e pr\'ec\'edente ne l'est pas. Cela d\'emontre l'\'egalit\'e de ces ensembles, d'o\`u la premi\`ere assertion de (15). 

Par d\'efinition, 
$J^+$ est l'ensemble des $j\geq1$ pour lesquels $\lambda_{1,j}$ et $\lambda_{2,j}$ sont de bonne parit\'e  et il existe  $\Delta\in \tilde{Int}(\lambda_{1})\cup \tilde{Int}(\lambda_{2})$ tel que $j=j_{min}(\Delta)$. En utilisant ce que l'on vient de d\'emontrer, il suffit d'appliquer (5) pour conclure $J^+=\mathfrak{J}^+$. On prouve de m\^eme que $J^-=\mathfrak{J}^-$. Alors $\xi=\mathfrak{r}$ par d\'efinition de ces fonctions. Cela prouve (15). 

On  a $Ind(\lambda_{1},\lambda_{2})=\lambda_{1}+\lambda_{2}+\xi$ par d\'efinition, d'o\`u $Ind(\lambda_{1},\lambda_{2})=\lambda$ d'apr\`es (15) et (9). Montrons que

(16) $\lambda_{1}$ et $\lambda_{2}$ induisent r\'eguli\`erement $\lambda$. 

Il s'agit de prouver que tout intervalle relatif est r\'eduit \`a un seul \'el\'ement. Soit $D$ un intervalle relatif. Si $J(D)$ est r\'eduit \`a un seul \'el\'ement, $D$ aussi.  Supposons que $J(D)$ a au moins deux \'el\'ements. Par d\'efinition, il existe un unique $d=1,2$ pour lequel il existe $\Delta_{d}\in Int(\lambda_{d})$ de sorte que $J(D)\subset J(\Delta_{d})$. Pour fixer la notation, on suppose $d=1$. Cela entra\^{\i}ne: pour $j,j+1\in J(D)$, il n'existe pas de $\Delta_{2}\in Int(\lambda_{2})$ tel que $\{j,j+1\}\subset J(\Delta_{2})$. En effet, les extr\'emit\'es $j_{min}(D)$ et $j_{max}(D)$ sont par d\'efinition des \'el\'ements cons\'ecutifs de l'ensemble ${\cal J}$ de 3.1.
Un  $\Delta_{2}$ comme ci-dessus v\'erifierait donc $j_{min}(\Delta_{2})\leq j_{min}(D)$ et $j_{max}(D)\leq j_{max}(\Delta_{2})$, donc $J(D)\subset J(\Delta_{2})$, ce qui est exclu.  On traduit d'apr\`es   (15): il existe $\mathfrak{I}_{1}\in \tilde{\mathfrak{Int}}_{1}$ tel que $J(D)\subset \mathfrak{I}_{1}$ et, pour $j,j+1\in J(D)$, $j$ et $j+1$ ne sont pas $2$-li\'es.  Soient $j,j+1\in J(D)$.  Les indices $j,j+1$ n'\'etant pas $2$-li\'es, ils ne v\'erifient pas les conditions (1)(b), (1)(c) et (1)(d) (cette derni\`ere \'etant de toute fa\c{c}on exclue puisque $\lambda_{j}$ et $\lambda_{j+1}$ sont pairs par d\'efinition des intervalles relatifs). Puisque $j$ et $j+1$ sont $1$-li\'es, ils v\'erifient forc\'ement la condition (1)(a) pour $d=1$. Donc $\lambda_{j}=\lambda_{j+1}$. Cela \'etant vrai pour tout couple $\{j,j+1 \}\subset J(D)$, $\lambda_{j}$ est constant pour $j\in J(D)$. Autrement dit, $D$ est r\'eduit \`a un seul \'el\'ement. 

Montrons que

(17) $\chi_{\lambda_{1},\lambda_{2}}=\chi$. 

On a $\chi_{\lambda_{1},\lambda_{2}}(0)=0 $ par d\'efinition et $\chi(0)=0$ par hypoth\`ese. Soit $i\in Jord^{bp}(\lambda)$. Si $mult_{\lambda}(i)=1$, $\chi_{\lambda_{1};\lambda_{2}}(i)=0 $ par d\'efinition et $\chi(i)=0$ par hypoth\`ese. Supposons $mult_{\lambda}(i)\geq2$. Comme dans la preuve de (16), il existe un unique $d=1,2$ de sorte qu'il existe $ \Delta_{d}\in  \tilde{Int}(\lambda_{d})$ tel que $J(i)\subset J(\Delta_{d})$. On a $\chi_{\lambda_{1},\lambda_{2}}(i)=d+1$ par d\'efinition. Toujours comme dans la preuve de (16), pour $j,j+1\in J(i)$, la condition (1)(a) est v\'erifi\'ee pour ce $d$. Alors $\chi(i)=d+1$. D'o\`u (17).

Montrons que  

(18) $\zeta(\lambda_{1})+\zeta(\lambda_{2})=\zeta(\lambda)+\xi$.

 On a d\'efini en 3.1 les ensembles  $P^+_{\lambda_{1},\lambda_{2}}(\lambda)$ et $P^-_{\lambda_{1},\lambda_{2}}(\lambda)$ et la suite $\zeta_{\lambda_{1},\lambda_{2}}(\lambda)$.  Puisque $\lambda_{1}$ et $\lambda_{2}$ induisent r\'eguli\`erement $\lambda$, on a les \'egalit\'es $P^+_{\lambda_{1},\lambda_{2}}(\lambda)=P^+(\lambda)$, $P^-_{\lambda_{1},\lambda_{2}}(\lambda)=P^-(\lambda)$. Donc $\zeta_{\lambda_{1},\lambda_{2}}(\lambda)=\zeta(\lambda)$. Alors le lemme 3.1 implique (18). 

 L'\'egalit\'e (18) entra\^{\i}ne
 $$\lambda_{1}+\zeta(\lambda_{1})+\lambda_{2}+\zeta(\lambda_{2})=\lambda_{1}+\lambda_{2}+\xi+\zeta(\lambda)=\lambda+\zeta(\lambda).$$
 Le lemme 2.7 et l'assertion 2.7(4) transforment cette \'egalit\'e en 
  $$^td(\lambda_{1})+{^td(\lambda_{2})}={^td(\lambda)},$$
 d'o\`u $d(\lambda_{1})\cup d(\lambda_{2})=d(\lambda)$. Cela ach\`eve la d\'emonstration. $\square$
 
 \bigskip

\subsection{Les fonctions $\tau^{\zeta},\delta^{\zeta}$}

 Soient $n_{1},n_{2}\in {\mathbb N}$ tels que $n_{1}+n_{2}=n$ et soient  $\lambda_{1}\in {\cal P}^{symp,sp}(2n_{1})$ et $\lambda_{2}\in {\cal P}^{orth,sp}(2n_{2})$. Soit $\lambda$ l'induite endoscopique de $\lambda_{1}$ et $\lambda_{2}$. On consid\`ere de plus des \'el\'ements $\iota_{1}=(\tau_{1},\delta_{1})\in {\cal F}am(\lambda_{1})$ et $\iota_{2}=(\tau_{2},\delta_{2})\in {\cal F}am(\lambda_{2})$. On pose $r_{1}=r(\tau_{1},\delta_{1})$, $r_{2}=r(\tau_{2},\delta_{2})$. 

Pour $d=1,2$ et $\Delta\in \tilde{Int}(\lambda_{d})$, on note $\Delta^+$ le plus petit $\Delta'\in Int(\lambda_{d})$ tel que $\Delta'>\Delta$, pour peu qu'il existe un tel $\Delta'$ (sinon, $\Delta^+$ n'existe pas). Pour $D\in \tilde{Int}_{\lambda_{1},\lambda_{2}}(\lambda)$, on d\'efinit $D^+$ de fa\c{c}on similaire. 

Pour $D\in Int_{\lambda_{1},\lambda_{2}}(\lambda)$ et pour $d=1,2$, consid\'erons l'ensemble des $\Delta\in \tilde{Int}(\lambda_{d})$ tels que $j_{max}(D)\leq j_{max}(\Delta)$ (ici, on pose par convention $j_{max}(\Delta_{1,min})=\infty$ o\`u $\Delta_{1,min}$ est le plus petit \'el\'ement de $\tilde{Int}(\lambda_{1})$). Si cet ensemble est non vide (ce qui est le cas si $d=1$ par la convention que l'on vient de poser), on note $\Delta_{d}(D)$ son plus grand \'el\'ement. On pose $\Delta_{1}(D_{min})=\Delta_{1,min}$ tandis que $\Delta_{2}(D_{min})$ n'existe pas. Si $\Delta_{2}(D)$ n'existe pas et si $Int(\lambda_{2})$ n'est pas vide, on note $\Delta_{2}(D)^+$ le plus petit \'el\'ement de $Int(\lambda_{2})$ (si $Int(\lambda_{2})$ est vide, $\Delta_{2}(D)$ et $\Delta_{2}(D)^+$ n'existent pas). 

Pour $\zeta=\pm$, on d\'efinit  une fonction $\delta^{\zeta}\in ({\mathbb Z}/2{\mathbb Z})^{Int_{\lambda_{1},\lambda_{2}}(\lambda)}$ par les formules ci-dessous. Soit $D\in Int_{\lambda_{1},\lambda_{2}}(\lambda)$. On pose $\Delta_{d}=\Delta_{d}(D)$ pour $d=1,2$. Ce terme existe toujours dans chaque cas ci-dessous. Par contre, $\Delta_{d}^+$ n'existe pas toujours. Dans ce cas, on consid\`ere que   $\delta_{d}(\Delta_{d}^+)=0$. 
On \'ecrit les formules comme des \'egalit\'es, en fait, il s'agit de congruences modulo $2{\mathbb Z}$. On pose

 si $j_{max}(D)\in J^+$, $\delta^+(D)=\tau_{1}(\Delta_{1})+\tau_{2}(\Delta_{2})+r_{1}+r_{2}+1$, $\delta^-(D)=\delta^+(D)+1$;
 
 si $j_{max}(D)\in J^-$, $\delta^+(D)=\delta^-(D)=\delta_{1}(\Delta_{1})+\delta_{2}(\Delta_{2})$;
 
 si $j_{max}(D)\not\in J^+\cup J^-$ et $J(D)\subset J(\Delta_{1})$, $\delta^+(D)=\delta^-(D)=\delta_{1}(\Delta_{1})+\delta_{2}(\Delta_{2}^+)$;
 
 si $j_{max}(D)\not\in J^+\cup J^-$ et $J(D)\subset J(\Delta_{2})$, $\delta^+(D)=\delta^-(D)=\delta_{1}(\Delta_{1}^+)+\delta_{2}(\Delta_{2})$.
 
 Avec les m\^emes notations, on d\'efinit une fonction $\tau^{\zeta}\in ({\mathbb Z}/2{\mathbb Z})^{\tilde{Int}_{\lambda_{1},\lambda_{2}}(\lambda)}$  par
 
 si $\vert J(D)\vert \geq2$ et $J(D)\subset J(\Delta_{1})$, $\tau^+(D)=\tau^-(D)=\tau_{1}(\Delta_{1})+\delta_{2}(\Delta_{2}^+)+r_{2}$;
 
 si $\vert J(D)\vert \geq2$ et $J(D)\subset J(\Delta_{2})$, $\tau^+(D) =\delta_{1}(\Delta_{1}^+)+\tau_{2}(\Delta_{2})+r_{1}$, $\tau^-(D)=\tau^+(D)+1$;
 
 si $\vert J(D)\vert =1$ et $j_{min}(D)=j_{max}(D)\in J^+$, $\tau^+(D)=\tau^-(D)=\tau_{1}(\Delta_{1})+\delta_{2}(\Delta_{2}^+)+r_{2}$;

 si $\vert J(D)\vert =1$ et $j_{min}(D)=j_{max}(D)\in J^-$, $\tau^+(D)=\tau^-(D)=\tau_{1}(\Delta_{1})+\delta_{2}(\Delta_{2})+r_{2}$.
 
 Tous ces cas sont exclusifs l'un de l'autre. 
On a \'evidemment
 
 (1) $\delta^-(D)=\delta^+(D)+1$ si et seulement si  $j_{max}(D)\in J^+$; $\tau^-(D)=\tau^+(D)+1$ si et seulement si $\vert J(D)\vert \geq2$ et $J(D)\subset J(\Delta_{2})$.
 
 On a aussi
 
(2)  $\tau^+(D_{min})=\tau^-(D_{min})=0$.

En effet, $J(D_{min})$ est infini. Il ne peut qu'\^etre contenu dans $J(\Delta_{1,min})$. Donc $\tau^+(D_{min})=\tau^-(D_{min})=\tau_{1}(\Delta_{1}(D_{min}))+\delta_{2}(\Delta_{2}(D_{min})^+)+r_{2}$. On a $\Delta_{1}(D_{min})=\Delta_{1,min}$ et $\Delta_{2}(D_{min})$ n'existe pas. On a $\tau_{1}(\Delta_{1,min})=0$. D'apr\`es 2.3(2) et nos conventions,   
$\delta_{2}(\Delta_{2}(D_{min})^+)=r_{2}$. D'o\`u (2).
  
 Pour $\zeta=\pm$, posons
  $$C^{\zeta}=\sum_{D\in Int_{\lambda_{1},\lambda_{2}}(\lambda)}(1-(-1)^{\tau^{\zeta}(D)})((-1)^{\delta^{\zeta}(D)}-(-1)^{\delta^{\zeta}(D^+)}).$$
   Ici encore, on  consid\`ere que $\delta^{\zeta}(D^+)=1$ si $D^+$ n'existe pas. 
  On a
 $$(3)\qquad C^{\zeta}=\left\lbrace\begin{array}{cc}2(r_{1}+\zeta r_{2}),&\,\,{\text si}\,\,r_{1}+r_{2}\text{\,\,est\,\,pair},\\ -2(r_{1}+\zeta r_{2}+1),&\,\,{\text si}\,\,r_{1}+r_{2}\text{\,\,est\,\,impair}.\\ \end{array}\right.$$
 Cela r\'esulte de \cite{W1} XI.24, \`a ceci pr\`es que les hypoth\`eses de cette r\'ef\'erence \'etaient plus restrictives que les n\^otres. On renvoie pour ce probl\`eme aux explications que l'on donnera apr\`es la proposition du paragraphe suivant.
 
 \bigskip
 
 \subsection{Le r\'esultat de \cite{W1}}
 Les donn\'ees sont les m\^emes que dans le paragraphe pr\'ec\'edent. Pour $d=1,2$,  le couple $\iota_{d}=(\tau_{d},\delta_{d})$  provient d'un symbole $\Lambda_{d}$ dans la famille de $\lambda_{d}$.  On note $(r_{d},\rho_{d})$ l'\'el\'ement de $\Sigma_{n_{1},imp}$ si $d=1$, $\Sigma_{n_{2},pair}$ si $d=2$, tel que $symb(r_{d},\rho_{d})=\Lambda_{d}$. On pose $N_{1}=n_{1}-r_{1}^2-r_{1}$, $N_{2}=n_{2}-r_{2}^2$. On fixe un \'el\'ement $\zeta\in \{\pm\}$. Si $\zeta=1$, on pose $h^+=r_{1}+\vert r_{2}\vert $, $h^-=sup(r_{1}-\vert r_{2}\vert,\vert r_{2}\vert -r_{1}-1)$. Si $\zeta=-1$, on pose $h^+=sup(r_{1}-\vert r_{2}\vert,\vert r_{2}\vert -r_{1}-1)$, $h^-=r_{1}+\vert r_{2}\vert $. On v\'erifie que $h^+(h^++1)/2+h^-(h^-+1)/2=r_{1}^2+r_{1}+r_{2}^2$. On fixe des entiers $n^+,n^-\in {\mathbb N}$ tels que $n^++n^-=n$, $n^+\geq h^+(h^++1)/2$, $n^-\geq h^-(h^-+1)/2$ et on pose $N^+=n^+-h^+(h^++1)/2$, $N^-=n^--h^-(h^-+1)/2$. On a $N^++N^-=N_{1}+N_{2}$. On d\'efinit un quadruplet d'entiers ${\bf a}=(a_{1}^+,a_{1}^-,a_{2}^+,a_{2}^-)$ par les formules suivantes
 
 ${\bf a}=(0,0,0,1)$ si $\zeta=1$ et $r_{1}\geq\vert r_{2}\vert $;
 
 ${\bf a}=(0,0,1,0)$ si $\zeta=-1$ et $r_{1}\geq \vert r_{2}\vert $;
 
$ {\bf a}=(0;1,0,0)$ si $\zeta=1$ et $r_{1}<\vert r_{2}\vert $;
 
 ${\bf a}=(1,0,0,0)$ si $\zeta=-1$ et $r_{1}<\vert r_{2}\vert $.
 
Avec les m\^emes notations qu'en  1.2, on d\'efinit une repr\'esentation $\Pi^{\zeta}(\iota_{1},\iota_{2})$ de $W_{N^+}\times W_{N^-}$ par la formule
$$\Pi^{\zeta}(\iota_{1},\iota_{2})=\oplus_{{\bf N}\in {\cal N}}ind_{W_{{\bf N}}} ^{W_{N^+}\times W_{N^-}}(sgn_{CD}^{{\bf a}}\otimes res_{W_{{\bf N}}}^{W_{N_{1}}\times W_{N_{2}}}(\rho_{1}\otimes \rho_{2})).$$

On  note ${\cal I}^{\zeta}(\iota_{1},\iota_{2})$ l'ensemble des quadruplets $(\lambda^+,\epsilon^+,\lambda^-,\epsilon^-)
 \in \boldsymbol{{\cal P}^{symp}}(2n^+) \times \boldsymbol{{\cal P}^{symp}}(2n^-)$ v\'erifiant les conditions suivantes:
 
 (1)  $k_{\lambda^+,\epsilon}=h^+$, $k_{\lambda^-,\epsilon^-}=h^-$;
 
 (2) la repr\'esentation $\rho_{\lambda^+,\epsilon^+}\otimes \rho_{\lambda^-,\epsilon^-}$ de $W_{N^+}\times W_{N^-}$ intervient dans $\Pi^{\zeta}(\iota_{1},\iota_{2})$ avec une multiplicit\'e strictement positive. 
 
 Pour poser la d\'efinition suivante, on a besoin d'introduire deux notations. Pour $D\in Int_{\lambda_{1},\lambda_{2}}(\lambda)$, notons $i_{min}(D)$ le plus petit \'el\'ement de $D$. On a $i_{min}(D)\geq1$ puisque $D\not=D_{min}$. Pour toute partition $\mu$, on pose $mult_{\geq D}(\mu)=\sum_{i'\in {\mathbb N},i'\geq i_{min}(D)}mult_{i'}(\mu)$. D'autre part, on pose $\nu=1$ si $r_{2}\geq0$, $\nu=-1$ si $r_{2}<0$. 
 
 On  note ${\cal I}^{\zeta,max}(\iota_{1},\iota_{2})$ l'ensemble des quadruplets $(\lambda^+,\epsilon^+,\lambda^-,\epsilon^-)
 \in \boldsymbol{{\cal P}^{symp}}(2n^+) \times \boldsymbol{{\cal P}^{symp}}(2n^-)$ v\'erifiant les conditions suivantes:
 
 (3) $\lambda^+\cup \lambda^-=\lambda$;
 
 (4) pour tout $D\in Int_{\lambda_{1},\lambda_{2}}(\lambda)$, on a
 $$mult_{\lambda^+}(\geq D)\equiv \delta^{\zeta\nu}(D)\,\mod\,2{\mathbb Z},\,\, {\text et\,\,}mult_{\lambda^-}(\geq D)\equiv \delta^{-\zeta\nu}(D)\,\mod\,2{\mathbb Z};$$
 
 (5) pour tout $D\in \tilde{Int}_{\lambda_{1},\lambda_{2}}(\lambda)$ et tout $i\in D$ tel que $i\not=0$ et $mult_{\lambda^+}(i)>0$, resp. $mult_{\lambda^-}(i)>0$, on a
 $$\epsilon^+_{i}=(-1)^{\tau^{\zeta\nu}(D)},\,\,{\text resp. \,\,}\epsilon^-_{i}=(-1)^{\tau^{-\zeta\nu}(D)}.$$
 
 Dans ces formules, on a \'evidemment identifi\'e les signes $\pm$ des d\'efinitions de $\tau^+$, $\tau^-$ etc... \`a des \'el\'ements de $\{\pm 1\}$. On a montr\'e en \cite{W1} XI.29 remarque 4 que, sous l'hypoth\`ese (3), les deux congruences de (4) \'etaient \'equivalentes. 
 
 \ass{Proposition}{(i) Soit $(\lambda^+,\epsilon^+,\lambda^-,\epsilon^-)\in {\cal I}^{\zeta}(\iota_{1},\iota_{2})$. Alors $\lambda^+\cup \lambda^-\leq \lambda$.
 
 (ii) L'ensemble ${\cal I}^{\zeta,max}(\iota_{1},\iota_{2})$ est \'egal au sous-ensemble des $(\lambda^+,\epsilon^+,\lambda^-,\epsilon^-)\in {\cal I}^{\zeta}(\iota_{1},\iota_{2})$ tels que $\lambda^+\cup \lambda^-=\lambda$. Pour $(\lambda^+,\epsilon^+,\lambda^-,\epsilon^-)\in {\cal I}^{\zeta,max}(\iota_{1},\iota_{2})$, la repr\'esentation $\rho_{\lambda^+,\epsilon^+}\otimes \rho_{\lambda^-,\epsilon^-}$ intervient avec multiplicit\'e $1$ dans $\Pi^{\zeta}(\iota_{1},\iota_{2})$.}
 
 Cela r\'esulte de \cite{W1} propositions XI.28 et XI.29, ainsi qu'on l'a expliqu\'e dans la preuve de la  proposition XII.7 de cette r\'ef\'erence. A ceci pr\`es qu'alors, les hypoth\`eses sur $\iota_{2}$ \'etaient restrictives: on supposait que $r_{2}$ \'etait pair et positif ou nul; dans le cas $r_{2}=0$, on supposait que le symbole $(X,Y)$ correspondant \`a $\iota_{2}$ v\'erifiait $X\geq Y$ pour l'ordre lexicographique. En fait, cette derni\`ere hypoth\`ese \'etait utilis\'ee dans  d'autres passages de \cite{W1} mais pas dans les d\'emonstrations des propositions utilis\'ees. Pour traiter le cas  o\`u $r_{2}$ est  impair et positif, il n'y a pas d'autre m\'ethode que de reprendre la d\'emonstration. C'est ce que l'on a fait mais elle est trop longue pour la r\'ecrire. Le cas o\`u $r_{2}<0$ se d\'eduit  du cas $r_{2}>0$ de la fa\c{c}on suivante. On suppose donc $r_{2}<0$. On a dit que $\iota_{2}$ correspondait \`a un symbole $\Lambda_{2}=(X_{2},Y_{2})$, puis \`a un couple $(r_{2},\rho_{2})$. Inversement, on voit que $(-r_{2},\rho_{2})$ correspond au symbole $\Lambda'_{2}=(Y_{2},X_{2})$, puis \`a un \'el\'ement $\iota'_{2} \in {\cal F}am(\lambda_{2})$.  Quand on remplace $\iota_{2}$ par $\iota'_{2}$ dans les constructions ci-dessus,  la repr\'esentation $\Pi^{\zeta}(\iota_{1},\iota_{2})$ ne change pas. Donc la proposition ci-dessus \'etant v\'erifi\'ee pour $\iota'_{2}$, elle le restera pourvu que l'on ait les \'egalit\'es ${\cal I}^{\zeta}(\iota_{1},\iota_{2})={\cal I}^{\zeta}(\iota_{1},\iota'_{2})$ et ${\cal I}^{\zeta,max}(\iota_{1},\iota_{2})={\cal I}^{\zeta,max}(\iota_{1},\iota'_{2})$. La premi\`ere \'egalit\'e est claire d'apr\`es (1) et (2). La deuxi\`eme ne l'est pas car les fonctions $\tau^{\pm}$ et $\delta^{\pm}$ d\'ependent de $\iota_{2}$. Mais, puisqu'on passe de $\Lambda_{2}$ \`a $\Lambda'_{2}$ en permutant $X_{2}$ et $Y_{2}$, on voit sur les formules de 2.2  que changer $\iota_{2}$ en $\iota'_{2}$ ne change pas $\delta_{2}$ et remplace $\tau_{2}$ par $\tau_{2}+1$. On voit ensuite sur les formules de 3.3 que cela \'echange les couples $(\tau^+,\delta^+)$ et $(\tau^-,\delta^-)$. Mais alors, parce qu'il figure dans les conditions (4) et (5) un signe terme $\nu$, qui vaut $1$ pour $\iota'_{2}$ et $-1$ pour $\iota_{2}$, on voit que ces conditions ne changent pas quand on remplace $\iota_{2}$ par $\iota'_{2}$. C'est ce qu'on voulait. 
 
  \bigskip
 
 \subsection{R\'eciproque de la construction des fonctions $\tau^{\zeta}$ et $\delta^{\zeta}$}
 
  Soient $n_{1},n_{2}\in {\mathbb N}$ tels que $n_{1}+n_{2}=n$ et soient  $\lambda_{1}\in {\cal P}^{symp,sp}(2n_{1})$ et $\lambda_{2}\in {\cal P}^{orth,sp}(2n_{2})$. Notons $\lambda$ l'induite endoscopique de $\lambda_{1}$ et $\lambda_{2}$. Pour $\iota_{1}=(\tau_{1},\delta_{1})\in {\cal F}am(\lambda_{1})$ et $\iota_{2}=(\tau_{2},\delta_{2})\in {\cal F}am(\lambda_{2})$, on a construit en 3.3 des fonctions $\tau^{\zeta}$ et $\delta^{\zeta}$ pour $\zeta=\pm$. Dans ce paragraphe, il convient de les noter plus pr\'ecis\'ement $\tau^{\zeta}_{\iota_{1},\iota_{2}}$ et $\delta^{\zeta}_{\iota_{1},\iota_{2}}$. On note aussi $C^{\zeta}_{\iota_{1},\iota_{2}}$ la somme d\'efinie en 3.3.
 
 Soient $r_{1}\in {\mathbb N}$, $r_{2}\in {\mathbb Z}$ et, pour $\zeta=\pm$, soient $\tau^{\zeta}\in ({\mathbb Z}/2{\mathbb Z})^{\tilde{Int}_{\lambda_{1},\lambda_{2}}(\lambda)}$ et $\delta^{\zeta}\in ({\mathbb Z}/2{\mathbb Z})^{Int_{\lambda_{1},\lambda_{2}}(\lambda)}$. On pose
  $$C^{\zeta}=\sum_{D\in Int_{\lambda_{1},\lambda_{2}}(\lambda)}(1-(-1)^{\tau^{\zeta}(D)})((-1)^{\delta^{\zeta}(D)}-(-1)^{\delta^{\zeta}(D^+)}).$$
 
 On suppose que ces donn\'ees v\'erifient les conditions
 
 (1) $\delta^-(D)=\delta^+(D)+1$ si et seulement si $j_{max}(D)\in J^+$; $\tau^-(\delta)=\tau^+(D)+1$ si et seulement si $\vert J(D)\vert \geq2$ et $J(D)\subset J(\Delta_{2}(D))$;
 
 (2) $\tau^+(D_{min})=\tau^-(D_{min})=0$;
   
$$ (3) \qquad C^{\zeta}=\left\lbrace\begin{array}{cc}2(r_{1}+\zeta r_{2}),&\,\,{\text si}\,\,r_{1}+r_{2}\text{\,\,est\,\,pair},\\ -2(r_{1}+\zeta r_{2}+1),&\,\,{\text si}\,\,r_{1}+r_{2}\text{\,\,est\,\,impair}.\\ \end{array}\right.$$
 
 \ass{Lemme}{Sous ces hypoth\`eses, il existe d'uniques $\iota_{1}=(\tau_{1},\delta_{1})\in {\cal F}am(\lambda_{1})$ et $\iota_{2}=(\tau_{2},\delta_{2})\in {\cal F}am(\lambda_{2})$ tels que, pour $\zeta=\pm$,  on ait les \'egalit\'es $\tau^{\zeta}=\tau^{\zeta}_{\iota_{1},\iota_{2}}$ et $\delta^{\zeta}=\delta^{\zeta}_{\iota_{1},\iota_{2}}$. 
 De plus, on a $r_{1}=r(\tau_{1},\delta_{1})$ et $r_{2}=r(\tau_{2},\delta_{2})$. }
 
 Preuve. S'il existe  $(\tau_{1},\delta_{1})$ et $(\tau_{2},\delta_{2})$ v\'erifiant la premi\`ere assertion de l'\'enonc\'e, les fonctions $\tau^{\zeta}$ et $\delta^{\zeta}$ sont donn\'ees par les formules du paragraphe 3.3, o\`u l'on remplace $r_{1}$ et $r_{2}$ par $r'_{1}= r(\tau_{1},\delta_{1})$ et $r'_{2}=r(\tau_{2},\delta_{2})$. Remarquons que ces formules ne d\'ependent que des images de $r'_{1}$ et $r'_{2}$ dans ${\mathbb Z}/2{\mathbb Z}$. On note symboliquement $(X_{r'_{1},r'_{2}})$ ces formules. 
 
 Commen\c{c}ons par prouver que, pour deux \'el\'ements donn\'es $r'_{1},r'_{2}\in {\mathbb Z}/2{\mathbb Z}$,  il existe d'uniques 
 $(\tau_{1},\delta_{1})\in ({\mathbb Z}/2{\mathbb Z})^{\tilde{Int}(\lambda_{1})}\times ({\mathbb Z}/2{\mathbb Z})^{Int(\lambda_{1})}$, $(\tau_{2},\delta_{2})\in ({\mathbb Z}/2{\mathbb Z})^{Int(\lambda_{2})}\times ({\mathbb Z}/2{\mathbb Z})^{Int(\lambda_{2})}$ telles que les formules 
 $(X_{r'_{1},r'_{2}})$ soient v\'erifi\'ees.
  Remarquons que l'on peut consid\'erer uniquement les formules exprimant $\tau^+$ et $\delta^+$: celles concernant $\tau^-$ et $\delta^-$ s'en d\'eduisent d'apr\`es l'hypoth\`ese (1).

Pour $D\in \tilde{Int}_{\lambda_{1},\lambda_{2}}(\lambda)$ et pour $d=1,2$, notons $\mathfrak{T}_{d}(\geq D)$ l'ensemble des $\Delta_{d}\in \tilde{Int}(\lambda_{d})$  tels que $j_{min}(\Delta_{d})\leq 
j_{max}(D)$ 
  (en convenant que $j_{max}(D_{min})=\infty$) et notons $\mathfrak{D}_{d}(\geq D)$ l'ensemble des $\Delta_{d}\in Int(\lambda_{d})$ tels que $j_{max}(\Delta_{d}) \leq j_{max}(D)$. Remarquons que $\mathfrak{T}_{d}(\geq D_{min}) = \tilde{Int}(\lambda_{d})$ et $\mathfrak{D}_{d}(\geq D_{min}) = Int(\lambda_{d})$. 
 Pour deux intervalles relatifs  $D>D'$, il est clair que $\mathfrak{T}_{d}(\geq D)$ est inclus dans $\mathfrak{T}_{d}(\geq D')$ et que $\mathfrak{D}_{d}(\geq D)$ est inclus dans $\mathfrak{D}_{d}(\geq D')$. On pose
  $$\mathfrak{T}_{d}(D)=\mathfrak{T}_{d}(\geq D)-\mathfrak{T}_{d}(\geq D^+),\,\, \mathfrak{D}_{d}(D)=\mathfrak{D}_{d}(\geq D)-\mathfrak{D}_{d}(\geq D^+),$$
  avec la convention $\mathfrak{T}_{d}(\geq D^+) =\mathfrak{D}_{d}(\geq D^+)=\emptyset$ si $D^+$ n'existe pas, c'est-\`a-dire si $D$ est l'intervalle relatif maximal. Cette d\'efinition entra\^{\i}ne:
  
  (4) pour deux intervalles relatifs $D\not=D'$, on a $\mathfrak{T}_{d}(D)\cap \mathfrak{T}_{d}(D')=\emptyset$ et $\mathfrak{D}_{d}(D)\cap \mathfrak{D}_{d}(D')=\emptyset$.
  
  Montrons que
  
  (5) $\mathfrak{T}_{d}(D)$ est l'ensemble des $\Delta_{d}\in \tilde{Int}(\lambda_{d})$ tels que 
  $j_{min}(\Delta_{d})\in \{j_{min}(D),j_{max}(D)\}$ et $\mathfrak{D}_{d}(D)$ est l'ensemble des $\Delta_{d}\in Int(\lambda_{d})$ tels que 
  $j_{max}(\Delta_{d})\in \{j_{min}(D),j_{max}(D)\}$; ces ensembles ont au plus un \'el\'ement. 
  
  Soit $\Delta_{d}\in \tilde{Int}(\lambda_{d}) $, supposons  $j_{min}(\Delta_{d})\in \{j_{min}(D),j_{max}(D)\}$. Alors $j_{min}(\Delta_{d})\leq j_{max}(D)$ et  $\Delta_{d}$ appartient \`a $\mathfrak{T}_{d}(\geq D)$. Si $D$ est l'intervalle relatif maximal, cela entra\^{\i}ne $\Delta_{d}\in \mathfrak{T}_{d}(D)$. Sinon, on a $j_{max}(D^+)<j_{min}(D)\leq j_{min}(\Delta_{d})$ donc $\Delta_{d}$ n'appartient pas \`a $\mathfrak{T}_{d}(\geq D^+)$. D'o\`u $\Delta_{d}\in \mathfrak{T}_{d}(D)$. R\'eciproquement, supposons $\Delta_{d}\in \mathfrak{T}_{d}(D)$. L'entier $j_{min}(\Delta_{d})$ appartient \`a l'ensemble ${\cal J}$ de 3.1. D'apr\`es 3.1(3), il existe $D'\in \tilde{Int}_{\lambda_{1},\lambda_{2}}(\lambda)$ tel que $j_{min}(\Delta_{d})\in \{j_{min}(D'),j_{max}(D')\}$. D'apr\`es ce que l'on vient de prouver, on a $\Delta_{d}\in \mathfrak{T}_{d}(D')$. Alors (4) entra\^{\i}ne $D'=D$, donc $j_{min}(\Delta_{d})\in \{j_{min}(D),j_{max}(D)\}$. Cela prouve la premi\`ere assertion de (4).  Supposons encore que $\Delta_{d}\in \mathfrak{T}_{d}(D)$ et consid\'erons un  intervalle $\Delta'_{d}\in \tilde{Int}(\lambda_{d})$ distinct de $\Delta_{d}$. Si $\Delta'_{d}> \Delta_{d}$, on a $ j_{max}(\Delta'_{d})<j_{min}(\Delta_{d})\leq j_{max}(D)$. Le nombre $j_{max}(\Delta'_{d})$ appartient \`a ${\cal J}$. Par d\'efinition des intervalles relatifs, $j_{min}(D)$ et $j_{max}(D)$ sont soit \'egaux, soit des \'el\'ements cons\'ecutifs de ${\cal J}$. Cela entra\^{\i}ne en tout cas l'in\'egalit\'e $j_{max}(\Delta'_{d})\leq j_{min}(D)$. Puisque $j_{min}(\Delta'_{d})<j_{max}(\Delta'_{d})$, on a donc $j_{min}(\Delta'_{d})\not\in \{j_{min}(D),j_{max}(D)\}$, d'o\`u $\Delta'_{d}\not\in \mathfrak{T}_{d}(D)$. Si maintenant $\Delta'_{d}<\Delta_{d}$, on   a $j_{min}(D)\leq j_{min}(\Delta_{d})<j_{max}(\Delta_{d})$. Comme ci-dessus, on en d\'eduit $j_{max}(D)\leq j_{max}(\Delta_{d})$, puis $j_{max}(D)<j_{min}(\Delta'_{d})$ et on conclut $\Delta'_{d}\not\in \mathfrak{T}_{d}(\geq D)$.  Donc $\mathfrak{T}_{d}(D)$ a au plus un \'el\'ement. Les assertions concernant $\mathfrak{D}_{d}(D)$ se d\'emontrent de la m\^eme fa\c{c}on. Cela prouve (5).
  
     On va montrer que, pour tout intervalle relatif $D$ les formules $(X_{r'_{1},r'_{2}})$ exprimant $\tau^+(D)$ et $\delta^+(D)$, d'une part ne font intervenir des $\tau_{d}(\Delta_{d})$ que pour des $\Delta_{d}\in \mathfrak{T}_{d}(\geq D)$ et des $\delta_{d}(\Delta_{d})$ que pour des $\Delta_{d}\in \mathfrak{D}_{d}(\geq D)$, d'autre part que, quand $\mathfrak{T}_{d}(D)$, resp. $\mathfrak{D}_{d}(D)$, est non vide, elles font intervenir 
      $\tau_{d}(\Delta_{d})$, resp. $\delta_{d}(\Delta_{d})$, pour l'unique \'el\'ement $\Delta_{d}$ de cet ensemble.   On \'etudie les diff\'erents cas possibles, pour $D\in Int_{\lambda_{1},\lambda_{2}}(\lambda)$. On suppose d'abord $D\not=D_{min}$. On pose simplement $\Delta_{d}=\Delta_{d}(D)$. 
  
  (a) Supposons que $\vert J(D)\vert=1$ et  que $j_{min}(D)=j_{max}(D)\in J^+$. Dans ce cas, on a $j_{max}(D)=j_{min}(\Delta_{1})=j_{min}(\Delta_{2})$. D'apr\`es (5), on a $\mathfrak{T}_{1}(D)=\{\Delta_{1}\}$, $\mathfrak{T}_{2}(D)=\{\Delta_{2}\}$. Si $\Delta'_{1}\in \mathfrak{D}_{1}(D)$, on a $j_{max}(D)=j_{min}(D)\in J(\Delta'_{1})$, donc $J(\Delta'_{1})\cap J(\Delta_{1})\not=\emptyset$, donc $\Delta'_{1}=\Delta_{1}$. Or $j_{max}(\Delta_{1})> j_{min}(\Delta_{1})=j_{max}(D)$, donc $\Delta_{1}\not\in \mathfrak{D}_{1}(D)$. Donc $\mathfrak{D}_{1}(D)=\emptyset$ et, de m\^eme, $\mathfrak{D}_{2}(D)=\emptyset$. Par ailleurs, si $\Delta_{2}^+$ existe, on a $j_{max}(\Delta_{2}^+)< j_{min}(\Delta_{2})=j_{max}(D)$, donc $\Delta_{2}^+\in \mathfrak{D}_{2}(\geq D)$. Enfin, les formules dans notre cas sont
  $$\delta^+(D)=\tau_{1}(\Delta_{1})+\tau_{2}(\Delta_{2})+r'_{1}+r'_{2}+1,$$
  $$\tau^+(D)=\tau_{1}(\Delta_{1})+\delta_{2}(\Delta_{2}^+)+r'_{2}.$$
  On voit que les propri\'et\'es requises sont v\'erifi\'ees.
  
  (b) Supposons que $\vert J(D)\vert=1$ et  que $j_{min}(D)=j_{max}(D)\in J^-$.  Ce cas est similaire au pr\'ec\'edent. On a cette fois $j_{max}(D)=j_{max}(\Delta_{1})=j_{max}(\Delta_{2})$. On a $\mathfrak{T}_{d}(D)=\emptyset$ pour $d=1,2$, $\mathfrak{D}_{1}(D)=\{\Delta_{1}\}$, $\mathfrak{D}_{2}(D)=\{\Delta_{2}\}$ et $\Delta_{1}\in \mathfrak{T}_{1}(\geq D)$.  Les formules sont
  $$\delta^+(D)=\delta_{1}(\Delta_{1})+\delta_{2}(\Delta_{2}),$$
  $$\tau^+(D)=\tau_{1}(\Delta_{1})+\delta_{2}(\Delta_{2})+r'_{2}.$$
  Les propri\'et\'es requises sont v\'erifi\'ees.
  
  (c) Supposons que $\vert J(D)\vert \geq2$, que $J(D)\subset J(\Delta_{1})$ et que $j_{min}(D)$ et $j_{max}(D)$ soient impairs. Puisque ces termes appartiennent \`a l'ensemble ${\cal J}$, l'imparit\'e impose qu'ils sont  de la forme $j_{min}(D)=j_{min}(\Delta'_{d'})$ et $j_{max}(D)=j_{min}(\Delta''_{d''})$ pour des entiers $d',d''=1,2$ et des intervalles $\Delta'_{d'}\in \tilde{Int}(\lambda_{d'})$ et $\Delta''_{d''}\in \tilde{Int}(\lambda_{d''})$. Si $d'=2$, puisque $j_{min}(D)$ et $j_{max}(D)$ sont des \'el\'ements cons\'ecutifs de ${\cal J}$, on a  $j_{max}(D)\leq j_{max}(\Delta'_{2})$, d'o\`u $J(D)\subset J(\Delta'_{2})$, ce qui est interdit par d\'efinition des intervalles et par l'hypoth\`ese $J(D)\subset J(\Delta_{1})$. Donc $d'=1$ et forc\'ement $\Delta'_{1}=\Delta_{1}$, c'est-\`a-dire $j_{min}(D)=j_{min}(\Delta_{1})$. Si $d''=1$, on a $J(\Delta''_{1})\cap J(\Delta_{1})\not=\emptyset$ donc $\Delta'_{1}=\Delta_{1}$. Mais $j_{min}(\Delta_{1})\leq j_{min}(D)$ par hypoth\`ese, donc $j_{min}(\Delta_{1})$ ne peut pas \^etre \'egal \`a $j_{max}(D)$. Donc $d''=2$ et forc\'ement $\Delta'_{2}=\Delta_{2}$. C'est-\`a-dire $j_{max}(D)=j_{min}(\Delta_{2})$. Alors  $\mathfrak{T}_{1}(D)=\{\Delta_{1}\}$, $\mathfrak{T}_{2}(D)=\{\Delta_{2}\}$. Pour $d=1,2$ et $\Delta'_{d}\in Int(\lambda_{d})$, on a $j_{max}(\Delta'_{d})\not=j_{min}(D)$, $j_{max}(\Delta'_{d})\not=j_{max}(D)$ par comparaison des parit\'es.  D'apr\`es (5), cela entra\^{\i}ne $ \Delta'_{d}\not\in \mathfrak{D}_{d}(D)$. Donc $\mathfrak{D}_{1}(D)=\mathfrak{D}_{2}(D)=\emptyset$.   Si $\Delta_{2}^+$ existe, on a $j_{max}(\Delta_{2}^+)<j_{min}(\Delta_{2})=j_{max}(D)$, d'o\`u
   $\Delta_{2}^+\in \mathfrak{D}_{2}(\geq D)$. Enfin, l'\'egalit\'e $j_{max}(D)=j_{min}(\Delta_{2})$ et la relation $j_{max}(D)\in J(\Delta_{1})$ entra\^{\i}nent $j_{max}(D)\in J^+$. Alors
  $$\delta^+(D)=\tau_{1}(\Delta_{1})+\tau_{2}(\Delta_{2})+r'_{1}+r'_{2}+1,$$
 $$\tau^+(D)=\tau_{1}(\Delta_{1})+\delta_{2}(\Delta_{2}^+)+r'_{2}.$$ 
  Les propri\'et\'es requises sont v\'erifi\'ees.

 (d)  Supposons que $\vert J(D)\vert \geq2$, que $J(D)\subset J(\Delta_{1})$,  que $j_{min}(D)$ soit pair et que  $j_{max}(D)$ soit impair. Comme en (c), on a $j_{max}(D)=j_{min}(\Delta_{2})$. On a $j_{min}(D)=j_{max}(\Delta'_{d})$ pour un $d=1,2$ et un $\Delta'_{d}\in Int(\lambda_{d})$. Si $d=1$, on a $J(\Delta'_{1})\cap J(\Delta_{1})\not=\emptyset $ donc $\Delta'_{1}=\Delta_{1}$. Mais c'est impossible puisque $j_{max}(\Delta_{1})\geq j_{max}(D)>j_{min}(D)$. Donc $d=2$ et forc\'ement $\Delta'_{2}=\Delta_{2}^+$ (ce raisonnement montre que $\Delta_{2}^+$ existe). D'o\`u $j_{min}(D)=j_{max}(\Delta_{2}^+)$. On voit que $\mathfrak{T}_{2}(D)=\{\Delta_{2}\}$ et $\mathfrak{D}_{2}(D)=\{\Delta_{2}^+\}$. Pour $\Delta'_{1}\in \tilde{Int}(\lambda_{1})$, on ne peut avoir $j_{min}(\Delta'_{1})\in J(D)$ ou $j_{max}(\Delta'_{1})\in J(D)$ que si $\Delta'_{1}=\Delta_{1}$. On sait que $j_{min}(\Delta_{1})\leq j_{min}(D)$ et $j_{max}(D)\leq j_{max}(\Delta_{1})$. Par comparaison des parit\'es, ces in\'egalit\'es sont strictes. Donc $j_{min}(\Delta_{1})$ et $j_{max}(\Delta_{1})$ n'appartiennent pas \`a $J(D)$ et, gr\^ace \`a (5),  on conclut $\mathfrak{T}_{1}(D)=\mathfrak{D}_{1}(D)=\emptyset$. Enfin, l'in\'egalit\'e $j_{min}(\Delta_{1})\leq j_{max}(D)$ montre que $\Delta_{1}\in \mathfrak{T}_{1}(\geq D)$. On a les m\^emes formules que dans le cas (c):
   $$\delta^+(D)=\tau_{1}(\Delta_{1})+\tau_{2}(\Delta_{2})+r'_{1}+r'_{2}+1,$$
 $$\tau^+(D)=\tau_{1}(\Delta_{1})+\delta_{2}(\Delta_{2}^+)+r'_{2}.$$ 
  Les propri\'et\'es requises sont v\'erifi\'ees.
  
   (e)  Supposons que $\vert J(D)\vert \geq2$, que $J(D)\subset J(\Delta_{1})$,  que $j_{min}(D)$ soit impair et que  $j_{max}(D)$ soit pair. Comme en (c), on a $j_{min}(D)=j_{min}(\Delta_{1})$. Un raisonnement similaire \`a ceux ci-dessus montre que $j_{max}(D)=j_{max}(\Delta_{1})$. Donc $\mathfrak{T}_{1}(D)=\mathfrak{D}_{1}(D)=\{\Delta_{1}\}$. Si $\Delta_{2}$, resp. $\Delta_{2}^+$,  existe, on a a forc\'ement $j_{max}(D)\leq j_{min}(\Delta_{2})$ et $j_{max}(\Delta_{2}^+)\leq j_{min}(D)$. Ces in\'egalit\'es sont strictes par comparaison des parit\'es. Cela entra\^{\i}ne qu'il n'existe pas de $\Delta'_{2}\in Int(\lambda_{2})$ tel que $j_{min}(\Delta'_{2})$ ou $j_{max}(\Delta'_{2})$ appartiennent \`a $J(D)$. Donc $\mathfrak{T}_{2}(D)=\mathfrak{D}_{2}(D)=\emptyset$. Par contre, si $\Delta_{2}^+$ existe, on a $\Delta_{2}^+\in \mathfrak{D}_{2}(\geq D)$. Puisque  $j_{max}(D)$  est pair, on a $j_{max}(D)\not\in J^+$. On a alors
   $$\delta^+(D)=\delta_{1}(\Delta_{1})+\delta_{2}(\Delta^+_{2}),$$
   $$\tau^+(D)=\tau_{1}(\Delta_{1})+\delta_{2}(\Delta_{2}^+)+r'_{2}.$$
   Les propri\'et\'es requises sont v\'erifi\'ees.

 (f) Supposons que $\vert J(D)\vert \geq2$, que $J(D)\subset J(\Delta_{1})$ et  que $j_{min}(D)$ et   $j_{max}(D)$ soient pairs. En utilisant des r\'esultats extraits de (d) et (e), on a $j_{min}(D)=j_{max}(\Delta_{2}^+)$ et $j_{max}(D)=j_{max}(\Delta_{1})$.  De plus, $j_{min}(\Delta_{1})<j_{min}(D)$ et $j_{max}(D)< j_{min}(\Delta_{2})$ si $\Delta_{2}$ existe. Donc $\mathfrak{T}_{1}(D)=\mathfrak{T}_{2}(D)=\emptyset$, $\mathfrak{D}_{1}(D)=\{\Delta_{1}\}$ et $\mathfrak{D}_{2}(D)=\{\Delta_{2}^+\}$.  On a encore $j_{max}(D)\not\in J^+$. Puisque $j_{min}(\Delta_{1})\leq j_{max}(D)$, on a $\Delta_{1}\in \mathfrak{T}_{1}(\geq D)$.  On a les m\^emes relations que dans le cas (e):
  $$\delta^+(D)=\delta_{1}(\Delta_{1})+\delta_{2}(\Delta^+_{2}),$$
 $$\tau^+(D)=\tau_{1}(\Delta_{1})+\delta_{2}(\Delta_{2}^+)+r'_{2}.$$ 
  
  On a des cas (g), (h), (i), (j) qui sont les sym\'etriques de (c), (d), (e), (f): on remplace la condition $J(D)\subset J(\Delta_{1})$ par $J(D)\subset J(\Delta_{2})$. Les formules que l'on obtient sont les exacts sym\'etriques de celles obtenues dans les cas trait\'es.   
  
  Comme on l'a dit, les formules ci-dessus supposaient $D\not=D_{min}$. Supposons maintenant $D=D_{min}$. On a $D_{1}=D_{1,min}$, $J(D_{min})\subset J(\Delta_{1,min})$ et $D_{2}$ n'existe pas.
  
  (k)  Supposons que $j_{min}(D_{min})$ soit impair. On a alors $j_{min}(D_{min})=j_{min}(\Delta_{1,min})$ comme en (c). On en d\'eduit $\mathfrak{T}_{1}(D_{min})=\{\Delta_{1,min}\}$ mais $\mathfrak{D}_{1}(D_{min})=\emptyset$ ($\Delta_{1,min}$ n'appartient pas \`a $\mathfrak{D}_{1}(D_{min})$ car, par d\'efinition, cet ensemble est un sous-ensemble de $Int(\lambda_{1})$, lequel ne contient pas $\Delta_{1,min}$).  Si $Int(\lambda_{2})\not=\emptyset$, on a $j_{max}(\Delta_{2}^+)>j_{min}(D)$, donc $\mathfrak{T}_{2}(D_{min})=\mathfrak{D}_{2}(D_{min})=\emptyset$. Par contre, $\Delta_{2}^+$ appartient \`a $\mathfrak{D}_{2}(\geq D_{min})$. L'unique formule est 
  $$\tau^+(D_{min})=\tau_{1}(\Delta_{1,min})+\delta_{2}(\Delta_{2}^+)+r'_{2}$$
  et les propri\'et\'es requises sont v\'erifi\'ees.
  
  (l) Supposons que $j_{min}(D_{min})$ soit pair. Alors $j_{min}(D_{min})=j_{max}(\Delta_{2}^+)$ comme  en (d). On voit que $\mathfrak{T}_{1}(D_{min})=\mathfrak{D}_{1}(D_{min})=\mathfrak{T}_{2}(D_{min})=\emptyset$ et  $\mathfrak{D}_{2}(D_{min})=\{\Delta_{2}^+\}$. On a aussi $\Delta_{1,min}\in \mathfrak{T}_{1}(\geq D_{min})$. La formule est la m\^eme que ci-dessus:
  $$\tau^+(D_{min})=\tau_{1}(\Delta_{1,min})+\delta_{2}(\Delta_{2}^+)+r'_{2}$$
  et les propri\'et\'es requises sont v\'erifi\'ees.

  On peut alors   prouver par r\'ecurrence descendante l'assertion suivante: pour $D\in \tilde{Int}_{\lambda_{1},\lambda_{2}}(\lambda)$, il existe pour $d=1,2$ d'uniques fonctions $\tau_{d}$, resp.  $\delta_{d}$, d\'efinies  sur  $\mathfrak{T}_{d}(\geq D)$, resp. $\mathfrak{D}_{d}(\geq D)$, de sorte que les formules  $(X_{r'_{1},r'_{2}})$ soient v\'erifi\'ees pour tout $D'\geq D$. En effet,    soit $D\in Int_{\lambda_{1},\lambda_{2}}(\lambda)$, supposons que l'assertion ci-dessus soit v\'erifi\'ee pour $D^+$ (la condition est vide si $D$ est maximal). Les fonctions $\tau_{d}$ et $\delta_{d}$ sont donc uniquement d\'efinies sur $\mathfrak{T}_{d}(\geq D^+)$, resp. $\mathfrak{D}_{d}(\geq D^+)$. Il faut montrer que l'on peut d\'efinir  d'une seule fa\c{c}on des termes $\tau_{d}(\Delta_{d})$ pour $\Delta_{d}\in \mathfrak{T}_{d}(D)$ et $\delta_{d}(\Delta_{d})$ pour $\Delta_{d}\in \mathfrak{D}_{d}(D)$ de sorte que les formules soient aussi v\'erifi\'ees pour l'intervalle $D$. Par exemple, traitons le cas (a). Le terme $\delta_{2}(\Delta_{2}^+)$ est d\'ej\`a d\'efini. On doit d\'efinir $\tau_{1}(\Delta_{1})$ et $\tau_{2}(\Delta_{2})$ de sorte que  
  $$\delta^+(D)=\tau_{1}(\Delta_{1})+\tau_{2}(\Delta_{2})+r'_{1}+r'_{2}+1,$$
  $$\tau^+(D)=\tau_{1}(\Delta_{1})+\delta_{2}(\Delta_{2}^+)+r'_{2}.$$
  Il est clair que ces \'equations ont une solution et que celle-ci est unique.  Les autres cas (b) \`a (l) sont similaires. L'assertion est donc d\'emontr\'ee par r\'ecurrence. Pour $D=D_{min}$, on obtient l'assertion voulue: pour deux \'el\'ements donn\'es $r'_{1},r'_{2}\in {\mathbb Z}/2{\mathbb Z}$,  il existe d'uniques 
 $(\tau_{1},\delta_{1})\in ({\mathbb Z}/2{\mathbb Z})^{\tilde{Int}(\lambda_{1})}\times ({\mathbb Z}/2{\mathbb Z})^{Int(\lambda_{1})}$, $(\tau_{2},\delta_{2})\in ({\mathbb Z}/2{\mathbb Z})^{Int(\lambda_{2})}\times ({\mathbb Z}/2{\mathbb Z})^{Int(\lambda_{2})}$ tels que   soient v\'erifi\'ees les formules $(X_{r'_{1},r'_{2}})$.

 Ces paires $(\tau_{1},\delta_{1})$ et $(\tau_{2},\delta_{2})$ ne  v\'erifient pas forc\'ement les conditions impos\'ees au d\'ebut de la d\'emonstration. Si $(\tau_{2},\delta_{2})$ est bien un \'el\'ement de ${\cal F}am(\lambda_{2})$, $(\tau_{1},\delta_{1})$ n'est pas forc\'ement un \'el\'ement de ${\cal F}am(\lambda_{1})$: c'en est un si et seulement si $\tau_{1}(\Delta_{1,min})=0$. D'autre part, en admettant que cette condition soit v\'erifi\'ee, nos paires v\'erifient les conditions requises si et seulement si $r'_{1}\equiv r(\tau_{1},\delta_{1})\,\,mod\,\,2{\mathbb Z}$ et $r'_{2}\equiv r(\tau_{2},\delta_{2})\,\,mod\,\,2{\mathbb Z}$.  Pour d\'emontrer la premi\`ere assertion du lemme, il suffit de prouver que ces conditions sont v\'erifi\'ees pour un seul couple $(r'_{1},r'_{2})$. 
 
 Continuons avec un couple quelconque $(r'_{1},r'_{2})$ et les paires $(\tau_{1},\delta_{1})$ et $(\tau_{2},\delta_{2})$ que l'on a construites ci-dessus. Posons $a=\tau_{1}(\Delta_{1,min})$. D\'efinissons $\underline{\tau}_{1}$ par $\underline{\tau}_{1}(\Delta)=\tau_{1}(\Delta)+a$. Alors $(\underline{\tau}_{1},\delta_{1})$ appartient bien \`a ${\cal F}am(\lambda_{1})$. On pose $\underline{r}_{1}=r(\underline{\tau}_{1},\delta_{1})$, $\underline{r}_{2}=r(\tau_{2},\delta_{2})$. Les conditions \`a v\'erifier sont 
 
  (6) $a=0$, $\underline{r}_{1}\equiv r'_{1}\,\,mod\,\,2{\mathbb Z}$, $\underline{r}_{2}\equiv r'_{2}\,\,mod\,\,2{\mathbb Z}$.
  
  Remarquons que la premi\`ere condition est redondante avec la troisi\`eme. En effet, comme on l'a vu dans la preuve de 3.1(2), on a par construction
  $$\tau^+(D_{min})=\tau_{1}(\Delta_{1,min})+\delta_{2}(\Delta_{2}(D_{min})^+)+r'_{2}.$$
  On sait que 
  $\delta_{2}(\Delta_{2}(D_{min})^+)=\underline{r}_{2}$, cf. 2.4(2). On a aussi $\tau^+(D_{min})=0$ par l'hypoth\`ese (2), d'o\`u $a+r'_{2}+\underline{r}_{2}\equiv 0\,\,mod\,\,2{\mathbb Z}$.
  
  Construisons les fonctions associ\'ees \`a $\underline{\iota}_{1}=(\underline{\tau}_{1},\delta_{1})$ et $\underline{\iota}_{2}=(\tau_{2},\delta_{2})$, que l'on note $\underline{\tau}^{\zeta}=\tau^{\zeta}_{\underline{\iota}_{1},\underline{\iota}_{2}}$ et $\underline{\delta}^{\zeta}=\delta^{\zeta}_{\underline{\iota}_{1},\underline{\iota}_{2}}$. Cela revient, dans la construction des fonctions $\tau^{\zeta}$ et $\delta^{\zeta}$ par les formules $(X_{r'_{1},r'_{2}})$, \`a changer $\tau_{1}$ en $\underline{\tau}_{1}$, $r'_{1}$ en $\underline{r}_{1}$ et $r'_{2}$ en $\underline{r}_{2}$.  On remarque que les  termes $\tau_{1}(\Delta_{1})$ et $r'_{2}$ n'interviennent que par leur somme $\tau_{1}(\Delta_{1})+r'_{2}$. Or, comme on vient de le voir, $\underline{\tau}_{1}(\Delta_{1})+\underline{r}_{2}=\tau_{1}(\Delta_{1})+a+\underline{r}_{2}=\tau_{1}(\Delta_{1})+r'_{2}$.  Changer $\tau_{1}$ en $\underline{\tau}_{1}$ et $r'_{2}$ en $\underline{r}_{2}$ ne change donc pas les fonctions $\tau^{\zeta}$ et $\delta^{\zeta}$. On remarque que $r'_{1}$ intervient exactement dans les expressions $\delta^{\zeta}(D)$ ou $\tau^{\zeta}(D)$ telles que $\delta^{-\zeta}(D)=\delta^{\zeta}(D)+1$ ou $\tau^{-\zeta}(D)=\tau^{\zeta}(D)+1$. Changer $r'_{1}$ en $\underline{r}_{1}$ change donc les fonctions $\tau^{\zeta}$ et $\delta^{\zeta}$ en multipliant \'eventuellement $\zeta$ par $-1$, en identifiant le signe $\zeta$ \`a un \'el\'ement de $\{\pm 1\}$. Pr\'ecis\'ement, posons $u=(-1)^{r'_{1}+\underline{r}_{1}}$. On obtient les \'egalit\'es
  $$\underline{\tau}^{\zeta}=\tau^{u\zeta},\,\,\underline{\delta}^{\zeta}=\delta^{u\zeta}.$$
  En posant $\underline{C}^{\zeta}=C^{\zeta}_{\underline{\iota}_{1},\underline{\iota}_{2}}$, ces \'egalit\'es entra\^{\i}nent  $\underline{C}^{\zeta}=C^{u\zeta}$. D'apr\`es 3.3(3),
on a les \'egalit\'es
$$\underline{C}^{\zeta}=\left\lbrace\begin{array}{cc}2(\underline{r}_{1}+\zeta \underline{r}_{2}),&\,\,{\text si}\,\,\underline{r}_{1}+\underline{r}_{2}\text{\,\,est\,\,pair},\\ -2(\underline{r}_{1}+\zeta \underline{r}_{2}+1),&\,\,{\text si}\,\,\underline{r}_{1}+\underline{r}_{2}\text{\,\,est\,\,impair}.\\ \end{array}\right.$$
 Par l'hypoth\`ese (3), on a aussi
$$C^{u\zeta}=\left\lbrace\begin{array}{cc}2(r_{1}+u\zeta r_{2}),&\,\,{\text si}\,\,r_{1}+r_{2}\text{\,\,est\,\,pair},\\ -2(r_{1}+u\zeta r_{2}+1),&\,\,{\text si}\,\,r_{1}+r_{2}\text{\,\,est\,\,impair}.\\ \end{array}\right.$$ 
 L'\'egalit\'e  de ces deux expressions est  \'equivalente aux \'egalit\'es suivantes
 
 si $\underline{r}_{1}+\underline{r}_{2}$ et $r_{1}+r_{2}$ sont pairs, $\underline{r}_{1}+\zeta\underline{r}_{2}=r_{1}+u\zeta r_{2}$ pour $\zeta=\pm 1$;
 
 si $\underline{r}_{1}+\underline{r}_{2}$ est pair et $r_{1}+r_{2}$ est impair, $\underline{r}_{1}+\zeta\underline{r}_{2}=-(r_{1}+u\zeta r_{2}+1)$ pour $\zeta=\pm 1$;
 
 si $\underline{r}_{1}+\underline{r}_{2}$ est impair et $r_{1}+r_{2}$ est pair, $-(\underline{r}_{1}+\zeta\underline{r}_{2}+1)=r_{1}+u\zeta r_{2}$ pour $\zeta=\pm 1$;
 
  si $\underline{r}_{1}+\underline{r}_{2}$ et $r_{1}+r_{2}$ sont impairs, $-(\underline{r}_{1}+\zeta\underline{r}_{2}+1)=-(r_{1}+u\zeta r_{2}+1)$ pour $\zeta=\pm 1$.
  
  En sommant en $\zeta=\pm 1$, le deuxi\`eme cas entra\^{\i}ne $\underline{r}_{1}=-(r_{1}+1)$. C'est impossible puisque $\underline{r}_{1}$ et $r_{1}$ sont positifs ou nuls. Ce cas ne se produit donc pas. Le troisi\`eme cas non plus, pour la m\^eme raison.  Cela montre que $\underline{r}_{1}+\underline{r}_{2}$ et $r_{1}+r_{2}$ sont de la m\^eme parit\'e. Dans ce cas, les \'egalit\'es ci-dessus  entra\^{\i}nent $\underline{r}_{1}=r_{1}$ et $\underline{r}_{2}=u r_{2}$. Alors les conditions (6) sont v\'erifi\'ees si et seulement si $r'_{1}\equiv r_{1}\,\,mod\,\,2{\mathbb Z}$ et $r'_{2}\equiv r_{2}\,\, mod\,\,2{\mathbb Z}$. Cela d\'emontre la premi\`ere assertion du lemme. Pour ce couple $(r'_{1},r'_{2})$ ainsi d\'etermin\'e, on vient de voir que $\underline{r}_{1}=r_{1}$. On a aussi  $u=(-1)^{r'_{1}+\underline{r}_{1}}=1$, donc $\underline{r}_{2}=u r_{2}=r_{2}$. Cela d\'emontre la seconde assertion de l'\'enonc\'e. 
  $\square$

  \bigskip
 
 \section{Le front d'onde de $\pi(\lambda^+,\epsilon^+,\lambda^-,\epsilon^-)$}

 \bigskip
 \subsection{Le r\'esultat de \cite{W5}}
 Soit $m\in {\mathbb N}$ et $(\lambda,\epsilon)\in \boldsymbol{{\cal P}^{symp}}(2m)$. On a introduit en 1.3 la repr\'esentation $\boldsymbol{\rho}_{\lambda,\epsilon}$ de $W_{N_{\lambda,\epsilon}}$. On  sait qu'elle se d\'ecompose en
 $$\boldsymbol{\rho}_{\lambda,\epsilon}=\oplus_{(\lambda',\epsilon')}mult(\lambda,\epsilon;\lambda',\epsilon')\rho_{\lambda',\epsilon'},$$
 o\`u $(\lambda',\epsilon')$ parcourt les \'el\'ements de $\boldsymbol{{\cal P}^{symp}}(2m)$ tels que $k_{\lambda',\epsilon'}=k_{\lambda,\epsilon}$ et les $mult(\lambda,\epsilon;\lambda',\epsilon')$ sont des entiers positifs ou  nuls. Le couple $(\lambda,\epsilon)$ est minimal dans cette d\'ecomposition, c'est-\`a-dire que l'on a

 si $mult(\lambda,\epsilon;\lambda',\epsilon')\not=0$, alors $\lambda'>\lambda$ ou $(\lambda',\epsilon')=(\lambda,\epsilon)$.
 
 \noindent De plus $mult(\lambda,\epsilon;\lambda,\epsilon)=1$.
 
 Pour tout couple $(\mu,\nu)\in \boldsymbol{{\cal P}^{symp}}(2m)$, notons $({^s\mu},{^s\nu})$ le couple tel que $k_{{^s\mu},{^s\nu}}=k_{\mu,\nu}$ et $\rho_{{^s\mu},{^s\nu}}=\rho_{\mu,\nu}\otimes sgn$.
 
 \ass{Proposition}{Supposons que tous les termes de $\lambda$ soient pairs. Alors il existe un  unique couple $(\lambda^{min},\epsilon^{min})\in \boldsymbol{{\cal P}^{symp}}(2m)$ v\'erifiant les propri\'et\'es suivantes:
 
 (i) $mult(\lambda,\epsilon;{^s\lambda}^{min},{^s\epsilon}^{min})=1$;
 
 (ii) pour tout \'el\'ement $(\lambda',\epsilon')\in \boldsymbol{{\cal P}^{symp}}(2m)$ tel que $mult(\lambda,\epsilon;{^s\lambda}',{^s\epsilon}')\not=0$, on a $\lambda^{min}<\lambda'$ ou $(\lambda',\epsilon')=(\lambda^{min},\epsilon^{min})$.}
 Cf. \cite{W5} th\'eor\`eme 5.7.

 \bigskip
 
  \subsection{Calcul de $M_{\pi}(\mu_{1},\eta_{1};\mu_{2},\eta_{2})$}
  
   On fixe d\'esormais un quadruplet $(\lambda^+,\epsilon^+,\lambda^-,\epsilon^-)\in \mathfrak{Irr}_{quad}^{bp}(2n)$. Rappelons que l'exposant $bp$ signifie que tous les termes de $\lambda^+$ et $\lambda^-$ sont pairs. On pose
 $$\pi=\pi(\lambda^+,\epsilon^+,\lambda^-,\epsilon^-)$$
 et on note $\sharp$ l'indice $iso$ ou $an $ tel que $\pi\in Irr_{tunip,\sharp}$.

 Soient $n_{1},n_{2}\in {\mathbb N}$ tels que $n_{1}+n_{2}=n$.  Soient $(\mu_{1},\eta_{1})\in \boldsymbol{{\cal P}^{orth}}(2n_{1}+1)_{k=1}$ et $(\mu_{2},\eta_{2})\in \boldsymbol{{\cal P}^{orth}}(2n_{2})_{k=0}$. On a d\'efini le nombre $M_{\pi}(\mu_{1},\eta_{1};\mu_{2},\eta_{2})$ en 1.4. On se propose de le calculer. 
 
 Le couple $(0,\rho_{\mu_{1},\eta_{1}})$ appartient \`a $\Sigma_{n_{1},imp}$ et son symbole appartient \`a $Fam(sp(\mu_{1},\eta_{1}))$ pour une partition sp\'eciale $sp(\mu_{1},\eta_{1})\in {\cal P}^{orth,sp}(2n_{1}+1)$. Posons $\lambda_{1}=d(sp(\mu_{1},\eta_{1}))$. On a $\lambda_{1}\in {\cal P}^{symp,sp}(2n_{1})$. Il r\'esulte de 2.6 que le symbole $\Lambda_{1}$ de $(0,\rho_{\mu_{1},\eta_{1}}\otimes sgn)$ appartient \`a $Fam(\lambda_{1})$. 
 
 Pour $\xi=\pm$, le couple $(0,\rho_{\mu_{2},\eta_{2}}^{\xi})$ appartient \`a $\Sigma_{n_{2},pair}$ et son symbole appartient \`a $Fam(sp(\mu_{2},\eta_{2}))$ pour une partition sp\'eciale $sp(\mu_{2},\eta_{2})\in {\cal P}^{orth,sp}(2n_{2})$. Celle-ci ne d\'epend pas du signe $\xi$: changer de signe revient \`a \'echanger les deux termes $X$ et $Y$ du symbole. Posons $\lambda_{2}=d(sp(\mu_{2},\eta_{2}))$. On a $\lambda_{2}\in {\cal P}^{orth,sp}(2n_{2})$. Le symbole  $\Lambda_{2}^{\xi}$ de $(0,\rho_{\mu_{2},\eta_{2}}^{\xi}\otimes sgn)$ appartient \`a $Fam(\lambda_{2})$.
 
 Signalons que l'on a les in\'egalit\'es
 
 (1) $\mu_{1}\leq sp(\mu_{1},\eta_{1})$, $\mu_{2}\leq sp(\mu_{2},\eta_{2})$,
 
 \noindent cf. \cite{W4} lemmes 1.4 et 1.5. 
 
Posons $\gamma_{0}=(0,0,n_{1},n_{2})$.  Par d\'efinition de la multiplicit\'e $m_{\pi}(  \rho_{\mu_{1},\eta_{1}}\otimes sgn,\rho_{\mu_{2},\eta_{2}}^{\xi}\otimes sgn)$ et d'apr\`es 1.5(4), cette multiplicit\'e est celle de $( \rho_{\mu_{1},\eta_{1}}\otimes sgn)\otimes (\rho_{\mu_{2},\eta_{2}}^{\xi}\otimes sgn)$ dans la composante dans ${\cal R}(\gamma_{0})$ de
$$\kappa_{\pi}={\cal F}( \Pi),$$
o\`u on a pos\'e
$$\Pi=\rho\iota((\boldsymbol{\rho}_{\lambda^+,\epsilon^+}\otimes sgn)\otimes 
 (\boldsymbol{\rho}_{\lambda^+,\epsilon^+}\otimes sgn)).$$
 En 1.3, on a associ\'e \`a  $(\lambda^+,\epsilon^+,\lambda^-,\epsilon^-)$ un \'el\'ement $\gamma=(r',r'',N^+,N^-)\in \Gamma$ et identifi\'e $(\boldsymbol{\rho}_{\lambda^+,\epsilon^+}\otimes sgn)\otimes 
 (\boldsymbol{\rho}_{\lambda^+,\epsilon^+}\otimes sgn)$ \`a un \'el\'ement de ${\cal R}(\gamma)$.  On pose $r_{1}=r'$, $r_{2}=(-1)^{r'}r''$. Par construction de $\rho\iota$, l'\'el\'ement 
 $ \Pi$
  n'a de composante non nulle que dans les composantes ${\cal R}(\gamma')$ pour $\gamma'$ de la forme $(r_{1},r_{2},N_{1},N_{2})$. Par d\'efinition de ${\cal F}$, pour un tel $\gamma'$ et pour $\varphi\in {\cal R}(\gamma')$,  l'\'el\'ement ${\cal F}(\varphi)$ n'a de composante non nulle dans ${\cal R}(\gamma_{0})$ que si $N_{1}+r_{1}^2+r_{1}=n_{1}$ et $N_{2}+r_{2}^2=n_{2}$. Cela entra\^{\i}ne
  
  (2) si $n_{1}<r_{1}^2+r_{1}$ ou $n_{2}< r_{2}^2$, on a $M_{\pi}(\mu_{1},\eta_{1};\mu_{2},\eta_{2})=0$.
  
  Supposons 
  
  (3) $n_{1}\geq r_{1}^2+r_{1}$ et $n_{2}\geq r_{2}^2$. 
  
  Posons $N_{1}=n_{1}-r_{1}^2-r_{1}$, $N_{2}=n_{2}-r_{2}^2$ et
  $\underline{\gamma}=(r_{1},r_{2},N_{1},N_{2})$. On peut se limiter \`a consid\'erer  la composante $\Pi_{\underline{\gamma}}$ de $\Pi$ dans ${\cal R}(\underline{\gamma})$. Plus pr\'ecis\'ement, pour $d=1,2$, notons ${\cal F}am_{r_{d}}(\lambda_{d})$ l'ensemble les $\iota_{d}=(\tau_{d},\delta_{d})\in {\cal F}am(\lambda_{d})$ tels que $r(\tau_{d},\delta_{d})=r_{d}$. Pour de tels \'el\'ements, notons $(r_{d},\rho_{\iota_{d}})$ l'\'el\'ement de $\Sigma_{n_{1},imp}$ si $d=1$ et $\Sigma_{n_{2},pair}$ si $d=2$ associ\'e \`a $\iota_{d}$. Posons $\Lambda_{\iota_{d}}=symb(r_{d},\rho_{\iota_{d}})$.  Notons $m(\Pi_{\underline{\gamma}},\rho_{\iota_{1}}\otimes \rho_{\iota_{2}})$ la multiplicit\'e de $\rho_{\iota_{1}}\otimes \rho_{\iota_{2}}$ dans $\Pi_{\underline{\gamma}}$. Alors, par d\'efinition de ${\cal F}$, on a l'\'egalit\'e
$$(4) \qquad m_{\pi}(\rho_{\mu_{1},\eta_{1}}\otimes sgn,\rho^{\xi}_{\mu_{2},\delta_{2}}\otimes sgn)=\vert Fam(\lambda_{1})\vert ^{-1/2}\vert Fam(\lambda_{2})\vert ^{-1/2}$$
$$\sum_{\iota_{1}\in {\cal F}am_{r_{1}}(\lambda_{1}),\iota_{2}\in {\cal F}am_{r_{2}}(\lambda_{2})}(-1)^{<\Lambda_{1},\Lambda_{\iota_{1}}>+<\Lambda^{\xi}_{2},\Lambda_{\iota_{2}}>}m(\Pi_{\underline{\gamma}},\rho_{\iota_{1}}\otimes \rho_{\iota_{2}}).$$  
 
 Pour $\zeta=\pm$, on pose $n^{\zeta}=S(\lambda^{\zeta})/2$,  $k^{\zeta}=k_{\lambda^{\zeta},\epsilon^{\zeta}}$. Notons  $\boldsymbol{{\cal P}^{symp}}(2n^{\zeta})_{k^{\zeta}}$  l'ensemble des $(\lambda',\epsilon')\in \boldsymbol{{\cal P}^{symp}}(2n^{\zeta})$ tels que $k_{\lambda',\epsilon'}=k^{\zeta}$.  On peut \'ecrire
 $$  (\boldsymbol{\rho}_{\lambda^+,\epsilon^+}\otimes sgn)\otimes 
 (\boldsymbol{\rho}_{\lambda^+,\epsilon^+}\otimes sgn)=$$
 $$\sum_{(\lambda^{'+},\epsilon^{'+})\in \boldsymbol{{\cal P}^{symp}}(2n^+)_{k^+},(\lambda^{'-},\epsilon^{'-})\in \boldsymbol{{\cal P}^{symp}}(2n^-)_{k^-}}x(\lambda^{'+},\epsilon^{'+},\lambda^{'-},\epsilon^{'-})\rho_{\lambda^{'+},\epsilon^{'+}}\otimes \rho_{\lambda^{'-},\epsilon^{'-}},$$
 o\`u les $x(\lambda^{'+},\epsilon^{'+},\lambda^{'-},\epsilon^{'-})$ sont des multiplicit\'es. Pr\'ecis\'ement, avec les notations de 4.1, on a
 $$(5) \qquad x(\lambda^{'+},\epsilon^{'+},\lambda^{'-},\epsilon^{'-})=mult(\lambda^+,\epsilon^+;{^s\lambda}^{'+},{^s\epsilon}^{'+})mult(\lambda^-,\epsilon^-;{^s\lambda}^{'-},{^s\epsilon}^{'-}).$$
  Pour $(\lambda^{'+},\epsilon^{'+})\in \boldsymbol{{\cal P}^{symp}}(2n^+)_{k^+}$ et $(\lambda^{'-},\epsilon^{'-})\in \boldsymbol{{\cal P}^{symp}}(2n^-)_{k^-}$, notons $\Pi_{\underline{\gamma}}(\lambda^{'+},\epsilon^{'+},\lambda^{'-},\epsilon^{'-})$ la composante dans ${\cal R}(\underline{\gamma})$ de
 $$\rho\iota(\rho_{\lambda^{'+},\epsilon^{'+}}\otimes \rho_{\lambda^{'-},\epsilon^{'-}}).$$
 Pour $\iota_{1}\in {\cal F}am_{r_{1}}(\lambda_{1})$ et $\iota_{2}\in {\cal F}am_{r_{2}}(\lambda_{2})$, notons $m(\Pi_{\underline{\gamma}}(\lambda^{'+},\epsilon^{'+},\lambda^{'-},\epsilon^{'-}),\rho_{\iota_{1}}\otimes \rho_{\iota_{2}})$ la multiplicit\'e de $\rho_{\iota_{1}}\otimes \rho_{\iota_{2}}$ dans $\Pi_{\underline{\gamma}}(\lambda^{'+},\epsilon^{'+},\lambda^{'-},\epsilon^{'-})$. On a
  $$m(\Pi_{\underline{\gamma}},\rho_{\iota_{1}}\otimes \rho_{\iota_{2}})=\sum_{(\lambda^{'+},\epsilon^{'+})\in \boldsymbol{{\cal P}^{symp}}(2n^+)_{k^+},(\lambda^{'-},\epsilon^{'-})\in \boldsymbol{{\cal P}^{symp}}(2n^-)_{k^-}}x(\lambda^{'+},\epsilon^{'+},\lambda^{'-},\epsilon^{'-})$$
  $$m(\Pi_{\underline{\gamma}}(\lambda^{'+},\epsilon^{'+},\lambda^{'-},\epsilon^{'-}),\rho_{\iota_{1}}\otimes \rho_{\iota_{2}}).$$
  En vertu de la d\'efinition pos\'ee en 1.4, on d\'eduit de (4) la formule finale
  $$(6) \qquad M_{\pi}(\mu_{1},\eta_{1};\mu_{2},\eta_{2})=\vert Fam(\lambda_{1})\vert ^{-1/2}\vert Fam(\lambda_{2})\vert ^{-1/2}$$
$$\sum_{\iota_{1}\in {\cal F}am_{r_{1}}(\lambda_{1}),\iota_{2}\in {\cal F}am_{r_{2}}(\lambda_{2})}
 (-1)^{<\Lambda_{1},\Lambda_{\iota_{1}}>}\left((-1)^{<\Lambda^{+}_{2},\Lambda_{\iota_{2}}>}+sgn_{\sharp}(-1)^{<\Lambda_{2}^-,\Lambda_{\iota_{2}}>}\right)$$
 $$\sum_{(\lambda^{'+},\epsilon^{'+})\in \boldsymbol{{\cal P}^{symp}}(2n^+)_{k^+},(\lambda^{'-},\epsilon^{'-})\in \boldsymbol{{\cal P}^{symp}}(2n^-)_{k^-}}x(\lambda^{'+},\epsilon^{'+},\lambda^{'-},\epsilon^{'-})m(\Pi_{\underline{\gamma}}(\lambda^{'+},\epsilon^{'+},\lambda^{'-},\epsilon^{'-}),\rho_{\iota_{1}}\otimes \rho_{\iota_{2}}).$$

 \bigskip
 
 \subsection{Comparaison entre  deux constructions}
 On conserve les notations du paragraphe pr\'ec\'edent et on impose l'hypoth\`ese (3) de ce paragraphe. Consid\'erons des \'el\'ements $\iota_{1}\in {\cal F}am_{r_{1}}(\lambda_{1})$, $\iota_{2}\in {\cal F}am_{r_{2}}(\lambda_{2})$, $(\lambda^{'+},\epsilon^{'+})\in \boldsymbol{{\cal P}^{symp}}(2n^+)_{k^+}$, $(\lambda^{'-},\epsilon^{'-})\in \boldsymbol{{\cal P}^{symp}}(2n^-)_{k^-}$. On a d\'efini la multiplicit\'e
  $$m(\Pi_{\underline{\gamma}}(\lambda^{'+},\epsilon^{'+},\lambda^{'-},\epsilon^{'-}),\rho_{\iota_{1}}\otimes \rho_{\iota_{2}}).$$ 
 
  Un \'el\'ement  $\zeta\in \{\pm 1\}$ \'etant fix\'e, on a  associ\'e  en 3.4 \`a  $(r_{1},r_{2})$ un couple $(h^+,h^-)$. En se reportant \`a la d\'efinition de 1.2 et en se rappelant que $(r_{1},r_{2})=(r',(-1)^{r'}r'')$, on   v\'erifie cas par cas qu'il est \'egal \`a $(k^+,k^-)$ pourvu que $\zeta=1$ si $k^+>k^-$, $\zeta=-1$ si $k^+<k^-$. Notons que $k^+>k^-$ \'equivaut \`a $(-1)^{r_{1}}r_{2}>0$ et $k^+<k^-$ \'equivaut \`a $(-1)^{r_{1}}r_{2}<0$.  Si $k^+=k^-$, ce qui \'equivaut \`a $r_{2}=0$, $\zeta$ est indiff\'erent, le couple $(h^+,h^-)$ ne d\'ependant pas de $\zeta$ et \'etant \'egal \`a $(k^+,k^-)$. On suppose que $\zeta$ v\'erifie ces conditions.   On peut donc appliquer la construction de 3.4 aux entiers $n^+$ et $n^-$. On en d\'eduit une repr\'esentation $\Pi^{\zeta}(\iota_{1},\iota_{2})$ de $W_{N^+}\otimes W_{N^-}$. On note $m(\Pi^{\zeta}(\iota_{1},\iota_{2}),\rho_{\lambda^{'+},\epsilon^{'+}}\otimes \rho_{\lambda^{'-},\epsilon^{'-}})$ la multiplicit\'e de $\rho_{\lambda^{'+},\epsilon^{'+}}\otimes \rho_{\lambda^{'-},\epsilon^{'-}}$ dans $\Pi^{\zeta}(\iota_{1},\iota_{2})$. Un jeu habituel avec les restrictions et inductions montre que cette multiplicit\'e est \'egale \`a celle de $\rho_{1}\otimes \rho_{2}$ dans la repr\'esentation
 $$\sum_{{\bf N}\in {\cal N}}ind^{W_{N_{1}}\times W_{N_{2}}}_{W_{{\bf N}}}\left(sgn_{CD}^{{\bf a}}\otimes res^{W_{N^+}\times W_{N^-}}_{W_{{\bf N}}}( \rho_{\lambda^{'+},\epsilon^{'+}}\otimes \rho_{\lambda^{'-},\epsilon^{'-}})\right),$$
 o\`u ${\bf a}$ est d\'efini comme en 3.4. Un calcul cas par cas montre que ce ${\bf a}$ est le m\^eme qu'en 1.2, pouvu que, dans le cas $r_{2}=0$, on choisisse $\zeta=1$ si $r_{1}$ est pair, $\zeta=-1$ si $r_{1}$ est impair.  Le signe $\zeta$ \'etant ainsi d\'etermin\'e en tout cas, la repr\'esentation ci-dessus n'est autre que la composante dans ${\cal R}(\underline{\gamma})$ de
 $$\rho\iota(\rho_{\lambda^{'+},\epsilon^{'+}}\otimes \rho_{\lambda^{'-},\epsilon^{'-}}).$$
 On conclut
 $$(1) \qquad  m(\Pi_{\underline{\gamma}}(\lambda^{'+},\epsilon^{'+},\lambda^{'-},\epsilon^{'-}),\rho_{1}\otimes \rho_{2})=m(\Pi^{\zeta}(\iota_{1},\iota_{2}),\rho_{\lambda^{'+},\epsilon^{'+}}\otimes \rho_{\lambda^{'-},\epsilon^{'-}}).$$
Dans les formules 3.4 (4) et (5) intervient le signe $\zeta\nu$, ou $\nu=1$ si $r_{2}\geq0$, $\nu=-1$ si $r_{2}<0$. Avec la d\'efinition de $\zeta$ ci-dessus, on a

(2) $\zeta\nu=(-1)^{r_{1}}$.

\bigskip

\subsection{D\'emonstration du (i) de la proposition 1.4}
On consid\`ere les donn\'ees de 4.2 et on suppose $M_{\pi}(\mu_{1},\eta_{1};\mu_{2},\eta_{2})\not=0$. La relation 4.2(2) entra\^{\i}ne que l'hypoth\`ese 4.2(3) est v\'erifi\'ee. D'apr\`es 4.2(6), on peut fixer des \'el\'ements $\iota_{1}\in {\cal F}am_{r_{1}}(\lambda_{1})$, $\iota_{2}\in {\cal F}am_{r_{2}}(\lambda_{2})$, $(\lambda^{'+},\epsilon^{'+})\in \boldsymbol{{\cal P}^{symp}}(2n^+)_{k^+}$, $(\lambda^{'-},\epsilon^{'-})\in \boldsymbol{{\cal P}^{symp}}(2n^-)_{k^-}$ v\'erifiant les conditions

(1) $x(\lambda^{'+},\epsilon^{'+},\lambda^{'-},\epsilon^{'-})\not=0$;

(2) $m(\Pi_{\underline{\gamma}}(\lambda^{'+},\epsilon^{'+},\lambda^{'-},\epsilon^{'-}),\rho_{\iota_{1}}\otimes \rho_{\iota_{2}})\not=0$.

En vertu de la d\'efinition 4.2(5) de $x(\lambda^{'+},\epsilon^{'+},\lambda^{'-},\epsilon^{'-})$ et de la proposition 4.1, la relation (1) entra\^{\i}ne

(3)  $\lambda^{+,min}\leq \lambda^{'+}$, $\lambda^{-,min}\leq \lambda^{'-}$.

Notons $\lambda$ l'induite endoscopique de $\lambda_{1}$ et $\lambda_{2}$. En vertu de 4.3(1) et de la proposition 3.4(i), la relation (2) entra\^{\i}ne

 (4) $\lambda^{'+}\cup \lambda^{'-}\leq \lambda$. 

De ces deux in\'egalit\'es, on d\'eduit
$$\lambda^{+,min}\cup \lambda^{-,min}\leq \lambda.$$
Posons
$$\mu=d(\lambda^{+,min}\cup \lambda^{-,min}).$$
La dualit\'e est une application d\'ecroissante. L'in\'egalit\'e pr\'ec\'edente entra\^{\i}ne
$d(\lambda)\leq \mu$.
D'apr\`es \cite{W4} proposition 1.9, on a aussi
$d(\lambda_{1})\cup d(\lambda_{2})\leq d(\lambda)$,
d'o\`u $d(\lambda_{1})\cup d(\lambda_{2})\leq \mu$.
Par construction, $d(\lambda_{1})=sp(\mu_{1},\eta_{1})$, $d(\lambda_{2})=sp(\mu_{2},\eta_{2})$. D'o\`u
 $sp(\mu_{1},\eta_{1})\cup sp(\mu_{2},\eta_{2})\leq \mu$.
 En appliquant 4.2(1), on obtient
 $$\mu_{1}\cup \mu_{2}\leq \mu.$$
 C'est l'assertion (i) de la proposition 1.4.
 
 \bigskip
 
 \subsection{D\'emonstration du (ii) de la proposition 3.4}
 La seule donn\'ee est ici le quadruplet $(\lambda^+,\epsilon^+,\lambda^-,\epsilon^-)\in \mathfrak{Irr}_{quad}^{bp}(2n)$. On pose $\lambda=\lambda^{+,min}\cup \lambda^{-,min}$. Fixons une fonction $\chi:Jord^{bp}(\lambda)\cup\{0\}\to {\mathbb Z}/2{\mathbb Z}$ v\'erifiant les conditions suivantes:
 
 $\chi(i)=0$ pour tout $i\in Jord^{bp}(\lambda)$ tel que $mult_{\lambda}(i)=1$;
 
 $\chi(0)=0$;
 
 pour tout $i\in Jord^{bp}(\lambda)$ tel que $mult_{\lambda^{+,min}}(i)\geq1$ et $mult_{\lambda^{-,min}}(i)\geq1$ (ce qui implique $mult_{\lambda}(i)\geq2$), $\chi(i)=0$ si $\epsilon^{+,min}_{i}\not=\epsilon^{-,min}_{i}$ et $\chi(i)=1$ si $\epsilon^{+,min}_{i}
=\epsilon^{-,min}_{i}$;

pour tout $i\in Jord^{bp}(\lambda)$ tel que $mult_{\lambda}(i)\geq2$ et que $mult_{\lambda^{+,min}}(i)=0$ ou $mult_{\lambda^{-,min}}(i)=0$, $\chi(i)=1$.

 On choisit $n_{1},n_{2}$ et $\lambda_{1},\lambda_{2}$ v\'erifiant les conditions de la proposition 3.2, pour ce choix   de la fonction $\chi$. C'est-\`a-dire que $\lambda_{1}\in {\cal P}^{symp,sp}(2n_{1})$, $\lambda_{2}\in {\cal P}^{orth,sp}(2n_{2})$, $\lambda_{1}$ et $\lambda_{2}$ induisent r\'eguli\`erement $\lambda$, $d(\lambda_{1})\cup d(\lambda_{2})=d(\lambda)=\mu$ et $\chi_{\lambda_{1},\lambda_{2}}=\chi$. On pose $\mu_{1}=d(\lambda_{1})$, $\mu_{2}=d(\lambda_{2})$. On a $\mu_{1}\in {\cal P}^{orth,sp}(2n_{1}+1)$, $\mu_{2}\in {\cal P}^{orth,sp}(2n_{2})$, et
 $$(1) \qquad \mu_{1}\cup \mu_{2}=\mu.$$
 
 On d\'efinit $r_{1}$ et $r_{2}$ comme en 4.2: $r_{1}=r'$, $r_{2}=(-1)^{r'}r''$, o\`u $r'$ et $r''$ sont d\'efinis en 1.3. Pour une partition $\nu$ et pour $i\in {\mathbb N}-\{0\}$, posons $mult_{\nu}(\geq i)=\sum_{i'\geq i}mult_{\nu}(i')$. Posons $\eta=(-1)^{r'}$. Pour $\zeta=\pm$, on d\'efinit une fonction $\delta^{\zeta}:Jord^{bp}(\lambda)\to {\mathbb Z}/2{\mathbb Z}$ par
 $$\delta^{\zeta}(i)\equiv mult_{\lambda^{\zeta\eta,min}}(\geq i)\,\mod\,\,2{\mathbb Z}.$$
 On d\'efinit une fonction $\tau^{\zeta}:Jord^{bp}(\lambda)\cup \{0\}\to {\mathbb Z}/2{\mathbb Z}$ par
 
 si $i\not=0$ et $mult_{\lambda^{\zeta\eta,min}}(i)>0$, $\epsilon^{\zeta \eta,min}_{i}=(-1)^{\tau^{\zeta}(i)}$;
 
 si $i\not=0$ et $mult_{\lambda^{\zeta\eta,min}}(i)=0$ (auquel cas $mult_{\lambda^{-\zeta\eta,min}}(i)>0$),  
 $\epsilon^{-\zeta \eta,min}_{i}=(-1)^{\tau^{\zeta}(i)}$;
 
 $\tau^{\zeta}(0)=0$.
 
 On peut consid\'erer que ces fonctions sont d\'efinies sur $Int_{\lambda_{1},\lambda_{2}}(\lambda)$, resp. $\tilde{Int}_{\lambda_{1},\lambda_{2}}(\lambda)$, puisque 
  $\lambda_{1}$ et $\lambda_{2}$ induisent r\'eguli\`erement $\lambda$.  Montrons que
 
 (2) ces fonctions v\'erifient les conditions de 3.5.
 
 Preuve. Soit $i\in Jord^{bp}(\lambda)$. D'apr\`es la d\'efinition ci-dessus, $\delta^-(i)=\delta^+(i)+1$ si et seulement si $mult_{\lambda}(i)$ est impair. Remarquons que l'on a l'\'egalit\'e $mult_{\lambda}(i)=j_{max}(i)$. Si $j_{max}(i)\in J^+$, $j_{max}(i)$ est impair. Inversement, supposons $j_{max}(i)$ impair. Si $mult_{\lambda}(i)=1$, $j_{max}(i)$ appartient \`a l'ensemble ${\cal J}^+\cup {\cal J}^-$ de 3.1 par d\'efinition des intervalles relatifs. L'imparit\'e impose $j_{max}(i)\in {\cal J}^+$. Or ${\cal J}^+\subset J^+$ par d\'efinition, donc $j_{max}(i)\in J^+$. Supposons $mult_{\lambda}(i)\geq2$. Par d\'efinition des intervalles relatifs, il existe $d=1,2$ et $\Delta_{d}\in \tilde{Int}(\lambda_{d})$ de sorte que $J(i)= \{j_{min}(i),...,j_{max}(i)\}\subset J(\Delta_{d})$. Pour fixer la notation, supposons que $d=1$, donc $ J(i)\subset J(\Delta_{1})$. Par d\'efinition des intervalles relatifs, $j_{max}(i) $ appartient \`a l'ensemble ${\cal J}$ de 3.1. L'imparit\'e  impose alors qu'il existe $d=1,2$ et $\Delta'_{d}\in \tilde{Int}(\lambda_{d})$ de sorte que $j_{max}(i)=j_{min}(\Delta'_{d})$. Si $d=1$, on a $j_{max}(i)\in J(\Delta_{1})\cap J(\Delta'_{1})$ donc $\Delta'_{1}=\Delta_{1}$. Mais $j_{min}(\Delta_{1})\leq j_{min}(i)<j_{max}(i)$, ce qui est contradictoire. Donc $d=2$. Alors $j_{max}(i)=j_{min}(\Delta'_{2})$. Puisque $j_{max}(i)\in J(\Delta_{1})$, $\lambda_{1,j}$ est pair. Alors, par d\'efinition de $J^+$, on a $j_{max}(i)\in J^+$. Cela prouve que les fonctions $\delta^{\zeta}$ v\'erifient la premi\`ere condition de la relation 3.5(1).
 
 Soit $i\in Jord^{bp}(\lambda)$. D'apr\`es la d\'efinition ci-dessus, $\tau^-(i)=\tau^+(i)+1$
si et seulement si $mult_{\lambda^{+,min}}(i)>0$, $mult_{\lambda^{-,min}}(i)>0$ et $\epsilon^{+,min}_{i}\not=\epsilon^{-,min}_{i}$. D'apr\`es la d\'efinition de $\chi$, ces conditions sont \'equivalentes \`a $mult_{\lambda}(i)\geq2$ et $\chi(i)=0$. La premi\`ere condition \'equivaut \`a $\vert J(i)\vert \geq2$. Sous cette condition, puisque $\chi=\chi_{\lambda_{1},\lambda_{2}}$, la seconde condition \'equivaut \`a $J(i)\subset J(\Delta_{2}(i))$ avec la notation de 3.5(1). Cela  ach\`eve de prouver cette condition 3.5(1). 

La condition 3.5(2) est claire.

Notons $i_{1}>...>i_{t}$ les entiers pairs $i\geq2$
tels que $mult_{\lambda^{+,min}}(i)$ soit impair. Pour $i\in Jord^{bp}(\lambda)$, on a $(-1)^{\delta^{\eta}(i)}-(-1)^{\delta^{\eta}(i^+)}\not=0$ si et seulement si $\delta^{\eta}(i)\not=\delta^{\eta}(i^+)$. Par d\'efinition de $\delta^{\eta}$, cela \'equivaut \`a  ce que  $mult_{\lambda^{+,min}}(i)$ soit impair, autrement dit \`a ce que $i=i_{h}$ pour un $h=1,...,t$. Pour un tel $i_{h}$, on a 
$$(-1)^{\delta^{\eta}(i_{h})}-(-1)^{\delta^{\eta}(i_{h}^+)}=2(-1)^{\delta^{\eta}(i_{h})}=2(-1)^{mult_{\lambda^{+,min}}(\geq i)}=2(-1)^h.$$
On a aussi
$$1-(-1)^{\tau^{\eta}(i_{h})}=1-\epsilon^{+,min}_{i_{h}}=\left\lbrace\begin{array}{cc}0,&\text{\,\,si\,\,}\epsilon^{+,min}_{i_{h}}=1,\\ 2,&\text{\,\,si\,\,}\epsilon^{+,min}_{i_{h}}=-1.\\ \end{array}\right.$$
On en d\'eduit
$$C^{\eta}=4\vert \{h=1,...,t; h\text{\,\,pair\,\,et\,\,}\epsilon^{+,min}_{i_{h}}=-1\}\vert -4\vert \{h=1,...,t; h\text{\,\,impair\,\,et\,\,}\epsilon^{+,min}_{i_{h}}=-1\}\vert.$$
En utilisant  1.3(1), on obtient
$$(3) \qquad C^{\eta}=\left\lbrace\begin{array}{cc}2k^+,&\text{\,\,si\,\,}k^+\text{\,\,est\,\,pair,}\\-2(k^++1),&\text{\,\,si\,\,}k^+\text{\,\,est\,\,impair.}\\ \end{array}\right.$$
On a une formule analogue pour $C^{-\eta}$, o\`u $k^+$ est remplac\'e par $k^-$. En reprenant les d\'efinitions de $r'$ et $r''$ donn\'ee en 1.3, un calcul cas par cas montre que (3) \'equivaut \`a
$$C^{\eta}=\left\lbrace\begin{array}{cc}2(r'+r''),&\text{\,\,si\,\,}r'+r''\text{\,\,est\,\,pair,}\\-2(r'+r''+1),&\text{\,\,si\,\,}r'+r''\text{\,\,est\,\,impair.}\\ \end{array}\right.$$
De m\^eme, l'\'egalit\'e analogue de (3) pour $C^{-\eta}$ \'equivaut \`a
$$C^{-\eta}=\left\lbrace\begin{array}{cc}2(r'-r''),&\text{\,\,si\,\,}r'+r''\text{\,\,est\,\,pair,}\\-2(r'-r''+1),&\text{\,\,si\,\,}r'+r''\text{\,\,est\,\,impair.}\\ \end{array}\right.$$
Par d\'efinition, $r'=r_{1}$ et $r''=\eta r_{2}$. Alors les formules ci-dessus sont la condition 3.5(3). Cela prouve (2).
 
 On peut appliquer le lemme 3.5. On note $\underline{\iota}_{1}$ et $\underline{\iota}_{2}$ les termes dont ce lemme affirme l'existence.  Avec les notations de 4.2, ils appartiennent \`a  ${\cal F}am_{r_{1}}(\lambda_{1})$, resp. ${\cal F}am_{r_{2}}(\lambda_{2})$. En cons\'equence, ces ensembles sont non vides. A fortiori, on a
 
 (4) $r_{1}^2+r_{1}\leq n_{1}$, $r_{2}^2\leq n_{2}$.

 Appliquons maintenant le calcul de 4.2 aux couples $(\mu_{1},1)\in \boldsymbol{{\cal P}^{orth}}(2n_{1}+1)_{k=1}$ et $(\mu_{2},1)\in \boldsymbol{{\cal P}^{orth}}(2n_{2})_{k=0}$. On a \'evidemment $sp(\mu_{1},1)=\mu_{1}$ et $sp(\mu_{2},1)=\mu_{2}$.  La condition (3) de ce paragraphe est v\'erifi\'ee: c'est (4) ci-dessus. Dans la formule 4.2(6), on peut limiter les sommations aux quadruplets $(\lambda^{'+},\epsilon^{'+},\lambda^{'-},\epsilon^{'-})$ et aux couples $(\iota_{1},\iota_{2})$ tels que $x(\lambda^{'+},\epsilon^{'+},\lambda^{'-},\epsilon^{'-})\not=0$ et
   $$m(\Pi_{\underline{\gamma}} (\lambda^{'+},\epsilon^{'+},\lambda^{'-},\epsilon^{'-}),\rho_{\iota_{1}}\otimes \rho_{\iota_{2}})\not=0.$$   
    Comme en 4.4, on d\'eduit de ces conditions les relations (3) et (4) de ce paragraphe:
 
  $\lambda^{+,min}\leq \lambda^{'+}$, $\lambda^{-,min}\leq \lambda^{'-}$, $\lambda^{'+}\cup \lambda^{'-}\leq \lambda$. 
  
  Mais ici $\lambda=\lambda^{+,min}\cup \lambda^{-,min}$ par d\'efinition. Les in\'egalit\'es ci-dessus sont donc des \'egalit\'es. D'apr\`es 4.1 et 4.2(5), les conditions  $\lambda^{+,min}= \lambda^{'+}$, $\lambda^{-,min}= \lambda^{'-}$ et $x(\lambda^{'+},\epsilon^{'+},\lambda^{'-},\epsilon^{'-})\not=0$ impliquent $(\lambda^{'+},\epsilon^{'+})=(\lambda^{+,min},\epsilon^{+,min})$ et $(\lambda^{'-},\epsilon^{'-})=(\lambda^{-,min},\epsilon^{-,min})$. Dans la somme 4.2(6), il ne reste que le quadruplet $(\lambda^{+,min},\epsilon^{+,min},\lambda^{-,min},\epsilon^{-,min})$ et on sait d'apr\`es 4.1 que, pour celui-l\`a, on a $x(\lambda^{+,min},\epsilon^{+,min},\lambda^{-,min},\epsilon^{-,min})=1$. 
  
  Il ne reste aussi que les couples $(\iota_{1},\iota_{2})$  tels que $m(\Pi_{\underline{\gamma}} (\lambda^{+,min},\epsilon^{+,min},\lambda^{-,min},\epsilon^{-,min}),\rho_{\iota_{1}}\otimes \rho_{\iota_{2}})\not=0$. 
   Ou encore, d'apr\`es 4.3(1),  tels que
   $m(\Pi^{\zeta}(\iota_{1},\iota_{2}),\rho_{\lambda^{+,min},\epsilon^{+,min}}\otimes \rho_{\lambda^{-,min},\epsilon^{-,min}})\not=0$, le signe $\zeta$ \'etant d\'etermin\'e comme en 4.3. 
   Cette condition \'equivaut \`a ce que 
   
   \noindent $(\lambda^{+,min},\epsilon^{+,min},\lambda^{-,min},\epsilon^{-,min})$ appartienne \`a l'ensemble ${\cal I}^{\zeta}(\iota_{1},\iota_{2})$ d\'efini en 3.4. 
   Puisque $\lambda^{+,min}\cup \lambda^{-,min}=\lambda$, la proposition 3.4(ii) nous dit qu'elle \'equivaut aussi \`a ce que $(\lambda^{+,min},\epsilon^{+,min},\lambda^{-,min},\epsilon^{-,min})$ appartienne \`a ${\cal I}^{\zeta,max}(\iota_{1},\iota_{2})$. En outre, on a dans ce cas 
   $$m(\Pi_{\underline{\gamma}} (\lambda^{+,min},\epsilon^{+,min},\lambda^{-,min},\epsilon^{-,min}),\rho_{\iota_{1}}\otimes \rho_{\iota_{2}})=1.$$
    La condition $(\lambda^{+,min},\epsilon^{+,min},\lambda^{-,min},\epsilon^{-,min}) \in {\cal I}^{\zeta,max}(\iota_{1},\iota_{2})$ \'equivaut \`a ce que les formules (4) et (5) de 3.4 soient v\'erifi\'ees, avec les modifications suivantes: les couples $(\lambda^+,\epsilon^+)$ et $(\lambda^-,\epsilon^-)$ de ce paragraphe sont remplac\'es par $(\lambda^{+,min},\epsilon^{+,min})$ et $(\lambda^{-,min},\epsilon^{-,min})$; les fonctions $\delta^+$, $\delta^-$, $\tau^+$ et $\tau^-$ sont remplac\'ees par $\delta_{\iota_{1},\iota_{2}}^+$ etc... La condition (4) d\'etermine enti\`erement les fonctions $\delta^{+}_{\iota_{1},\iota_{2}}$ et $\delta^{-}_{\iota_{1},\iota_{2}}$. En se rappelant que le signe $\zeta\nu$ qui intervient vaut pr\'ecis\'ement $\eta$ (cf. 4.3(2)), on voit que ces fonctions co\"{\i}ncident avec les fonctions $\delta^+$ et $\delta^-$ construites ci-dessus. Les fonctions $\tau^+_{\iota_{1},\iota_{2}}$ et $\tau^-_{\iota_{1},\iota_{2}}$ ne sont pas \`a premi\`ere vue enti\`erement d\'etermin\'ees par la relation (5) de 3.4. Toutefois, pour tout $i\in Jord^{bp}(\lambda)$, l'une au moins des valeurs $\tau^+_{\iota_{1},\iota_{2}}(i)$ ou $\tau^-_{\iota_{1},\iota_{2}}(i)$ est d\'etermin\'ee et co\"{\i}ncide avec la valeur de $\tau^+(i)$ ou $\tau^-(i)$. Puisque les couples $(\tau^+,\tau^-)$ et $\tau_{\iota_{1},\iota_{2}}^+,\tau_{\iota_{1},\iota_{2}}^-)$ v\'erifient tous deux la condition 3.5(1), cela suffit \`a conclure que ces deux couples sont \'egaux. Alors le lemme 3.5 nous dit que $(\iota_{1},\iota_{2})$ est \'egal au couple $(\underline{\iota}_{1},\underline{\iota}_{2})$ introduit ci-dessus. Inversement, pour ce dernier couple, les conditions (4) et (5) de 3.4 sont bien v\'erifi\'ees. Autrement dit, dans la somme 4.2(6), il ne reste plus que le couple  $(\underline{\iota}_{1},\underline{\iota}_{2})$ et on a  $m(\Pi_{\underline{\gamma}} (\lambda^{+,min},\epsilon^{+,min},\lambda^{-,min},\epsilon^{-,min}),\rho_{\underline{\iota}_{1}}\otimes \rho_{\underline{\iota}_{2}})=1$. 
   
 Cette formule 4.2(6) devient
 $$(5) \qquad M_{\pi}(\mu_{1},1;\mu_{2},1)=\vert Fam(\lambda_{1})\vert ^{-1/2}\vert Fam(\lambda_{2})\vert ^{-1/2}(-1)^{<\Lambda_{1},\Lambda_{\underline{\iota}_{1}}>}$$
 $$\left( (-1)^{<\Lambda_{2}^+,\Lambda_{\underline{\iota}_{2}}>}+sgn_{\sharp}(-1)^{<\Lambda_{2}^-,\Lambda_{\underline{\iota}_{2}}>}\right).  $$
 
  Rappelons que $\Lambda_{2}^+$ et $\Lambda_{2}^-$ sont les symboles des couples $(0,\rho_{\mu_{2},1}^+\otimes sgn)$ et $(0,\rho_{\mu_{2},1}^-\otimes sgn)$. Ils se d\'eduisent l'un de l'autre par permutation des deux termes $X$ et $Y$ de chaque symbole. D'apr\`es 2.5(1), on a donc
  $$(6) \qquad (-1)^{<\Lambda_{2}^-,\Lambda_{\underline{\iota}_{2}}>}=(-1)^{r_{2}}(-1)^{<\Lambda_{2}^+,\Lambda_{\underline{\iota}_{2}}>}.$$
  Consid\'erons la formule 1.5(1). Notons $i_{1}>...>i_{t}$ les entiers pairs $i\geq2$ tels que $mult_{\lambda^+}(i)$ soit impair. Le premier produit de la formule vaut $(-1)^{X^+}$, o\`u
  $$X^+=\vert \{h=1,...,t; \epsilon^+_{i_{h}}=-1\}.$$
  On a
  $$X^+\equiv \vert \{h=1,...,t; h\text{\,\,pair\,\,et\,\,}\epsilon^{+}_{i_{h}}=-1\}\vert -\vert \{h=1,...,t; h\text{\,\,impair\,\,et\,\,}\epsilon^{+}_{i_{h}}=-1\}\vert\,\, mod\,\,2{\mathbb Z}.$$
  D'apr\`es  1.3(1),  le membre de droite vaut $k^+/2$ si $k^+$ est pair, $-(k^++1)/2$ si $k^+$ est impair. D'apr\`es le m\^eme calcul cas par cas qui a calcul\'e $C^{\eta}$ ci-dessus, c'est aussi $(r'+r'')/2$ si $r'+r''$ est pair, $-(r'+r''+1)/2$ si $r'+r''$ est impair.
  On obtient
  $$(-1)^{X^+}=\left\lbrace\begin{array}{cc}(-1)^{(r'+r'')/2},&\text{\,\,si\,\,}r'+r''
\text{\,\,est\,\,pair,}\\  (-1)^{(r'+r''+1)/2},&\text{\,\,si\,\,}r'+r''
\text{\,\,est\,\,impair.}\\ \end{array}\right.$$
  Le deuxi\`eme facteur de 1.5(1) se calcule de m\^eme, $r''$ \'etant remplac\'e par $-r''$. Le produit de ces termes vaut $(-1)^{r''}$, ou encore $(-1)^{r_{2}}$. La formule 1.5(1) nous dit donc que
  $$(7) \qquad sgn_{\sharp}=(-1)^{r_{2}}.$$ 
  Gr\^ace \`a (6) et (7), (5) se simplifie en
  $$M_{\pi}(\mu_{1},1;\mu_{2},1)=2\vert Fam(\lambda_{1})\vert ^{-1/2}\vert Fam(\lambda_{2})\vert ^{-1/2}(-1)^{<\Lambda_{1},\Lambda_{\underline{\iota}_{1}}>+<\Lambda_{2}^+,\Lambda_{\underline{\iota}_{2}}>}.$$
  Donc $M_{\pi}(\mu_{1},1;\mu_{2},1)\not=0$.  Alors, en vertu de (1), les couples $(\mu_{1},1)$ et $(\mu_{2},1)$ v\'erifient le (ii) de la proposition 1.4, \`a ceci pr\`es que l'on doit de plus prouver que $n_{2}\geq1$ si $\sharp=an$. Mais, si  $\sharp=an$, (7) implique que $r_{2}$ est impair et (4) implique alors que $n_{2}\geq1$.
  
  \bigskip
  
  \subsection{Conclusion}
  
  On a prouv\'e que $\mu$ v\'erifiait les conditions de la proposition 1.4. Celle-ci implique que $\mu$ est le front d'onde de $\pi(\lambda^+,\epsilon^+,\lambda^-,\epsilon^-)$. Cela d\'emontre le deuxi\`eme th\'eor\`eme de l'introduction. Comme on l'a dit dans celle-ci, le premier th\'eor\`eme s'en d\'eduit gr\^ace \`a \cite{W4} 3.4.
  
  \bigskip
  
  \section{Sur le calcul effectif du front d'onde}
  
  \bigskip
  
  \subsection{Le couple $(\lambda^{max},\epsilon^{max})$}
  Soit $(\lambda,\epsilon)\in \boldsymbol{{\cal P}^{symp}}(2n)$, supposons que tous les termes de $\lambda$ sont pairs. On lui associe un couple $(\lambda^{max},\epsilon^{max})\in \boldsymbol{{\cal P}^{symp}}(2n)$ par r\'ecurrence sur $n$, selon la construction qui suit et qui est extraite de \cite{W5} 6.1 et 6.2. On repr\'esente $\lambda$ sous la forme $\lambda=(\lambda_{1},...,\lambda_{2r+1})$, avec $\lambda_{2r+1}=0$. On associe \`a $\epsilon$ une fonction encore not\'ee $\epsilon$ sur l'ensemble d'indices $\{1,...,2r+1\}$ par $\epsilon(j)=\epsilon_{\lambda_{j}}$ pour $j\in \{1,...,2r+1\}$, avec la convention $\epsilon_{0}=1$. On note $\mathfrak{S}$ la r\'eunion de $\{1\}$ et de l'ensemble des $j\in \{2,...,2r+1\}$ tels que $\epsilon(j)(-1)^j\not=\epsilon(j-1)(-1)^{j-1}$. On note $s_{1}=1<s_{2}<...<s_{S}$ les \'el\'ements de $\mathfrak{S}$. Pour $\zeta\in \{\pm 1\}$, notons $J^{\zeta}=\{j=1,...,2r+1; (-1)^{j+1}\epsilon(j)=\zeta\}$. On pose
  $$\lambda_{1}^{max}=(\sum_{h=1,...,S}\lambda_{s_{h}}) +2[S/2]-2\vert J^{-\epsilon(1)}\vert .$$
On pose $n'=n-\lambda_{1}^{max}/2$.   On note $\lambda'$ la r\'eunion des $\lambda_{j}$ pour $j\in J^{\epsilon(1)}-(J^{\epsilon(1)}\cap \mathfrak{S})$ et des $\lambda_{j}+2$ pour $j\in J^{-\epsilon(1)}-(J^{-\epsilon(1)}\cap \mathfrak{S})$. Pour $i\in Jord^{bp}(\lambda')$, on a $i=\lambda_{j}$ ou $i=\lambda_{j}+2$ pour un $j$ comme ci-dessus. On note $h[j]$ le plus grand entier $h\in \{1,...,S\}$ tel que $s_{h}<j$ et on pose $\epsilon'_{i}=(-1)^{h[j]+1}\epsilon(j)$ ($j$ n'est pas uniquement d\'etermin\'e par $i$ mais on montre que cette d\'efinition ne d\'epend pas du choix de $j$).  On montre que  $n'<n$  (si $n\not=0$), que le couple $(\lambda',\epsilon')$ appartient \`a $\boldsymbol{{\cal P}^{symp}}(2n')$ et que tous les termes de $\lambda'$ sont pairs. Par r\'ecurrence, on dispose d'un couple $(\lambda^{'max},\epsilon^{'max})$. On pose $\lambda^{max}=\{\lambda_{1}^{max}\}\cup \lambda^{'max}$. On d\'efinit $\epsilon^{max}$ par $\epsilon^{max}_{\lambda_{1}^{max}}=\epsilon_{\lambda_{1}}$ et $\epsilon^{max}_{i}=\epsilon^{'max}_{i}$ pour $i\in Jord^{bp}(\lambda^{'max})$ (c'est possible, c'est-\`a-dire que, si $\lambda_{1}^{max}$ appartient \`a $Jord^{bp}(\lambda^{'max})$, on a l'\'egalit\'e $\epsilon_{\lambda_{1}}=\epsilon^{'max}_{\lambda_{1}^{max}}$). Cela d\'efinit le couple $(\lambda^{max},\epsilon^{max})$. Les termes de $\lambda^{max}$ sont pairs et $\lambda_{1}^{max}$ est bien le plus grand terme de $\lambda^{max}$.

\bigskip

\subsection{La partition $^t\lambda^{min}$}
On conserve les m\^emes hypoth\`eses. On pose $k=k_{\lambda,\epsilon}$. On a $k_{\lambda^{max},\epsilon^{max}}=k$. On a rappel\'e en 1.3(1) comment se calculait l'entier $k$.  
On \'ecrit $\lambda^{max}=(\lambda^{max}_{1},...,\lambda^{max}_{2R+1})$ avec $\lambda^{max}_{2R+1}=0$. On note $j_{1}^+<...<j_{N}^+$ les $j=1,...,2R+1$ tels que $\epsilon^{max}( j)(-1)^{j+1}=(-1)^k$ (en consid\'erant comme dans le paragraphe pr\'ec\'edent que $\epsilon^{max}$ se d\'efinit sur l'ensemble d'indices). On note $j_{1}^-<...<j_{N}^-$ les $j=1,...,2R+1$ tels que $\epsilon^{max}( j)(-1)^{j}=(-1)^k$. 
 On v\'erifie que $N=R+[k/2]+1$, $M=R-[k/2]$. Notons $\nu'$ la r\'eunion disjointe des partitions suivantes:

$\{2R+3u-k-1+\lambda^{max}_{j^+_{u}}-2j^+_{u}; u=1,...,N\}$;

$\{2R+3v+k+\lambda^{max}_{j^-_{v}}-2j^-_{v}; v=1,...,M\}$;

$\{R+[(k-1)/2],R+[(k-1)/2]-1,...,0\}$;

$\{R-[(k+3)/2],R-[(k+3)/2]-1,...,0\}$.

On note $\nu'=(\nu'_{1},...,\nu'_{4R+1})$. 
Pour $j=1,...,4R+1$, on pose $\nu_{j}=\nu'_{j}-2R+[j/2]$. Cela d\'efinit une partition $\nu$ et on a l'\'egalit\'e $^t\lambda^{min}=\nu$ (cette \'egalit\'e se d\'eduit de  \cite{W5} 5.6 et 5.7).

\bigskip

\subsection{Exemples}
Soit $(\lambda^+,\epsilon^+,\lambda^-,\epsilon^-)\in \mathfrak{Irr}^{bp}_{quad}(2n)$. Les formules des deux paragraphes pr\'ec\'edents permettent de calculer les transpos\'ees des  partitions $\lambda^{+,min}$ et $\lambda^{-,min}$. Le front d'onde de $\pi(\lambda^+,\epsilon^+,\lambda^-,\epsilon^-)$ est $d(\lambda^{+,min}\cup \lambda^{-,min})$. Cette partition duale se calcule ainsi: on 
 note $\nu$ la partition obtenue en ajoutant $1$ au plus grand terme de ${^t\lambda}^{+,min}+{^t\lambda}^{-,min}$; alors $d(\lambda^{+,min}\cup \lambda^{-,min})$
est la plus grande partition orthogonale $\mu$ de $2n+1$ telle que $\mu\leq \nu$. Le moins que l'on puisse dire est que ce calcul n'est pas simple. 

Signalons le cas particulier rassurant o\`u $\epsilon^+=1$, c'est-\`a-dire $\epsilon^+_{i}=1$ pour tout $i\in Jord^{bp}(\lambda^+)$, et $\epsilon^- =1$. Dans ce cas, on voit que $\lambda^{+,max}=(2n^+)$, $\epsilon^{+,max}_{2n^+}=1$, $\lambda^{-,max}=(2n^-)$ et $\epsilon^{-,max}_{2n^-}=1$. On a $k^+=k^-=0$. On calcule $^t\lambda^{+,min}=(2n^+)$, $^t\lambda^{-,min}=(2n^-)$, puis  $d(\lambda^{+,min}\cup \lambda^{-,min})=(2n+1)$. Autrement dit, notre repr\'esentation $\pi(\lambda^+,1,\lambda^-,1)$ admet un mod\`ele de Whittaker usuel, ce qui est bien connu.

Un autre cas particulier est celui o\`u, pour $\zeta=\pm$, $n^{\pm}$ est de la forme $h^{\pm}(h^{\pm}+1)$, $\lambda^{\zeta}$ est \'egal \`a $(2h^{\zeta},2h^{\zeta}-2,...,2)$ et o\`u $\epsilon^{\zeta}$ est altern\'e, c'est-\`a-dire $\epsilon^{\zeta}_{2i}=(-1)^{i}$ pour $i=1,...,h^{\zeta}$. Dans ce cas, on v\'erifie que $\lambda^{\zeta,max}=\lambda^{\zeta,min}=\lambda^{\zeta}$. Le front d'onde de   $\pi(\lambda^+,\epsilon^+,\lambda^-,\epsilon^-)$ est alors $d(\lambda^+\cup \lambda^-)$. On retrouve le r\'esultat de \cite{M} et \cite{W4} car notre repr\'esentation est ici cuspidale donc  \'egale  \`a son image par l'involution d'Aubert-Zelevinsky. 

Donnons enfin comme exemple le calcul du front d'onde de $\pi(\lambda^+,\epsilon^+,\lambda^-,\epsilon^-)$ dans le cas o\`u $\lambda^-$ est vide et o\`u $\lambda^+$ a au plus trois termes non nuls. On pose simplement $\lambda=\lambda^+$, $\epsilon=\epsilon^+$, $\mu=d(\lambda^{min})$. On identifie $\epsilon$ au triplet $(\epsilon(1),\epsilon(2),\epsilon(3))$ que l'on note comme un triplet de signes $\pm$. Evidemment,  certains triplets ne sont autoris\'es que sous certaines hypoth\`eses sur $\lambda$: si $\epsilon(j)\not=\epsilon(j+1)$, on doit avoir $\lambda_{j}>\lambda_{j+1}$; si $\epsilon(j)=-$, on doit avoir $\lambda_{j}>0$. On note de m\^eme $\epsilon^{max}$ comme une famille de signes. Les r\'esultats sont les suivants:

$$\begin{array}{ccccc}\epsilon& k& \lambda^{max}&\epsilon^{max}& {^t\lambda}^{min}\\ (+,+,+)&0&(\lambda_{1}+\lambda_{2}+\lambda_{3})&(+)&(\lambda_{1}+\lambda_{2}+\lambda_{3})\\ (+,+,-)&1&(\lambda_{1}+\lambda_{2}-4,\lambda_{3}+2,2)&(+,+,-)&(\lambda_{1}+\lambda_{2}-2,\lambda_{3},1,1)\\ (+,-,+)& 2& (\lambda_{1},\lambda_{2},\lambda_{3})&(+,-,+)&(\lambda_{1}-2,\lambda_{2},\lambda_{3}+1,1)   \\ (+,-,-)&0&(\lambda_{1}+\lambda_{3}-2,\lambda_{2},2)&(+,-,-)& (\lambda_{1}+\lambda_{3}-2,\lambda_{2}+2) \\ (-,+,+)&1&(\lambda_{1}+\lambda_{3},\lambda_{2})&(-,+)&(\lambda_{1}+\lambda_{3}-1,\lambda_{2}+1) \\ (-,+,-)&3&(\lambda_{1},\lambda_{2},\lambda_{3})&(-,+,-)&(\lambda_{1}-3,\lambda_{2}-1,\lambda_{3},2,1,1)  \\ (-,-,+)& 0&(\lambda_{1}+\lambda_{2}-2,\lambda_{3}+2) & (-,-)&(\lambda_{1}+\lambda_{2}-1,\lambda_{3}+1)\\ (-,-,-)& 1&(\lambda_{1}+\lambda_{2}+\lambda_{3})&(-)&(\lambda_{1}+\lambda_{2}+\lambda_{3}-1,1)\\ \end{array}$$
\bigskip

$$\begin{array}{cc}\epsilon&\mu\\ (+,+,+)&(\lambda_{1}+\lambda_{2}+\lambda_{3}+1)\\ (+,+,-)&(\lambda_{1}+\lambda_{2}-1,\lambda_{3}-1,1,1,1)\\ (+,-,+)&(\lambda_{1}-1,\lambda_{2}-1,\lambda_{3}+1,1,1)\\ (+,-,-)&(\lambda_{1}+\lambda_{3}-1,\lambda_{2}+1,1)\\ (-,+,+)&(\lambda_{1}+\lambda_{3}-1,\lambda_{2}+1,1)\\ (-,+,-)&(\lambda_{1}-3,\lambda_{2}-1,\lambda_{3}+1,1,1,1,1)\\ (-,-,+)&(\lambda_{1}+\lambda_{2}-1,\lambda_{3}+1,1)\\ (-,-,-)&(\lambda_{1}+\lambda_{2}+\lambda_{3}-1,1,1)\\ \end{array}$$

CNRS Institut de Math\'ematiques de Jussieu Paris Rive Gauche

4 place Jussieu

75005 Paris

jean-loup.waldspurger@imj-prg.fr

 \end{document}